\documentclass[12pt]{article}
\usepackage{color}
\usepackage{amsmath}
\usepackage{amssymb}
\usepackage{amsthm}
\usepackage[totalwidth=17cm, totalheight=25cm]{geometry}
\usepackage{graphicx}
\usepackage{mathabx}
\usepackage{tabularx}
	\newcolumntype{C}[1]{>{\centering\arraybackslash}m{#1}} 
	\newcolumntype{R}[1]{>{\raggedleft\arraybackslash}m{#1}} 

\newtheoremstyle{boldplain}
{9pt}
{9pt}
{\itshape}
{}
{\bfseries}
{.}
{.5em}
{\thmname{#1}\thmnumber{ #2}\thmnote{ (#3)}}%

\newtheoremstyle{bolddefinition}
{9pt}
{9pt}
{}
{}
{\bfseries}
{.}
{.5em}
{\thmname{#1}\thmnumber{ #2}\thmnote{ (#3)}}%

\theoremstyle{boldplain}
\newtheorem{add}[equation]{Addendum}

\newtheorem{cor}[equation]{Corollary}

\newtheorem{lem}[equation]{Lemma}

\newtheorem{prop}[equation]{Proposition}

\newtheorem{thm}[equation]{Theorem}

\theoremstyle{bolddefinition}
\newtheorem{dfn}[equation]{Definition}

\newtheorem{example}[equation]{Example}

\newtheorem{ques}[equation]{Question}

\newtheorem{rem}[equation]{Remark}

\newfont{\bigbf}{cmbx10 scaled\magstep1}
\normalbaselines
\numberwithin{equation}{section}

\def\no{\noindent}

\def\C{{\mathbb C}}
\def\R{{\mathbb R}}
\def\F{{\mathbb F}}

\def\H{{\mathbb H}}

\def\N{{\mathbb N}}
\def\P{{\mathbb P}}
\def\Z{{\mathbb Z}}

\def\al{\alpha}
\def\ga{\gamma}
\def\Ga{\Gamma}
\def\de{\delta}
\def\De{\Delta}
\def\eps{\epsilon}
\def\la{\lambda}
\def\La{\Lambda}
\def\si{\sigma}
\def\Si{\Sigma}
\def\ups{\upsilon}

\def\Om{\Omega}

\def\3{\ss}

\def\acc{\;\operatorname{acc}\;}
\def\Acc{\operatorname{Acc}}
\def\acts{\curvearrowright}
\def\amod{a_{mod}}
\def\Aut{\mathop{\hbox{Aut}}}

\def\D{\partial}
\def\DF{\partial_{F\ddot u}}

\def\Flag{\operatorname{Flag}}
\def\Flagn{\Flag_{\nu_{mod}}}
\def\Flags{\Flag_{\si_{mod}}}
\def\Flagt{\Flag_{\tau_{mod}}}
\def\Flagit{\Flag_{\iota\tau_{mod}}}
\def\Flagmt{\Flag_{-\tau_{mod}}}
\def\Flagpmt{\Flag_{\pm\tau_{mod}}}

\def\Flagu{\Flag_{\ups_{mod}}}

\def\geo{\partial_{\infty}}

\def\Homeo{\operatorname{Homeo}}
\def\id{\mathop{\hbox{id}}}

\def\interior{\operatorname{int}}
\def\inte{\operatorname{int}}
\def\Isom{\mathop{\hbox{Isom}}}

\def\Lan{\La_{\nu_{mod}}}
\def\Lap{\La_{\pi_{mod}}}
\def\Las{\La_{\si_{mod}}}

\def\Lat{\La_{\tau_{mod}}}
\def\Latm{\La_{\tau_{mod}}^-}
\def\Latp{\La_{\tau_{mod}}^+}
\def\Latpm{\La_{\tau_{mod}}^{\pm}}
\def\Lamt{\La_{-\tau_{mod}}}

\def\lra{\longrightarrow}

\def\numod{\nu_{mod}}
\def\oa{\overrightarrow}
\def\ol{\overline}

\def\pihalf{\frac{\pi}{2}}
\def\pithird{\frac{\pi}{3}}
\def\2pithird{\frac{2\pi}{3}}
\def\piforth{\frac{\pi}{4}}
\def\pisixth{\frac{\pi}{6}}
\def\pos{\operatorname{pos}}
\def\cpos{\operatorname{c-pos}}
\def\ccpos{\operatorname{c-c-pos}}
\def\prect{\prec_{\taumod}}
\def\proj{\operatorname{proj}}

\def\rank{\mathop{\hbox{rank}}}

\def\Ra{\Rightarrow}
\def\simod{\si_{mod}}
\def\epsmod{\eps_{mod}}
\def\pimod{\pi_{mod}}
\def\taumod{\tau_{mod}}
\def\htaumod{\hat\tau_{mod}}
\def\tiltaumod{\tilde\tau_{mod}}
\def\pmtaumod{\pm\tau_{mod}}
\def\upsmod{\ups_{mod}}

\def\st{\operatorname{st}}
\def\ost{\mathop{\hbox{ost}}}

\def\tangle{\angle_{Tits}}

\def\tits{\partial_{Tits}}
\def\Th{\operatorname{Th}}

\def\ThF{\operatorname{Th}_{F\ddot u}}

\def\Wn{W_{\nu_{mod}}}
\def\Wp{W_{\pi_{mod}}}
\def\Wt{W_{\tau_{mod}}}

\def\8{\infty}
\def\<{\langle}
\def\>{\rangle}

\def\BI{\begin{itemize}}
\def\EI{\end{itemize}}

\title{Dynamics on flag manifolds:\\ domains of proper discontinuity and cocompactness}
\author{Michael Kapovich, Bernhard Leeb, Joan Porti}
\date{March 3, 2017}

\begin{document}

\maketitle
\hfill{\em To Guiomar}

\begin{abstract}
For noncompact semisimple Lie groups $G$
with finite center 
we study the dynamics of the actions of their discrete subgroups $\Ga<G$ 
on the associated partial flag manifolds $G/P$.
Our study is based on the observation
that they exhibit also in higher rank 
a certain form of convergence type dynamics.
We identify geometrically  domains of proper discontinuity in all partial flag manifolds.
Under certain dynamical assumptions 
equivalent to the Anosov subgroup condition,
we establish the cocompactness of the $\Ga$-action on various domains of proper discontinuity,
in particular on domains in the full flag manifold $G/B$. 
In the regular case (e.g.\ of $B$-Anosov subgroups)
we prove the nonemptiness of such domains 
if $G$ has (locally) at least one noncompact simple factor 
not of the type $A_1, B_2$ or $G_2$,
by showing the nonexistence of certain ball packings of the visual boundary.
\end{abstract}

\tableofcontents

\section{Introduction}

Let $G$ be a noncompact semisimple Lie group with finite center.
In this paper, 
we study the natural actions 
$$ \Ga \acts G/P$$
of discrete subgroups $\Ga<G$
on the (partial) flag manifolds $G/P$ associated to $G$.
(Here, $P<G$ denotes a parabolic subgroup.) 
We are interested in aspects of the  
topological dynamics of the action of $\Ga$,
notably in domains of proper discontinuity and criteria for the cocompactness 
on these
(see Theorems~\ref{thm:pdmainintro}, \ref{thm:cocomainintro} and~\ref{thm:centintro} below).
Our approach relies on the geometry of the associated symmetric space $X=G/K$ 
of noncompact type.
The connection is established by the fact that
the flag manifolds occur as the $G$-orbits in the visual boundary $\geo X$,
that is, the boundary at infinity of the visual compactification $\ol X=X\sqcup\geo X$ of $X$.
The results are essentially generalizations of the main results 
in the first version of this paper \cite{coco13}, 
however we adopt here a more dynamical viewpoint.

If $\rank(X)=1$, equivalently, if $X$ has strictly negative sectional curvature,
then the only flag manifold is $\geo X$ itself 
and the transitive action $G\acts\geo X$ has {\em convergence dynamics}.
This means that divergent sequences in $G$ exhibit a certain attraction-repulsion behavior,
namely they subconverge on the complement of one point in $\geo X$ locally uniformly to a constant map.
More precisely,
for a sequence $g_n\to\infty$ in $G$ there exist a subsequence $(g_{n_k})$
and (not necessarily distinct) points $\xi_{\pm}\in\geo X$
such that 
$$ g_{n_k}|_{\geo X-\{\xi_-\}} \to \xi_+$$
uniformly on compacta.

As a consequence of convergence dynamics, 
for a discrete subgroup $\Ga<G$,
there is a clean $\Ga$-invariant {\em dynamical decomposition} 
$$\geo X=\Om_{disc}\sqcup\La$$
into the open {\em domain of discontinuity} or {\em wandering set} $\Om_{disc}$ and the compact {\em limit set} $\La$.
The latter consists of all points occuring as limits $\xi_+$ as above
for sequences $\ga_n\to\infty$ in $\Ga$,
and the $\Ga$-action on $\Om_{disc}$ is even {\em properly} discontinuous. 
In order for this action to be cocompact,
one needs to impose further conditions on the group.
The action $\Ga\acts\Om_{disc}$ is {\em cocompact} if (but not only if) $\Ga$ is {\em convex-cocompact}.

If $\rank(X)\geq2$,
then the action $G\acts\geo X$ is no longer transitive.
The $G$-orbits are compact and, as $G$-spaces, copies of flag manifolds. 
They are parametrized by the spherical Weyl chamber $\simod$ associated to $G$;
a $G$-orbit for an interior point of 
a face $\taumod\subseteq\simod$
(which we will also refer to as a {\em face type})
is naturally identified with the flag manifold $\Flagt\cong G/P_{\tau}$,
the conjugacy class of the parabolic subgroup $P_{\tau}$ corresponding to the face $\taumod$.
In particular, 
the {\em regular} $G$-orbits in $\geo X$,
i.e.\ those corresponding to interior points of $\simod$ itself,
are identified with the full flag manifold,
$\DF X\cong\Flags\cong G/B_{\si}$,
the space of Weyl chambers at infinity,
also called the Furstenberg boundary;
here, $B_{\si}$ denotes a minimal parabolic subgroup of $G$.

Our study is based on the observation 
that a {\em weak} form of {\em convergence dynamics} persists 
for the action $G\acts\geo X$
in higher rank (cf.\ sections~\ref{sec:accac} and~\ref{sec:accfl}):
Sequences $g_n\to\infty$ in $G$ 
sub{\em accumulate} outside a compact exceptional subset locally uniformly at another compact subset,
meaning that 
there exist a subsequence $(g_{n_k})$
and compact subsets $A_{\pm}\subset\geo X$
such that 
$$ g_{n_k}|_{\geo X-A_-} \hbox{ accumulates at } A_+$$
uniformly on compacta. 
We briefly say in this case that $(g_{n_k})$ is {\em $(A_-,A_+)$-accumulating}.
There is a certain flexibility in the choice of the pair of compact subsets $(A_-,A_+)$ 
and a trade-off (``uncertainty relation"):
If one shrinks one of the subsets $A_{\pm}$, one must enlarge the other.

For instance, one can make the following {\em metric} choice for the pair of compact subsets:
For a suitable subsequence $(g_{n_k})$
there exist points $\xi_{\pm}\in\geo X$ such that 
$(g_{n_k})$ is $(\ol B(\xi_-,\pi-r),\ol B(\xi_+,r))$-accumulating
for all radii $r\in (0,\pi)$, 
where $B(\xi,r)$ denotes a ball in $\geo X$ with respect to the Tits angle metric.
This kind of convergence type dynamical behavior had been observed,
in the general setting of proper CAT(0) spaces, 
by Karlsson \cite[Thm.\ 1]{Karlsson} and Papasoglu and Swenson \cite[Thm.\ 4]{PapaSwen}, 
see also the first version of this paper \cite[Thm.\ 1.1 and \S 6.1]{coco13}.

In our setting of CAT(0) {\em model} spaces with their rich geometric structure,
one can make more flexible ``combinatorial" choices for the pair of compact subsets
which can be described in terms of the (partial) Bruhat order $\prec$ on the Weyl group $W$.
These will enable us to construct {\em larger} domains of proper discontinuity for discrete group actions
than those obtained from metric choices.
To explain the combinatorial choices, we need some preliminary considerations.

It will be useful for us to interpret the Bruhat order geometrically
and we give a 
geometric description of it and its generalization as the {\em folding order}
(see section~\ref{sec:foldord}).

Since any two Weyl chambers at infinity 
$\si,\si'\subset\geo X$
are contained in an apartment, that is, the visual boundary $\geo F\subset\geo X$ of a maximal flat $F\subset X$,
we can define a combinatorial {\em relative position} 
$$ \pos(\si',\si)\in W$$
of $\si'$ with respect to $\si$
(section~\ref{sec:relpos}). 
The larger the position is with respect to $\prec$, the more generic it is.
We note that 
the sublevels of $\pos(\cdot,\si)$ in $\DF X$ are precisely the Schubert cycles with respect to $\si$,
that is, the $B_{\si}$-orbit closures.

We define a {\em thickening} 
$$ \Th\subset W $$
of the neutral element inside the Weyl group 
as a union of sublevels for the Bruhat order
(section~\ref{sec:chthick}).
Each thickening $\Th\subset W$
gives rise to corresponding thickenings of chambers inside the Furstenberg boundary,
$$ \ThF(\si) :=\{ \pos(\cdot,\si)\in\Th \} \subset\DF X ,$$
which can also be regarded as thickenings 
$\Th(\si)\subset\geo X$ inside the visual boundary, by taking the union of the chambers contained in them.
The thickenings of chambers in $\DF X$ are finite unions of Schubert cycles and hence projective subvarieties.
Thickenings $\ThF(A)\subset\DF X$ and $\Th(A)\subset\geo X$
of compact subsets $A\subset\DF X$ 
are defined as the union of the thickenings of the individual chambers $\si\in A$;
they are again compact.

This discussion generalizes:
There is a well-defined $\Wt\backslash W$-valued position $\pos(\cdot,\tau)$ 
of chambers relative to a {\em simplex} $\tau\in\Flagt$.
Here, $\Wt<W$ denotes the stabilizer of the face
$\taumod\subseteq\simod$.
If the thickening $\Th\subset W$ is $\Wt$-left invariant,
then it yields well-defined compact thickenings 
$\ThF(A)\subset\DF X$ and $\Th(A)\subset\geo X$
of compact subsets $A\subset\Flagt$.
Even more generally,
there is a well-defined $\Wt\backslash W/\Wn$-valued position $\pos(\nu,\tau)$ 
of {\em simplices} $\nu\in\Flagn$ relative to simplices $\tau\in\Flagt$,
and a $\Wt$-left and $\Wn$-right invariant thickening inside $W$
yields thickenings of subsets of $\Flagt$ inside $\Flagn$.

For every thickening $\Th\subset W$, there is the {\em complementary} thickening $\Th^c\subset W$
defined by 
$$ W=\Th\sqcup w_0\Th^c $$
where $w_0$ denotes the longest element of the Weyl group.
(This partition of the Weyl group generalizes the decomposition $\pi=r+(\pi-r)$
of the maximal distance in the unit sphere, cf.\ the metric choices above.)
We call a thickening $\Th$ {\em slim}, {\em fat} or {\em balanced}
if $\Th\subseteq\Th^c$, 
$\Th\supseteq\Th^c$, 
respectively, $\Th=\Th^c$.
Existence results for balanced thickenings with different invariance properties
are stated in Proposition~\ref{prop:exbalintro} below,
examples are given in section~\ref{sec:chthick}.

\medskip
Returning to the dynamics of sequences $g_n\to\infty$ in $G$ on $\geo X$,
we show (cf.\ Lemma~\ref{lem:obspureg} and Proposition~\ref{prop:regimplcontr})
that there always exist a subsequence $(g_{n_k})$,
a face $\taumod\subseteq\simod$ and a pair of simplices $\tau_{\pm}\in\Flagpmt$ 
such that
\begin{equation}
\label{eq:atrrepg}
g_{n_k}|_{C(\tau_-)} \to \tau_+
\end{equation}
uniformly on compacta 
in the open Schubert cell $C(\tau_-)\subset\Flagt$ of simplices opposite to $\tau_-$.
Here $-\taumod:=\iota\taumod$ for the canonical involution $\iota=-w_0$ of $\simod$. 
This locally uniform convergence property implies 
(Corollary~\ref{cor:orbacc})
that, more generally, 
\begin{equation}
\label{eq:accdivg}
\hbox{$(g_{n_k})$ is $(\Th^c(\tau_-),\Th(\tau_+))$-accumulating on $\geo X$}
\end{equation}
for {\em all} $\Wt$-left invariant thickenings $\Th\subset W$.
Note that (\ref{eq:atrrepg}) is equivalent to (\ref{eq:accdivg})
for the slimmest nonempty choice $\Th=\{e\}$. 
We call sequences satisfying (\ref{eq:atrrepg}) {\em $\taumod$-contracting}
(cf.\ Definition~\ref{def:contracting_sequence}).
An equivalent notion had been introduced by Benoist in \cite{Benoist},
see in particular part (5) of his Lemma 3.5,
cf.\ also Remark~\ref{rem:convregdyn}.

\medskip
We now turn to discussing the {\em dynamics of discrete subgroups} $\Ga<G$ on $\geo X$.

For a face $\taumod\subseteq\simod$
we define the ``small" {\em $\taumod$-limit set}
$$ \Lat\subset\Flagt$$
as the set of all simplices in $\Flagt$
which occur as limits $\tau_+$ as in (\ref{eq:atrrepg})
for sequences $\ga_n\to\infty$ in $\Ga$
(cf.\ Definition~\ref{dfn:lims}).
\begin{rem}
Benoist introduced in \cite[\S 3.6]{Benoist} 
a notion of limit set $\La_{\Ga}$ for Zariski dense subgroups $\Ga$ of reductive algebraic groups over local fields
which in the case of real semisimple Lie groups 
is equivalent to our concept of $\simod$-limit set $\Las$.\footnote{Benoist's limit set $\La_{\Ga}$ is contained in the flag manifold $Y_{\Ga}$
which in the case of real Lie groups is the full flag manifold $G/B$, see the beginning of \S 3 of his paper.
It consists of the limit points of sequences contracting on $G/B$, cf.\ his Definitions 3.5 and 3.6.}
What we call the $\taumod$-limit set $\Lat$ for other face types $\taumod\subsetneq\simod$ 
is mentioned in his Remark 3.6(3),
and his work implies that, in the Zariski dense case, 
$\Lat$ is the image of $\Las$ under the natural projection $\Flags\to\Flagt$ of flag manifolds.
\end{rem}
By choosing a $\Wt$-left invariant thickening $\Th\subset W$,
we obtain from these small limit sets the ``large" {\em thickened} limit sets $\Th(\Lat)\subset\geo X$.
Our main results concern the proper discontinuity and cocompactness 
of the $\Ga$-action on the complements 
$$ \Om_{\Th} := \geo X - \Th(\Lat) ,$$
respectively, 
on their intersections with the $G$-orbits $G\eta\subset\geo X$.
We obtain the strongest results for the dynamics on the Furstenberg boundary.
This is reasonable because the latter fibers with compact fiber over all partial flag manifolds,
and cocompact domains of proper discontinuity in any flag manifold 
pull back to such domains in $\DF X$.
Our results are of the kind, in the spirit of Mumford's {\em Geometric Invariant Theory},
that the $\Ga$-actions become properly discontinuous 
when removing a sufficiently ``fat" thickening of $\Lat$,
and remain cocompact when removing a sufficiently ``slim" one.
(See Example~\ref{ex:config} and the discussion in section~\ref{sec:git} 
for a concrete connection with configuration spaces and GIT.)

The accumulation property (\ref{eq:accdivg}) is the key step in constructing {\em domains of proper discontinuity} 
for all discrete subgroups 
(see Propositions~\ref{prop:pddisc} and~\ref{prop:pddisclar} in section~\ref{sec:discr}).
We obtain the most useful results for subgroups
which satisfy a certain generalization of the convergence property
(see Definition~\ref{dfn:convact}):
\begin{dfn}[Weak convergence subgroup]
We call a discrete subgroup $\Ga<G$ a {\em $\taumod$-convergence subgroup}
with respect to a face type $\taumod\subseteq\simod$,
if every sequence $\ga_n\to\infty$ in $\Ga$ has a subsequence
satisfying (\ref{eq:atrrepg}) with this particular face type $\taumod$,
equivalently, 
has a subsequence
satisfying (\ref{eq:accdivg}) for any choice of $\Wt$-left invariant
thickenings $\Th\subset W$.
\end{dfn}

\begin{rem}[Convergence dynamics versus regularity]
\label{rem:convregdyn}
We note that the $\taumod$-con\-ver\-gen\-ce property of a subgroup $\Ga< G$,
formulated in terms of the dynamics of the action on the visual boundary $\geo X$,
can be equivalently described in terms of the asymptotic behavior of $\Ga$-orbits 
in the symmetric space $X$. 
Namely, $\Ga$ is a {\em $\taumod$-convergence subgroup} if and only if it is a {\em $\taumod$-regular  subgroup} of $G$, 
see section~\ref{sec:wreg}.
The notion of $\taumod$-regularity was introduced in our earlier paper \cite{morse}
where also the equivalence of the two notions was established.
In the present paper 
we only need (and verify)  that $\taumod$-convergence implies $\taumod$-regularity.

The notions of regularity and contraction for sequences
and their essential equivalence 
can be found in some form already in the work of Benoist,
see \cite[\S 3]{Benoist}. 
For the sake of completeness we give independent proofs in our setting 
of discrete subgroups of semisimple Lie groups.
Also our methods are rather different.
We give a geometric treatment
and present the material
in a form suitable for the further development
of our theory of discrete isometry groups acting on Riemannian symmetric spaces and euclidean buildings,
such as in our papers \cite{mlem,bordif}.
\end{rem}

Our main result on proper discontinuity includes (cf.\ Theorem~\ref{thm:pdwconv}):
\begin{thm}[Proper discontinuity outside fat thickenings]
\label{thm:pdmainintro}
Let $\taumod\subseteq\simod$ be $\iota$-invariant.
If $\Ga<G$ is a $\taumod$-con\-ver\-gen\-ce subgroup,
then for any fat $\Wt$-left invariant thickening $\Th\subset W$
the action 
$$ \Ga\acts\geo X-\Th(\Lat) $$
is properly discontinuous. 
\end{thm}

In order to obtain {\em cocompactness} for actions of $\taumod$-convergence subgroups,
we must impose further conditions,
as it is the case for convergence actions,
compare the situation in rank one. 
Our main requirement is that the action 
$\Ga\acts\Flagt$ should be {\em expanding at $\Lat$} in the sense of Sullivan \cite[\S 9]{Sullivan},
cf.\ Definition~\ref{dfn:expact}. 
Moreover, 
if $\taumod$ is $\iota$-invariant,
we call 
the limit set 
$\Lat$ {\em antipodal}
if the simplices in it are pairwise opposite
(see Definition~\ref{dfn:antip}(ii)).

\begin{dfn}[CEA subgroup]
\label{dfn:ceaintro}
For a $\iota$-invariant face $\taumod\subseteq\simod$ 
we call a $\taumod$-con\-ver\-gen\-ce subgroup $\Ga<G$ a {\em $\taumod$-CEA subgroup}
(convergence, expanding, antipodal)
if $\Lat$ is antipodal and if the action $\Ga\acts\Flagt$ is expanding at $\Lat$.
\end{dfn}
The restricted action $\Ga\acts\Lat$ is then a convergence action in the traditional sense.
Such subgroups are higher-rank generalizations of convex-cocompact subgroups of rank one Lie groups.
In fact, the CEA condition is only one of various equivalent dynamical and (coarse) geometric conditions 
which can be used to characterize this class of discrete subgroups, 
see \cite{morse} 
and also \cite{manicures,anosov,anolec} 
for a detailed study of these conditions and their equivalence.
In particular:
\begin{rem}[CEA versus Anosov]
The class of $\taumod$-CEA subgroups coincides with the class of $P_{\taumod}$-Anosov subgroups,
see \cite[\S 6.5]{morse}.
Here, $P_{\taumod}$ refers to the conjugacy class of parabolic subgroups of $G$ 
corresponding to the face $\taumod$ of the spherical Weyl chamber $\simod$.
We recall that 
the notion of Anosov subgroup had first been introduced in \cite{Labourie} 
using the language of geodesic flows, and further extended in \cite{GW}.
We gave the first flow-free definitions in \cite{morse}.
We note that Labourie's original definition did not require $\taumod$ to be $\iota$-invariant,
instead he worked with $(P_{\taumod},P_{\iota\taumod})$-Anosov subgroups.
However, as already observed in \cite{GW}, 
the general case readily reduces to the $\iota$-invariant one. 
\end{rem}

Our main result regarding cocompactness includes 
(cf.\ Theorem~\ref{thm:cocomain} and Corollary~\ref{cor:cocomain}):

\begin{thm}[Cocompactness outside slim thickenings]
\label{thm:cocomainintro}
Suppose that $\Ga<G$ is a $\taumod$-CEA subgroup. 
Then for each slim $\Wt$-left invariant thickening $\Th\subset W$,
the action 
$$\Ga\acts \DF X-\ThF(\Lat)$$ 
is cocompact. 

More generally, 
suppose that $\numod\subseteq\simod$ is another face type 
and that the thickening $\Th$ is also $\Wn$-right invariant.
Then for any $G$-orbit $G\eta\subset\geo X$ 
corresponding to an interior point of $\numod$,
the action 
$$\Ga\acts G\eta-\Th(\Lat)$$ 
is cocompact. 
\end{thm}

By combining the two theorems, we obtain the central result of this paper:
\begin{thm}[Cocompact domains of proper discontinuity]
\label{thm:centintro}
Let $\Ga<G$ and the data $\taumod,\numod,\Th$ be as in the previous theorem
with the additional requirement that the thickening $\Th$ be balanced.
Then the respective actions are properly discontinuous and cocompact. 
\end{thm}

We note that topology of the quotient space $(\DF X-\ThF(\Lat))/\Ga$, in general, depends on 
the balanced thickening,  see Example \ref{ex:product_actions}.  

Balanced thickenings do not exist for all invariance requirements,
but for many they do.
For instance, one can impose arbitrary left invariance and,
as a consequence, one has balanced thickenings of $\taumod$-limits sets inside $\DF X$
for all $\iota$-invariant $\taumod$,
as the first part of the next result shows
(see section~\ref{sec:chthick} for more general results):

\begin{prop}[Existence of balanced thickenings]
\label{prop:exbalintro}
For every $\iota$-invariant face type $\taumod$ 
there exists a $\Wt$-left invariant balanced thickening $\Th\subset W$.

For an arbitrary face type $\numod$,
a $\Wn$-right invariant balanced thickening 
exists if and only if 
left multiplication by $w_0$ has no fixed point on $W/\Wn$.
This is the case, for instance, if $w_0=-\id$,
equivalently, 
if all irreducible factors of the symmetric space are of type 
$A_1$, $B_{n\geq2}$, $D_{2k\geq4}$, $E_{7,8}$, $F_4$ or $G_2$. 
\end{prop}

The {\em nonemptiness} of the domains found in Theorem~\ref{thm:centintro} 
(and Theorem~\ref{thm:pdmainintro})
is an issue. For instance, uniform lattices in rank one Lie groups have empty domains of discontinuity at infinity (and such lattices are CEA). See also Example~\ref{ex:product_actions} for empty domains in the reducible case.
If for a $\taumod$-convergence subgroup $\Ga$ 
with antipodal $\taumod$-limit set
all domains given by Theorem~\ref{thm:pdmainintro} were empty,
it would follow that the visual boundary of $X$ admits a {\em packing}
by a compact family (with respect to the visual topology) 
of $\pihalf$-balls (with respect to the Tits metric),
cf.\ Proposition~\ref{prop:neccdd}.
However, 
the existence of such packings can be ruled out
for most Weyl groups (Theorem~\ref{thm:nopack}),
and we conclude
(see Theorem~\ref{thm:nedoman}):
\begin{thm}[Nonemptiness of domains of proper discontinuity]
Suppose that $X$ has at least one de Rham factor not of the type $A_1, B_2$ or $G_2$,
and let $\Ga<G$ be a $\simod$-convergence subgroup
with antipodal limit set $\Las$.
Then for some balanced thickening $\Th\subset W$
the domain of proper discontinuity $\DF X-\ThF(\Las)$ for the $\Ga$-action 
provided by Theorem~\ref{thm:pdmainintro}
is nonempty.
\end{thm}

Note that the theorem covers the case of CEA subgroups,
but is more general.

The possible balanced thickenings can be described more precisely, cf.\ Theorem~\ref{thm:nedoman}.
In the $B_2$-case, we have partial nonemptiness results
for the groups $G=O(2k+1,2)$ with $k\geq1$
(see Addendum~\ref{add:nedomanb2orth}).
The $G_2$-case is not discussed in this paper.

The above results yield for the dynamics of CEA subgroups on the Furstenberg boundary:
\begin{cor}[Dynamics on the Furstenberg boundary]
Suppose that $\Ga< G$ is a $\taumod$-CEA subgroup. 
There exist $\Wt$-left invariant balanced thickenings $\Th\subset W$,
and for every such thickening the action $$\Ga\acts \DF X-\ThF(\Lat)$$
is properly discontinuous and cocompact.

If $\Ga< G$ is a $\simod$-CEA subgroup,  
and if $X$ has at least one de Rham factor not of the type $A_1, B_2$ or $G_2$,
then for some balanced thickening $\Th\subset W$ 
the cocompact domain of proper discontinuity $\DF X-\ThF(\Las)$
is nonempty.
\end{cor}

Again, the possible thickenings occuring in the $\simod$-case can be described more precisely. 

\begin{rem}[Dynamics on Finsler compactifications]
Our results regarding domains of proper discontinuity and cocompactness 
for discrete group actions on flag manifolds have analogs 
for the actions of the same classes of subgroups on a {\em Finsler compactification} 
$\ol X^{Fins}$ of $X$.
This is done in our paper \cite{bordif}.
The compactification $\ol X^{Fins}$ is obtained from $X$ geometrically by applying the horoboundary construction 
to suitable $G$-invariant {\em regular polyhedral Finsler} metrics on $X$ 
rather than to $G$-invariant Riemannian metrics (which yields the visual compactification
$\ol X=X\sqcup\geo X$),
and it coincides with the maximal Satake compactification from algebraic group theory,
see also \cite{Parreau}.
Note that the Furstenberg boundary $\DF X$ naturally embeds into $\ol X^{Fins}$
as a $G$-orbit, namely as the only compact one. 
Some of the results become easier in the Finsler setting, for instance, 
the nonemptiness of domains of proper discontinuity at infinity is no longer an issue:
Each $\taumod$-convergence subgroup with antipodal limit set $\Lat$ 
has a nonempty domain of proper discontinuity in the Finsler ideal boundary 
(defined using an arbitrary $W_{\taumod}$-left invariant balanced thickening), once $\rank(X)\geq2$, 
see \cite[Lemma 9.19]{bordif}.  
\end{rem}

\begin{rem}
There is overlap of our results with \cite{GW}.
There, cocompact domains of proper discontinuity are constructed 
for Anosov subgroups of various semisimple Lie groups 
acting on various partial flag manifolds.
However,
in the general case of arbitrary semisimple Lie groups $G$,
such domains are constructed only in 
$G$-homogeneous spaces fibering over $\DF X\cong G/B$
with compact fiber \cite[Thm.\ 1.9]{GW}.
Nonemptiness of these domains is proven 
for $P$-Anosov subgroups of small cohomological dimension 
\cite[Thms.\ 1.11, 1.12 and 9.10]{GW},
while our nonemptiness results apply to $\simod$-convergence subgroups with antipodal $\simod$-limit set,
which includes $B$-Anosov subgroups, without restriction on the cohomological dimension.

Observe also that our treatment is intrinsic, 
while in \cite{GW} first a theory for Anosov subgroups of Lie groups of the type $\Aut(F)$ is developed 
(where the $F$'s are certain bilinear and hermitian forms), 
and then generalized to other semisimple Lie groups by embedding these into the groups $\Aut(F)$. 
The intrinsic approach is more uniform and seems to provide better control,
e.g.\ it allows us to get the domains, 
for general semisimple Lie groups,
in flag manifolds instead of only in bundles over these as in \cite{GW}.
While in some low rank cases the outcomes of the two constructions of thickened limit sets are the same, 
our construction appears to be more general.
\end{rem}

The earlier version \cite{coco13} of this preprint 
written in 2013 covered only the $\simod$-regular case. 
Some of the material of \cite{coco13}, 
dealing with equivalent characterizations of $\taumod$-CEA actions, 
was later moved to our paper \cite{morse}.
Most of the rest of the material of \cite{coco13} 
was generalized and moved into this paper.

\medskip 
{\bf Acknowledgements.} 
The first author was supported by the NSF grants DMS-09-05802 and DMS-12-05312. 
The last author was supported by the grants Mineco MTM2012-34834, AGAUR  SGR2009-1207 
and the Icrea Acad\`emia Award 2008. The three authors are also grateful to 
the GEAR grant which partially supported the IHP trimester in Winter of 2012 (DMS 1107452,
1107263, 1107367 ``RNMS: Geometric structures and representation varieties'' 
(the GEAR Network)), 
and the Max Planck Institute for Mathematics in Bonn, 
where some of this work was done.

\section{Geometric preliminaries}
\label{sec:prelim}

In this section we collect some standard material on  
Coxeter complexes, the geometry of nonpositively curved symmetric spaces and associated 
spherical Tits buildings; we refer the reader to \cite{qirigid} and \cite{habil} 
for more detailed discussion of symmetric spaces and buildings.

\subsection{General notation}

We will use the notation
$B(a,r)$ and $\ol B(a,r)$ 
for the open, respectively, closed $r$-ball, centered at $a$ in a metric space $Z$.
We will denote the nearest point distance of a point $z\in Z$ 
to a subset $A\subset Z$ by 
$d(z, A):= \inf d(z,\cdot)|_A$. 
The Hausdorff distance between two subsets $A, B\subset Z$ will be denoted $d_H(A,B)$.
A {\em geodesic} in a metric space is an isometric embedding from a (possibly infinite) interval $I\subset\R$.

\subsection{Coxeter complexes}
\label{sec:cox}

A {\em spherical Coxeter complex} $\amod$ is  a pair $(S,W)$ 
consisting of a unit sphere $S$ in a Euclidean vector space $V$ and 
a finite group $W$ which acts isometrically on $S$ 
and is generated by reflections at hyperplanes. A Coxeter complex is {\em reducible} if $W$ splits as a (nontrivial) direct product $W_1\times W_2$ and $V$ admits a $W$-invariant (nontrivial) orthogonal direct sum decomposition $V=V_1\oplus V_2$ such that $W_i$ fixes $V_{3-i}, i=1,2$. In this case, we obtain two induced Coxeter complexes $(S_i, W_i)$ on the unit spheres $S_i\subset V_i$. A Coxeter complex which is not reducible is called {\em irreducible}.

We will use the notation $\angle$ for the angular metric on $S$. 
Throughout the paper, we assume that $W$ does not fix a point in $S$ 
and is associated with a root system $R$. 
Spherical Coxeter complexes will occur as {\em model apartments} 
of spherical buildings, 
mostly of Tits boundaries of symmetric spaces, 
and will in this context usually be denoted by $\amod$. 

A {\em wall} $m_{\rho}$ in $S$ 
is the fixed point set of a hyperplane reflection $\rho$ in $W$. 
A {\em half-apartment} in $S$ is 
a closed hemisphere  bounded by a wall. 
A point $\xi\in S$ is called {\em singular} 
if it belongs to a wall and {\em regular} otherwise. 

The action $W\acts S$ determines on $S$ a structure 
as a simplicial complex 
whose facets, called {\em chambers}, 
are the closures of the connected components of
\begin{equation*}
S-\bigcup_{\rho} m_{\rho}
\end{equation*}
where the union is taken over all reflections $\rho$ in $W$. 
We will refer to the simplices in this complex as {\em faces}. 
(If one allows fixed points for $W$ on $S$, 
then $S$ carries only a structure as a cell complex.) 
Codimension one faces of this complex are called {\em panels}. The {\em interior} $\interior(\tau)$ 
of a face $\tau$ is the complement in $\tau$ to the union of walls not containing $\tau$. The interiors  
$\interior(\tau)$ are called {\em open simplices}. 
A geodesic sphere in $S$ is called {\em singular}
if it is simplicial, equivalently, if it equals an intersection of walls. 

Each chamber is a fundamental domain for the action $W\acts S$. 
We define the {\em spherical model Weyl chamber} 
as the quotient $\simod=S/W$. 
The natural projection $\theta:S\to\simod$ 
restricts to an isometry on every chamber. 
An important elementary property of the chamber $\simod$ is 
that its diameter (with respect to the spherical metric) is $\le \pihalf$. 

For a face $\taumod$ of $\simod$, we define the subgroup $\Wt\subset W$
as the stabilizer of $\taumod$ in $W$. 
Accordingly, 
for a point $\bar\xi\in\simod$,
we define $W_{\bar\xi}\subset W$ as the stabilizer of $\bar\xi$ in $W$. 
Then $W_{\bar\xi}=\Wt$ where $\taumod$ is the face of $\simod$ spanned by $\bar\xi$,
i.e.\ which contains $\bar\xi$ as an interior point.
Note that $W_{\simod}=1$ and $W_{\bar\xi}=1$ for $\bar\xi\in\interior(\simod)$.

It is convenient, and we will frequently do so, to identify  $\simod$ with a chamber $\si\subset S$ 
(traditionally called the {\em positive chamber}). 
Such an identification 
determines a generating set of $W$, 
namely the reflections at the walls bounding $\simod$, 
and hence a word metric on $W$; 
the longest element with respect to this metric is denoted $w_0$. 
This element sends 
$\simod$ to the opposite chamber in $S$. 
We say that two points $\xi, \hat\xi\in S$ are {\em Weyl antipodes} if $\hat\xi=w_0 \xi$. 
We define the {\em standard} or {\em opposition} {\em involution}  
$$
\iota=\iota_S: S\to S$$ 
as the composition $-w_0$. This involution preserves $\simod$ 
and equals the identity if and only if $-\id_S\in W$ 
because then $w_0=-\id_S$. 

A point $\xi$ in $S$ is called a {\em root} 
if the hemisphere centered at $\xi$ is simplicial, 
equivalently, is bounded by a wall. If $(S,W)$ is associated with a root system $R$, then $\xi\in S$ is a root if and only if it has the direction of a  {\em coroot}. Note that irreducible root systems correspond to irreducible Coxeter complexes and vice versa.  
 
 \begin{rem}
 We will be assuming in what follows that $(S,W)$ is associated with a root system $R$ which spans $V^*$. Equivalently, $W$ is isomorphic to the linear part of an affine {\em crystallographic} Coxeter group, i.e., one  acting cocompactly on the affine space underlying the vector space $V$. The root system $R$ in this situation can be assumed to be {\em reduced}, i.e.,   if roots $\al, \beta$ have the same kernel then $\al=\pm \beta$.
In what follows we will be assuming that $R$ is reduced.  
\end{rem}

Note that each root type $\bar\zeta\in \simod$ is $\iota$-invariant, since the reflection $w\in W$ corresponding to the root $\bar\zeta$ sends  $\bar\zeta$ to $-\bar\zeta$. 
 
Each {\em irreducible} root system $R$ has one or two distinct root types, 
i.e.\ $W$ acts on $R$ with one or two orbits. Geometrically speaking, 
this means that $W$ acts on the set of walls with one or two orbits. 
We refer the reader to \cite{Bourbaki} for details.

\medskip 
Suppose that $S$ is identified with the 
sphere at infinity of a Euclidean space $F$, 
$S\cong\geo F$, 
where $\geo F$ is equipped with the angular metric. 
For a closed subset $A\subset S$ and a point $x\in F$ 
we define 
$V(x,A)\subset F$ 
as the complete cone over $A$ with tip $x$, 
that is, 
as the union of rays emanating from $x$ and asymptotic to $A$. 
If $\tau\subset S$ is a face, we call the cone $V(x, \tau)$ 
a {\em Weyl sector},
and if $\si\subset S$ is a chamber,
we call $V(x, \si)$ a {\em euclidean Weyl chamber}.

After fixing an origin $o\in F$,  
the group $W$ lifts to a group of isometries of $F$ fixing $o$. 
The euclidean Weyl chambers $V(o,\si)$ are then fundamental domains 
for the action of $W\acts F$. 

We define the {\em euclidean model Weyl chamber} 
as the quotient 
$V_{mod}=F/W$; 
we will also denote it by $\Delta$ or $\De_{euc}$. 
It is canonically isometric to 
the complete euclidean cone over $\simod$. 
The natural projection $$\proj:F\to V_{mod}=\De_{euc}=\De$$ 
restricts to an isometry on every euclidean Weyl chamber $V(o,\si)$. 

For a closed subset $\bar A\subset\simod$
we define $V(0,\bar A)\subset V_{mod}$
as the complete cone over $\bar A$ with tip $0$.
In particular,
a face $\taumod$ of $\simod$ corresponds
to a face $V(0,\taumod)$ of $V_{mod}$. 

We define the {\em $\Delta$-valued distance function} 
or {\em $\Delta$-distance} 
$d_{\De}$ on $F$ by:
\begin{equation*}
d_\Delta(x,y)= \proj(y-x)\in\De
\end{equation*}
Note the symmetry property:
\begin{equation}
\label{eq:symprop}
d_{\De}(x,y)= \iota d_{\De}(y,x) 
\end{equation}
The Weyl group is precisely the group of isometries  
for the $\Delta$-valued distance on $F_{mod}$ 
which fix the origin. 

\begin{lem}\label{lem:iota}
Suppose that the Coxeter complex $(S,W)$ is irreducible. Then $\iota=\id$ if and only if the root system of 
 $(S,W)$ is of type $A_1, B_\ell, C_\ell$, $D_{2k}, E_{7,8}$, $F_4$ or $G_2$. 
If $\iota\ne \id$, then $\simod$ contains exactly one root which, therefore, is $\iota$-invariant. 
\end{lem}
\proof The proof is by examination of the irreducible root systems, see e.g. \cite{Bourbaki}: 
$w_0=- \id$ if and only if the root system is of type $A_1, B_\ell, C_\ell$, $D_{2k}, E_{7,8}$, $F_4$ or $G_2$. 
All the remaining irreducible root systems are simply-laced; equivalently, $W$ acts transitively on roots. \qed 

\subsection{Hadamard manifolds}
\label{sec:had}

In this section only,
$X$ denotes a Hadamard manifold, 
i.e.\ a simply connected complete Riemannian manifold 
with nonpositive sectional curvature.  We will use the notation $\Isom(X)$ for the full isometry group of $X$. 

Any two points in $X$ are connected by a unique geodesic segment. 
We will use the notation $xy$ for the oriented geodesic 
segment connecting $x$ to $y$. 
We will often regard 
geodesic segments, geodesic rays and complete geodesics as parameterized with unit speed
and treat them as isometric maps of intervals to $X$.

We will denote by $\angle_x(y,z)$ 
the angle between the geodesic segments $xy$ and $xz$ at the point $x$. 
For $x\in X$ we let $\Si_xX$ denote 
the {\em space of directions} of $X$ at $x$, 
i.e.\ the unit sphere in the tangent space $T_xX$, 
equipped with the angle metric. 

The {\em ideal} or {\em visual boundary} of $X$, 
denoted $\geo X$, is the set of asymptote classes of geodesic rays in $X$, 
where two rays are {\em asymptotic} 
if and only if they have finite Hausdorff distance. 
Points in $\geo X$ are called {\em ideal points}. 
For $\xi\in \geo X$ and $x\in X$ we denote by $x\xi$ the geodesic ray 
emanating from $x$ and asymptotic to $\xi$, 
i.e.\ representing the ideal point $\xi$. For $x\in X$ we have a natural map
$$
\log_x: \geo X\to \Si_xX
$$
sending $\xi\in \geo X$ to the velocity vector at $x$ of the geodesic ray $x\xi$. 
The {\em cone} or {\em visual topology} on $\geo X$ 
is characterized by the property that all  the maps $\log_x$ are homeomorphisms;  
with respect to this topology, $\geo X$ is homeomorphic to 
the sphere of dimension $\dim(X)-1$. The visual topology extends to $\bar{X}=X\cup \geo X$ as follows: 
A sequence $(x_n)$ converges to an ideal point $\xi\in \geo X$ 
if the sequence of geodesic segments $xx_n$ 
emanating from some (any) base point $x$ 
converges to the ray $x\xi$ pointwise (equivalently, uniformly on compacta in $\R$). 
This topology makes $\bar{X}$ into a closed ball. 
We define the visual boundary of a subset $A\subset X$ 
as the set $\geo A=\bar A\cap\geo X$ 
of its accumulation points at infinity. 

The visual boundary $\geo X$ carries the natural 
{\em Tits (angle) metric} $\tangle$, 
defined as 
\begin{equation*}
\tangle(\xi,\eta)= \sup_{x\in X} \angle_{x}(\xi,\eta)
\end{equation*}
where $\angle_{x}(\xi,\eta)$ 
is the angle between the geodesic rays $x\xi$ and $x\eta$. 
The {\em Tits boundary} $\tits X$ is the metric space $(\geo X,\tangle)$. 
The Tits metric is lower semicontinuous with respect to the visual topology 
and, accordingly, 
the {\em Tits topology} induced by the Tits metric 
is finer than the visual topology. 
It is discrete if there is an upper negative curvature bound, 
and becomes nontrivial if flat directions occur. 
For instance, 
the Tits boundary of flat $r$-space is the unit $(r-1)$-sphere, 
$\tits\R^r\cong S^{r-1}(1)$. 
An isometric embedding $X\to Y$ of Hadamard spaces 
induces an isometric embedding $\tits X\to\tits Y$ of Tits boundaries. 

A subset $A$ of $\tits X$ is called {\em convex} if for any two points $\xi, \eta\in A$ with $\tangle(\xi,\eta)<\pi$, the (unique) geodesic 
$\xi\eta$ connecting $\xi$ and $\eta$ in $\tits X$ is entirely contained in $A$.

\subsection{Symmetric spaces of noncompact type}
\label{sec:symm}

The standard references for this and the following section are \cite{Eberlein} and \cite{Helgason}. 
Our treatment of this standard material is more geometric than the one presented in these books. 

A symmetric space, denoted by $X$ throughout this paper, 
is said to be of {\em noncompact type} 
if it is nonpositively curved, simply connected 
and has no Euclidean factor. 
In particular, it is a Hadamard manifold. 
We will identify $X$ with the quotient $G/K$ 
where $G$ is a semisimple Lie group with finite center 
acting isometrically and transitively on $X$, 
and $K$ is a maximal compact subgroup of $G$.  
We will assume that $G$ is commensurable with the isometry group $\Isom(X)$ 
in the sense that we allow compact kernel and finite cokernel 
for the natural map $G\to\Isom(X)$. 
In particular, 
the image of $G$ in $\Isom(X)$ contains the identity component $\Isom(X)_o$. 
The Lie group $G$ carries a natural structure as a real algebraic group. 

A {\em point reflection} (also known as a {\em Cartan involution}) at a point $x\in X$ 
is an isometry $\si_x$ which fixes $x$ and has differential $-\id_{T_xX}$ in $x$. 
In a symmetric space,
point reflections exist in all points (by definition).
A {\em transvection} of $X$ is an isometry 
which is the product $\si_x\si_{x'}$ of two point reflections; 
it preserves the oriented geodesic through $x$ and $x'$ 
and the parallel vector fields along it. 
The transvections preserving a unit speed geodesic $c(t)$ 
form a one parameter subgroup $(T^c_t)$ of $\Isom(X)_o$ 
where $T^c_t$ denotes the transvection mapping 
$c(s)\mapsto c(s+t)$. 
A nontrivial isometry $\phi$ of $X$ is called {\em axial} 
if it preserves a geodesic $l$ and shifts along it. 
(It does not have to be a transvection.)
The geodesic $l$ is called an {\em axis} of $\phi$. 
Axes are in general not unique. 
They are parallel to each other. 

A {\em flat} in $X$ is a totally geodesic flat submanifold, 
equivalently, 
a convex subset isometric to a Euclidean space. 
A maximal flat in $X$ is a flat 
which is not contained in any larger flat; 
we will use the notation $F$ for maximal flats. 
The group $\Isom(X)_o$ acts transitively on the set of maximal flats; 
the common dimension of maximal  flats is called the {\em rank} of $X$. 
The space $X$ has rank one if and only if it has 
strictly negative sectional curvature. 

A maximal flat $F$ 
is preserved by all transvections along geodesic lines contained in it. 
In general, there exist nontrivial isometries of $X$ fixing $F$ pointwise. 
The subgroup of isometries of $F$ 
which are induced by elements of $G$ 
is isomorphic to a semidirect product $ \R^r \rtimes W$, 
where $r$ is the rank of $X$. 
The subgroup $\R^r$ acts simply transitively on $F$ by translations. 
The linear part $W$ 
is a finite reflection group, called the {\em Weyl group} of $G$ and $X$. 
Since maximal flats are equivalent modulo $G$, 
the action $W\acts F$ is well-defined up to isometric conjugacy. 

We will think of the Weyl group as 
acting on a {\em model flat} $F_{mod}\cong \R^r$ 
and on its visual boundary sphere at infinity, 
the {\em model apartment} $\amod=\tits F_{mod}\cong S^{r-1}$. 
The pair $(\amod,W)$ is the {\em spherical Coxeter complex} 
associated with $X$. 
We identify the spherical model Weyl chamber $\simod$ 
with a (fundamental) chamber in the model apartment, 
$\simod\subset \amod$. 
Accordingly, 
we identify the {\em euclidean model Weyl chamber} $V_{mod}$ 
with the sector in $F_{mod}$ 
with tip in the origin and visual boundary $\simod$, 
$V_{mod}\subset F_{mod}$. 

The $\Delta$-valued distance naturally extends from $F_{mod}$ to $X$ 
because every pair of points lies in a maximal flat. 
In order to define the $\De$-distance $d_\Delta(x,y)$ of two points $x,y\in X$ 
one chooses a maximal flat $F$ containing $x,y$ 
and identifies it isometrically with $F_{mod}$ 
so that the Weyl group actions correspond.
The resulting quantity $d_\Delta(x,y)$ is independent of the choices. 
We refer the reader to \cite{KLM} for the detailed discussion of 
{\em metric properties} of $d_\Delta$.

For every maximal flat $F\subset X$, 
we have a Tits isometric embedding $\geo F\subset\geo X$ 
of its visual boundary sphere. 
There is an identification $\geo F\cong \amod$ 
with the model apartment, 
unique up to composition with elements in $W$. 
The Coxeter complex structure on $\amod$ 
induces a simplicial structure on $\geo F$. 
The visual boundaries of maximal flats cover $\geo X$ 
because every geodesic ray in $X$ is contained in a maximal flat.
Moreover, their intersections are simplicial. 
One thus obtains a $G$-invariant 
piecewise spherical {\em simplicial structure} on $\geo X$ 
which makes $\geo X$ into a {\em spherical building} and, 
also taking into account the visual topology, 
into a topological spherical building. 
It is called the {\em spherical} or {\em Tits building} 
associated to $X$. 
The Tits metric is the path metric 
with respect to the piecewise spherical structure. 
We will refer to the simplices as {\em faces}.

The visual boundaries $\geo F\subset\geo X$ 
of the maximal flats $F\subset X$ 
are precisely the {\em apartments} 
with respect to the spherical building structure at infinity, 
which in turn are precisely the convex subsets 
isometric to the unit $(r-1)$-sphere 
with respect to the Tits metric. 
Any two points in $\geo X$ lie in a common apartment. 

The action $G\acts\geo X$ on ideal points is not transitive if $X$ has rank $\ge 2$. 
Every $G$-orbit meets every chamber exactly once. 
The quotient can be identified with the spherical model chamber, 
$\geo X/G\cong\simod$. 
We call the projection 
\begin{equation*}
\theta:\geo X\to\geo X/G \cong\simod
\end{equation*}
the {\em type} map. 
It restricts to an isometry on every chamber $\si\subset\geo X$. 
We call the inverse 
$\kappa_{\si}=(\theta|_{\si})^{-1}: \simod\to\si$
the {\em (chamber) chart} for $\si$. 
Consequently, 
$\theta$ restricts to an isometry on every face $\tau\subset\geo X$. 
We call $\theta(\tau)\subset\simod$ 
the {\em type} of the face $\tau$ and 
$\kappa_{\tau}=(\theta|_{\tau})^{-1}:\theta(\tau)\to\tau$
its {\em chart}. 
We define the {\em type} of an ideal point $\xi\in\geo X$ 
as its image $\theta(\xi)\in\simod$. 
A point $\xi\in \geo X$ is called {\em regular} 
if its type is an interior point of $\simod$, 
and {\em singular} otherwise. 
We denote by $\geo^{reg}X\subset\geo X$ 
the set of regular ideal boundary points. 
A point $\rho\in \tits X$ is said to be of {\em root type} if $\theta(\rho)$ is a root in $\simod\subset S$. 
Equivalently, the closed $\pihalf$-ball centered at $\rho$ (with respect to the Tits metric) 
is simplicial, i.e.\ is a simplicial subcomplex of $\tits X$. 
If $a\subset\geo X$ is an apartment,
we call a type preserving isometry $\kappa_a:\amod\to a$ an {\em apartment chart} for $a$.

A geodesic segment $xy$ in $X$ is called {\em regular} if $x\ne y$ and for the unique geodesic ray $x\xi$ extending $xy$ 
the point $\xi\in \tits X$ is regular. 
Equivalently, the vector $d_\Delta(x,y)$ belongs  to the interior of $V_{mod}$. 

\begin{dfn}[Antipodal]
\label{dfn:antip}
(i) 
Two ideal points $\xi,\eta\in\geo X$ are 
{\em antipodal} 
if $\tangle(\xi, \eta)=\pi$. 
A subset of $\geo X$ is called {\em antipodal} if the points in it are pairwise antipodal.

(ii)
Two simplices $\tau_1,\tau_2\subset\geo X$ 
are {\em opposite} (or {\em antipodal}) {\em with respect to a point $x\in X$} if $\tau_2=\si_x\tau_1$, 
where $\si_x$ denotes the reflection at the point $x$. 
Two simplices $\tau_1,\tau_2\subset\geo X$ 
are {\em opposite} (or {\em antipodal}) 
if they are opposite simplices in the apartments containing both of them.
\end{dfn}
Note that the last property holds 
iff some (every) interior point of $\tau_1$ has an antipode in the interior of $\tau_2$,
equivalently, iff $\tau_1$ and $\tau_2$ are opposite with respect to some point $x\in X$.
Their types are then related by $\theta(\tau_2)=\iota(\theta(\tau_1))$.
We will frequently use the notation $\tau, \hat\tau$ and $\tau_+, \tau_-$ for pairs of antipodal simplices. 

\medskip
A pair of opposite chambers $\si_+,\si_-\subset\geo X$ 
is contained in a unique apartment,  
which we will denote by $a(\si_+,\si_-)$; the apartment  $a(\si_+,\si_-)$  
is the visual boundary of a unique maximal flat $F(\si_+, \si_-)$ in $X$.

For a point $x\in X$ and a simplex $\tau\subset\geo X$ 
we define the {\em (Weyl) sector} $V=V(x,\tau)\subset X$ 
as the union of rays $x\xi$ for all ideal points $\xi\in\tau$. 
Weyl sectors are contained in flats.
They are isometric images of faces $V(0, \taumod)\subset V_{mod}$ of the euclidean model Weyl chamber 
under isometric embeddings $F_{mod}\to X$ which are type preserving at infinity. 
More generally, 
for a point $x\in X$ and a closed subset $A\subset\geo X$,
we define the {\em Weyl cone} $V(x,A)$
as the union of all rays $x\xi$ for $\xi\in A$. 
Weyl cones are in general not flat.

\medskip
The stabilizers $B_{\si}\subset G$ of the chambers $\sigma\subset \geo X$ 
are the {\em minimal parabolic subgroups} of $G$. 
After choosing a reference chamber 
$\si_0\subset\geo X$, 
we call $B=B_{\si_0}$ the {\em positive} minimal parabolic subgroup. 
The group $G$ acts transitively on the set of chambers in $\geo X$, 
which we will then identify with $G/B$, 
the {\em  full flag manifold} of $G$. 
The minimal parabolic subgroups are algebraic subgroups of $G$, 
and $G/B$ is a real projective variety. 
The set $\DF X\cong G/B$ of chambers in $\geo X$ 
is called the {\em Furstenberg boundary} of $X$; 
we will equip it with the visual topology 
(as opposed to the Zariski topology coming from $G/B$) 
which coincides with its manifold topology 
as a compact homogeneous $G$-space. 
Every regular $G$-orbit $G\xi\subset\geo X$, $\xi \in\interior(\si_0)$, 
is $G$-equivariantly and homeomorphically identified with $\DF X$ 
by assigning to the (regular) ideal point $g\xi$ 
the unique chamber $g\si_0$ containing it. 

The stabilizers $P_{\tau}\subset G$ of simplices $\tau\subset\geo X$ 
are the {\em parabolic subgroups} of $G$. 
The group $G$ acts transitively on simplices of the same type. 
The set 
$\Flagt\cong G/P_{\taumod}$ 
of the simplices $\tau$ of type $\theta(\tau)=\taumod\subset\simod$ 
is called the {\em  partial flag manifold} of type $\taumod$. 
In particular, 
$\Flags=\DF X$. 
Again, we equip the flag manifolds with the visual topology; 
it agrees with their topology as compact homogeneous $G$-spaces. 
Every $G$-orbit $G\xi\subset\geo X$ of type $\theta(\xi)\in\interior(\taumod)$
is $G$-equivariantly homeomorphic to $\Flagt$.

For a flag manifold $\Flagt$ and a simplex $\hat\tau$ of type $\iota\taumod$ 
we define the {\em open Schubert stratum} $C(\hat\tau)\subset \Flagt$ as the subset of simplices opposite to $\hat\tau$
in the sense of Definition~\ref{dfn:antip}.
It follows from semicontinuity of the Tits distance that 
the subset $C(\hat\tau)\subset \Flagt$ is indeed open. Furthermore, this subset is also 
dense in $\Flagt$. We note that for rank 1 symmetric spaces, 
the only flag manifold associated to $G$ is $\geo X$ 
and the open Schubert strata are the complements of points. 

If $\taumod$ is $\iota$-invariant,
we say as in Definition~\ref{dfn:antip} 
that a subset of $\Flagt$ is {\em antipodal}
if the simplices in it are pairwise opposite.

\section{Geometry of visual boundaries}
\label{sec:geomsymmsp}

In this section 
we introduce definitions and prove some properties of symmetric spaces of noncompact type 
and their visual boundaries 
of more specific nature 
which are needed for our study of discrete group actions at infinity. 

\subsection{Stars at infinity and regular points}
\label{sec:star}

For a simplex $\tau\subset\geo X$,
the {\em star} $\st(\tau)\subset\geo X$ is the union of all closed chambers $\si\supseteq\tau$.
It is proven in \cite[Proposition 2.14]{morse} that for each face $\tau\subset\geo X$, 
the {\em Weyl cone} $V(x,\st(\tau))$ is a closed convex subset of $X$. 

\medskip
For a face type $\taumod\subseteq\simod$,
we define the {\em open star} $$\ost(\taumod)\subset\simod$$
as the union of all open faces of $\simod$ whose closure contains $\taumod$.
Its complement $$\D\st(\taumod):=\simod-\ost(\taumod)$$ 
is the union of all (closed) faces of $\simod$ which do not contain $\taumod$.

For a simplex $\tau\subset\geo X$,
we define the {\em open star} $$\ost(\tau)\subset\st(\tau)\subset\geo X$$
as the union of all open simplices in $\geo X$ whose closure contains $\tau$.
Then $$\D\st(\tau):=\st(\tau)-\ost(\tau)$$
is the union of all (closed) simplices in $\st(\tau)$ which do not contain $\tau$,
\begin{dfn}
An ideal point $\xi\in\geo X$ is said to be {\em $\taumod$-regular} 
if $\theta(\xi)\in\ost(\taumod)$, 
and {\em $\taumod$-singular} if $\theta(\xi)\in\D\st(\taumod)$.
\end{dfn}
We will call $\simod$-regular points simply {\em regular}.
Note that $\ost(\simod)=\inte(\simod)$,
and the regular points in $\geo X$ are precisely the interior points of chambers.

Note that $\ost(\tau)$ is the subset of $\theta(\tau)$-regular points in $\st(\tau)$
and $\D\st(\tau)$ is subset of $\theta(\tau)$-singular points.
The {\em $\taumod$-regular part} 
\begin{equation*}
\geo^{\taumod-reg} X=\theta^{-1}(\ost(\taumod))\subset\geo X 
\end{equation*}
of the visual boundary 
contains all open chambers and is in particular dense in $\geo X$
(also with respect to the Tits topology). 
For a $\taumod$-regular point $\xi\in\geo X$ 
there is a unique closest (with respect to the Tits metric) simplex $\tau\subset\geo X$ of type $\taumod$,
namely the one with $\xi\in\ost(\tau)$.

The notion of regularity extends to oriented geodesic segments, rays and lines in $X$:
A geodesic ray $x\eta\subset X$ is $\taumod$-regular if its ideal endpoint $\eta\in\geo X$ is.
An oriented geodesic segment $xy\subset X$ is $\taumod$-regular 
if the geodesic ray $x\eta$ extending it is.

The geometric significance of $\simod$-regularity of geodesic segments 
comes  from the fact that a geodesic segment (or ray) in $X$ 
is $\simod$-regular iff it is contained in a unique maximal flat.

\subsection{Folding order}
\label{sec:foldord}

In this section, we discuss natural partial orders on Weyl orbits in the model apartment 
and give different equivalent geometric definitions for them.

By a {\em folding map} $\amod\to \amod$ 
we mean a type preserving continuous map 
which sends chambers isometrically onto chambers. 

We will be working with folding maps which fix some reference face 
and think of them as moving points ``closer'' towards this face.
For instance, for a simplicial hemisphere $h\subset\amod$ (containing the reference face)
there is the folding map fixing $h$ and reflecting the complementary hemisphere onto it,
see the discussion of special foldings below.

\begin{dfn}[Folding order]
\label{dfn:fold}
For a face type $\taumod\subseteq\simod$, 
we define the {\em $\taumod$-folding order} $\prec_{\taumod}$ on $\amod$ as follows:
For distinct points $\bar\xi_1,\bar\xi_2\in \amod$ 
we say that $\bar\xi_1\prec_{\taumod}\bar\xi_2$ 
if and only if there exists a folding map 
$f: \amod\to \amod$ 
such that $f|_{\taumod}=\id_{\taumod}$ and $f(\bar\xi_2)=\bar\xi_1$. 

We will use the notation $\prec$ for $\prec_{\simod}$.
\end{dfn}
The relations $\prect$ are {\em transitive},
because the composition of folding maps are folding maps.
\begin{rem}
\label{rem:fold}
(o) Our folding order inequalities are non-strict inequalities allowing equality.

(i) The relation $\prect$ is closed.
It compares only points which lie in the same Weyl orbit,
i.e.\ if $\bar\xi_1\prec_{\taumod}\bar\xi_2$ then $W\bar\xi_1=W\bar\xi_2$. 

(ii) The relation $\prect$ on singular Weyl orbits is the closure of the relation $\prect$ on regular ones: 
It holds that $\bar\xi_1\prec_{\taumod}\bar\xi_2$ if and only if there exist sequences of regular points 
$\bar\xi_i^n\to\bar\xi_i$ 
with $\bar\xi_1^n\prec_{\taumod}\bar\xi_2^n$.

(iii) Any isometry of $\amod$ preserving $\taumod$ as a set 
preserves the relation $\prect$.

(iv) The relations $\prec_{\taumod}$ and $\prec$ are closely related:
Clearly,
$\prec$ is stronger than $\prect$,
i.e.\ if $\bar\xi_1\prec\bar\xi_2$ then $\bar\xi_1\prect\bar\xi_2$.
More precisely, 
note that a folding map $f$ fixing $\taumod$ 
is the composition $w\circ f'$ of a folding map $f'$ fixing $\simod$ 
with an element $w\in\Wt$.
Thus $\bar\xi_1\prec_{\taumod}\bar\xi_2$ if and only if 
there exist $\bar\xi'_i\in \Wt\bar\xi_i$
such that $\bar\xi'_1\prec\bar\xi'_2$. 

(v) If $\bar\xi_1\prec_{\taumod}\bar\xi_2$ then $w_1\bar\xi_1\prec_{\taumod}  w_2\bar\xi_2$
for all $w_1,w_2\in\Wt$,
because $w_1fw_2^{-1}$ is again a folding map fixing $\taumod$.
Hence, $\prect$ descends to a relation on 
the quotient $\Wt\backslash \amod$ 
which we also denote by $\prec_{\taumod}$. 
\end{rem}

There is a {\em metric estimate} for the folding order,
because folding maps are 1-Lipschitz:
\begin{lem}
\label{lem:metestfold}
If $\bar\vartheta\in\taumod$, 
then we have the implication:
\begin{equation*}
\Wt\bar\xi_1\prec_{\taumod}\Wt\bar\xi_2 
\quad\Rightarrow\quad
\tangle(\bar\xi_1,\bar\vartheta)\leq\tangle(\bar\xi_2,\bar\vartheta) 
\end{equation*}
Moreover, if $\bar\vartheta\in\interior(\taumod)$, 
then equality holds on the right hand side only if 
$\Wt\bar\xi_1=\Wt\bar\xi_2$.
\end{lem}
\proof
Suppose that $\Wt\bar\xi_1\prec_{\taumod}\Wt\bar\xi_2$.
Then there exists a folding map $\amod\to\amod$ fixing $\taumod$
with $f(\bar\xi_2)=\bar\xi_1$.
It maps the geodesic segment $\bar\vartheta\bar\xi_2$
to a broken geodesic segment $\beta$ from $\bar\vartheta$ to $\bar\xi_1$ 
of the same length,
whence the implication of inequalities.

Suppose now in addition that 
$\bar\vartheta\in\interior(\taumod)$ 
and $\Wt\bar\xi_1\neq \Wt\bar\xi_2$.
The initial segments of $\beta$ and $\bar\vartheta\bar\xi_2$ 
have the same type. 
Therefore there exists $w\in W_{\bar\vartheta}=\Wt$ 
such that $w\beta$ and $\bar\vartheta\bar\xi_2$
have a common initial segment. 
Since $w\bar\xi_1\neq\bar\xi_2$, 
the broken geodesic segment $\beta$ cannot be a true geodesic segment,
and we obtain the strict metric inequality 
$\tangle(\bar\xi_1,\bar\vartheta)<\tangle(\bar\xi_2,\bar\vartheta)$.
\qed

\medskip
As a consequence, we can justify our terminology of ``partial order":
\begin{cor}
\label{cor:ord}
$\prec_{\taumod} $ is a partial order on $\Wt\backslash \amod$. 
\end{cor}
\proof 
We must verify antisymmetry.

Suppose that $\Wt \bar\xi_1 \prec_{\taumod} \Wt \bar\xi_2\prec_{\taumod} \Wt \bar\xi_1$.
By Lemma \ref{lem:metestfold}, we have
$\tangle(\bar\xi_1,\bar\vartheta)=\tangle(\bar\xi_2,\bar\vartheta)$
for all $\bar\vartheta\in\taumod$.
The equality part implies that $\Wt \bar\xi_1 =\Wt \bar\xi_2$.
\qed

\medskip
We discuss next the structure of folding maps and decompositions into simple ones.

Each wall $m$ splits $\amod$ into two  hemispheres, 
the inner hemisphere $h^+$ containing $\simod$ 
and the outer hemisphere $h^-$. This decomposition  gives rise to the folding map 
which fixes $h^+$ and reflects $h^-$ onto it. 
We call a composition of such folding maps at walls $m_i$ 
a {\em special folding}. 
The intersection $\cap_ih^+_i$ of inner hemispheres 
is fixed by the special folding. 
In particular, special foldings fix the model chamber $\simod$.

In general, there are folding maps fixing $\simod$ which are not special. 
However, this makes no difference for the folding order $\prec$:
\begin{lem}
[{Cf. \cite[page\ 441, Thm.\ 4.9]{KM}}]
If for points $\bar\xi_1,\bar\xi_2\in \amod$ 
there exists a folding map fixing $\simod$ 
and mapping $\bar\xi_2\mapsto\bar\xi_1$, 
then there exists a special folding with this property. 
\end{lem}
\proof
We may assume that $\bar\xi_1$ and $\bar\xi_2$ are regular 
and different. We connect a point $\bar\eta$ in the interior of $\simod$ 
to $\bar\xi_2$ by a geodesic segment $\bar\ga$ 
which avoids faces of codimension at least two.
Let $f$ be a folding map fixing $\simod$ with 
$f(\bar\xi_2)=\bar\xi_1$.
Then $\bar\beta=f\circ\bar\ga$ is a broken geodesic path 
which connects $\bar\eta$ to $\bar\xi_1$ 
and has the same length and initial direction as $\bar\ga$. 
Its bending points are interior points of panels 
and $\bar\beta$ is locally ``reflected'' 
at the walls containing these panels. 
The assertion follows if we can replace $\bar\beta$ by a broken geodesic path 
from $\bar\eta$ to $\bar\xi_1$,
which is the image of $\bar\ga$ under a {\em special} folding. 

Let $\bar\eta_1$ denote the first bending point of $\bar\beta$ 
starting from $\bar\eta$. 
It lies in a wall $m_1$. 
If $\bar\beta$ crosses $m_1$ again in some point $\bar\eta_2$, 
then we replace the subpath $\bar\eta_1\bar\eta_2$ 
by its reflection at $m_1$. 
The modified broken path $\bar\beta'$ has again reflection folds, 
the same initial direction and the same endpoint. 
Moreover, its initial segment is strictly longer. 
After finitely many such modifications, 
we may assume that $\bar\beta'$ stays inside $h_1^+$. 
(The wall $m_1$ has changed in the process.) 
We then can obtain $\bar\beta'$ as the image of another broken path $\bar\beta''$ 
under the special folding $s_1$ at $m_1$, 
i.e.\ $\bar\beta'=s_1\circ\bar\beta''$, 
such that $\bar\beta''$ has a strictly longer initial segment than $\bar\beta'$. 

Thus, 
we can replace $\bar\beta$ 
by another broken path $\bar\beta''$ with reflection bends, 
with the same length and initial direction as $\bar\ga$,
with a strictly longer initial segment than $\bar\beta$,
and such that some special folding $s_1$ maps the endpoint of $\bar\beta''$ to the endpoint of $\bar\beta$.
It follows by induction that $\bar\beta$ can be replaced 
by another broken path with the same endpoint 
and which is the image of $\bar\ga$ under a special folding. 
\qed

\begin{cor}[Alternative definition of $\simod$-folding order]
\label{cor:altdeffold}
$\bar\xi_1\prec\bar\xi_2$ 
if and only if there exists a special folding 
which maps $\bar\xi_2\mapsto\bar\xi_1$. 
\end{cor}
We note that the partial order $\prec$ had been defined exactly in this way by P.~Littelmann, 
see \cite[p.\ 509]{Littelmann}. 

\medskip
The folding orders on the Weyl orbits in the model apartment 
correspond to orders on the Weyl group and its (double) coset spaces,
as we explain now.

We can regard the regular Weyl orbits as copies of $W$ 
by identifying the orbit point in the chamber $w\simod$ 
with the element $w\in W$. 
Under this identification, it holds that 
$$ w_1\prec w_2$$
for elements $w_1,w_2\in W$ 
if and only if there exists a folding map $\amod\to\amod$ fixing $\simod$ and mapping
$w_2\simod\mapsto w_1\simod$,
and $$ w_1\prect w_2$$
if and only if there exists such a folding map fixing only $\taumod$.
Again, $\prect$ descends to an order on $\Wt\backslash W$,
also denoted $\prect$.
\begin{rem}[Bruhat order]
The corollary shows 
that the folding order $\prec$ on $W$ coincides with the 
{\em Bruhat order}, 
see \cite[ch.\ 5.9]{Humphreys} or \cite[ch.\ 2]{BjoernerBrenti} 
for a definition;
hence, the folding order gives a {\em geometric interpretation} of the Bruhat order.
To verify this, 
one observes that if the chambers $w\simod$ and $w'\simod$ 
are symmetric with respect to a wall 
and if $w\simod$ lies in the inner hemisphere, 
then we have the inequality $l(w)<l(w')$ for word lengths. 
Here the word length on $W$ is defined 
using as generators the reflections at the walls of $\simod$. 
\end{rem}
More generally, 
if $\bar\xi$ is an interior point of the face $\numod\subseteq\simod$,
$\bar\xi\in\interior(\numod)$,
then $W_{\bar\xi}=\Wn$ and $W\bar\xi\cong W/\Wn$.
Under this identification,
the order $\prec_{\taumod}$ on the Weyl orbit quotient 
$\Wt\backslash W\bar\xi\subset \Wt\backslash\amod$ 
becomes a partial order on the double quotient 
$\Wt\backslash W/\Wn$, 
compare Remark~\ref{rem:fold}(v). 
It holds that 
$$\Wt w_1\Wn\prect\Wt w_2\Wn$$
if and only if 
there exist $w'_i\in \Wt w_i\Wn$
such that $w'_1\prec w'_2$,
cf.\ \cite{Mitchell2008} for a slightly different description of this order.

\medskip
We next describe the effect of the longest element $w_0\in W$ on the folding order.
Recall that $w_0$ is the involution sending $\simod$ 
to the opposite chamber $\hat\si_{mod}$ in $\amod$.
\begin{lem}
\label{lem:compposrev}
Left multiplication with $w_0$ reverses the $\simod$-folding order.
\end{lem}
\proof 
Suppose that the special folding $s_m$ at the wall $m$ maps 
$\bar\xi_2$ to $\bar\xi_1$, 
i.e.\ $s_m\bar\xi_2=\bar\xi_1$. 
When applying $w_0$, 
the inner hemisphere bounded by $m$ becomes the outer hemisphere 
bounded by $w_0m$ and vice versa, 
$w_0h^{\pm}_m=h^{\mp}_{w_0m}$. 
Hence $s_{w_0m}w_0\bar\xi_1=w_0\bar\xi_2$. 
The assertion follows by applying Corollary~\ref{cor:altdeffold} 
and induction. 
\qed

\medskip
Regarding the analogous fact for the orders $\prec_{\taumod}$,
note that 
$w_0\Wt w_0^{-1}=W_{w_0\taumod}=W_{\iota\taumod}$
and $w_0$ maps $\Wt$-orbits to $W_{\iota\taumod}$-orbits. 
The action of $w_0$ therefore induces a natural map 
\begin{equation}
\label{eq:involrelpos}
\Wt\backslash \amod
\buildrel w_0 \over\lra 
W_{\iota\taumod}\backslash \amod ,
\quad
\Wt\bar\xi\mapsto w_0\Wt\bar\xi= W_{\iota\taumod} w_0\bar\xi
\end{equation}
and, correspondingly, 
\begin{eqnarray}
\label{eq:involrelposcos}
\begin{aligned}
\Wt\backslash W/\Wn
&\buildrel w_0 \over\lra &
W_{\iota\taumod}\backslash W/\Wn &&\\
\Wt w\Wn &\mapsto& w_0\Wt w\Wn &=&W_{\iota\taumod} w_0w\Wn ,
\end{aligned}
\end{eqnarray}
and the lemma implies that these maps are order reversing:
\begin{equation}
\label{eq:ordrev}
\Wt\bar\xi_1\prec_{\taumod} \Wt\bar\xi_2 
\quad\Leftrightarrow\quad
W_{\iota\taumod}w_0\bar\xi_1\succ_{\iota\taumod} W_{\iota\taumod}w_0\bar\xi_2
\end{equation}
respectively
\begin{equation*}
\Wt w\Wn \prect \Wt w'\Wn
\quad\Leftrightarrow\quad
W_{\iota\taumod}w_0 w\Wn \succ_{\iota\taumod} W_{\iota\taumod}w_0 w'\Wn
\end{equation*}

\subsection{Relative position at infinity}
\label{sec:relpos}

Let $\si_0,\si\subset\geo X$ be chambers. 
There exists an (in general non-unique) apartment $a\subset\geo X$ containing these chambers,
$\si_0,\si\subset a$,
and a unique apartment chart $\al:\amod\to a$ such that 
$\si_0=\al(\simod)$.
We define the {\em position of $\si$ relative to $\si_0$} as the chamber 
$$ \pos(\si,\si_0) := \al^{-1}(\si)\subset\amod .$$
Abusing notation, it can be regarded algebraically as the unique element 
$$ \pos(\si,\si_0)\in W $$
such that 
$$ \si=\al\bigl(\pos(\si,\si_0)\simod\bigr) .$$
The relative position does not depend on the choice of the apartment $a$. 
To see this, choose regular points $\xi_0\in\interior(\si_0)$ and $\xi\in\interior(\si)$
which are not antipodal, $\tangle(\xi,\xi_0)<\pi$.
Then the segment $\xi_0\xi$ is contained in $a$ by convexity,
and its image $\al^{-1}(\xi_0\xi)$ in $\amod$ is independent of the chart $\al$
because its initial portion $\al^{-1}(\xi_0\xi\cap\si_0)$ in $\simod$ is.

More generally, we define the position of a chamber $\si$ 
relative to an arbitrary simplex $\tau_0$ of type $\taumod$ as follows.
Let again $a\subset\geo X$ be an apartment containing $\tau_0$ and $\si$,
and let $\al:\amod\to a$ be a chart such that 
$\tau_0=\al(\taumod)$.
It is unique up to precomposition with an element in $\Wt$.
We define the {\em position of $\si$ relative to $\tau_0$} as the $\Wt$-orbit of the chamber 
$\al^{-1}(\si)\subset\amod$. 
It can be interpreted algebraically as a coset 
$$\pos(\si,\tau_0)\in\Wt\backslash W .$$
Even more generally, 
we define the position of a simplex $\nu\subset\geo X$ relative $\tau_0$. 
Let $a\subset\geo X$ be an apartment containing $\tau_0$ and $\nu$,
and let $\al:\amod\to a$ be a chart such that 
$\tau_0=\al(\taumod)$.
We define the {\em position of $\nu$ relative to $\tau_0$} as the $\Wt$-orbit 
of the simplex $\al^{-1}(\nu)\subset\amod$.
It can be interpreted algebraically as a double coset 
$$ \pos(\nu,\tau_0)\in \Wt\backslash W/\Wn $$ 
where $\numod=\theta(\nu)$ is the type of $\nu$.
Finally, we define the position of an ideal point $\xi\in\geo X$ relative $\tau_0$
as the relative position of the simplex $\nu_{\xi}\subset\geo X$ spanned by $\xi$
(i.e.\ containing $\xi$ as an interior point),
$$ \pos(\xi,\tau_0) := \pos(\nu_{\xi},\tau_0)\in \Wt\backslash W/\Wn $$ 
where $\numod=\theta(\nu_{\xi})$.
In particular, $\pos(\xi,\tau_0)\in\Wt\backslash W$ if $\xi$ is regular.

\begin{lem}
Two ideal points $\xi_1, \xi_2$ in the same $G$-orbit $G\xi\subset\geo X$
have the same position relative to a simplex $\tau\subset\geo X$
iff they lie in the same orbit
of the parabolic subgroup $P_{\tau}<G$, 
\begin{equation*}
\pos(\xi_1,\tau)=\pos(\xi_2,\tau) \quad\Leftrightarrow\quad P_{\tau}\xi_1=P_{\tau}\xi_2 .
\end{equation*}
\end{lem}
\proof
The implication ``$\Leftarrow$'' is clear.
For ``$\Rightarrow$'',
let $a_i\subset\geo X$ be apartments containing $\tau$ and $\xi_i$.
There exists $p\in P_{\tau}$ such that $a_1=pa_2$.
Then $\pos(\xi_1,\tau)=\pos(\xi_2,\tau)=\pos(p\xi_2,\tau)$ 
iff $\xi_1$ and $p\xi_2$ span the same simplex in $a_1$.
In view of $\theta(\xi_1)=\theta(\xi_2)$,
the latter is equivalent to $\xi_1=\xi_2$.
\qed

\medskip
The positions relative $\tau$ 
thus correspond to the orbits of $P_{\tau}$
and we have the identification 
$$P_{\tau}\backslash G\xi \cong P_{\tau}\backslash\Flagn
\cong P_{\tau}\backslash G/P_{\nu_{\xi}}
\cong \Wt\backslash W/\Wn $$
with $\taumod=\theta(\tau)$ and $\numod=\theta(\nu_{\xi})$.

In particular,
for regular orbits, which are copies of the Furstenberg boundary,
we obtain the identification
$$P_{\tau}\backslash\DF X \cong \Wt\backslash W .$$
The positions relative to a chamber $\si$ 
correspond to the orbits of the minimal parabolic subgroup $B_{\si}$, 
and we have
$$B_{\si}\backslash G\xi \cong B_{\si}\backslash \Flagn
\cong B_{\si}\backslash G/P_{\nu_{\xi}} \cong W/\Wn 
\quad\hbox{ and }\quad
B_{\si}\backslash\DF X \cong W .$$

The $G$-orbits $G\xi$, respectively, the flag manifolds $\Flagn$ 
thus decompose into finitely many $P_{\tau}$-orbits 
which we call {\em Schubert strata} relative $\tau$
or {\em $\tau$-Schubert strata},
and their closures {\em (generalized) Schubert cycles}. 
(We will see below that the cycles are unions of strata.)
The level sets of $\pos(\cdot,\si)$,
i.e.\ the $B_{\si}$-orbits, 
are called {\em Schubert cells} relative $\si$.

Note that the Schubert cycles in the flag manifolds are {\em projective subvarieties}.

We will use the following notation. 
For a simplex $\tau_-\in\Flagit$, 
we denote by 
$$C_{\taumod}(\tau_-):=\{ \tau : \tau\hbox{ opposite to }\tau_- \bigr\} \subset\Flagt$$
the {\em open Schubert stratum} associated with $\tau_-$ in $\Flagt$,

\medskip
We can now use the folding order 
to compare the positions of points in a $G$-orbit $G\xi\subset\geo X$, respectively, a flag manifold $\Flagn$
relative to simplices $\tau$ of a fixed type $\taumod$.

We begin by proving a useful monotonicity property for the folding order under folding maps.
It is a direct consequence of the definition of the folding order
that folding maps $\amod\to\amod$ decrease the relative positions of pairs of simplices.
We will need the same fact for {\em folding maps} 
$\amod\to\geo X$ and $\geo X\to\geo X$
by which we mean, as before, type preserving continuous maps
sending chambers isometrically onto chambers. 
\begin{lem}[Monotonicity]
\label{lem:monot}
(i)
For a folding map $f:\geo X\to\geo X$ and simplices $\tau,\nu\subset\geo X$ it holds that 
\begin{equation*}
\pos(f(\nu),f(\tau)) \prec_{\theta(\tau)} \pos(\nu,\tau) .
\end{equation*}

(ii)
For a folding map $\al:\amod\to\geo X$, a face type $\taumod$ and a simplex $\ol\nu\subset\amod$ 
it holds that 
\begin{equation*}
\pos(\al(\ol\nu),\al(\taumod)) \prec_{\taumod} \pos(\ol\nu,\taumod) .
\end{equation*}
\end{lem}
\proof
Part (i) reduces to (ii) 
by choosing an apartment $a\supset\tau\cup\nu$ 
and a chart $\kappa:\amod\to a$ with $\kappa(\taumod)=\tau$ 
for the face type $\taumod=\theta(\tau)$. 
Then apply (ii) to $\al=f\circ\kappa$
and $\bar\nu=\kappa^{-1}(\nu)$.

To verify (ii),
consider the composition
$$\amod\stackrel{\al}{\lra}\geo X\lra\geo X/B_{\al(\simod)}\cong\amod$$
where the second map is the natural projection.
It is a folding map $\bar\al:\amod\to\amod$ 
fixing $\simod$,
and therefore
$$ \pos(\bar\al(\ol\nu),\taumod) \prec_{\taumod} \pos(\ol\nu,\taumod) .$$
Since the $B_{\al(\simod)}$-action on $\geo X$ preserves positions relative to faces of $\al(\simod)$,
we also have
$$ \pos(\al(\ol\nu),\al(\taumod)) = \pos(\bar\al(\ol\nu),\taumod) .$$
The assertion follows.
\qed

\begin{lem}[Semicontinuity of relative position]
\label{lem:semcontrelpos}
If $\xi_n\to\xi$ in $G\xi\subset\geo X$ 
and $\tau_n\to\tau$ in $\Flagt$ 
are sequences such that the sequence of relative positions $\pos(\xi_n,\tau_n)$ is constant, 
$\pos(\xi_n,\tau_n)=p\in\Wt\backslash W/W_{\theta(\nu_{\xi})}$ for all $n$,
then $\pos(\xi,\tau)\prec_{\taumod}p$.

In particular,
the sublevels of $\pos(\cdot,\tau)$ in $G\xi$ are closed.
\end{lem}
\proof
There exist apartment charts $\al_n:\amod\to\geo X$ 
with $\al_n|_{\taumod}=\kappa_{\tau_n}$ and $\al_n(\bar\xi)=\xi_n$. 
The charts subconverge to a folding map $\al$ 
with $\al|_{\taumod}=\kappa_{\tau}$ and $\al(\bar\xi)=\xi$. 
The assertion follows from monotonicity, cf.\ Lemma~\ref{lem:monot}(ii).
\qed

\medskip
It follows that the {\em suplevels}
$\{\pos(\cdot,\tau)\succ_{\taumod}p\}$ in $G\xi$ are {\em open},
because their complements are finite unions of sublevels 
$\{\pos(\cdot,\tau)\prec_{\taumod}p'\}$. 

\medskip
We show now that the folding order coincides with the inclusion order 
on Schubert cycles. 

We start with the chamber case, 
where the relation between
closures and the Bruhat order is well known: 
In the case of complex Lie groups, it goes back to the work of Chevalley in 1950s \cite{Chevalley}; 
for the proofs in the general case (including reductive groups over local fields), see \cite{Borel-Tits} and \cite{Mitchell}. 
(We are grateful to James Humphreys and Shrawan Kumar for the references.) 

\begin{prop}
\label{prop:relposint}
For a chamber $\si\subset\geo X$ 
and ideal points $\xi_1,\xi_2$ in the same $G$-orbit $G\xi\subset\geo X$, 
we have:
$$\pos(\xi_1,\si)\prec\pos(\xi_2,\si) \Leftrightarrow 
B_{\si}\xi_1\subseteq\ol{B_{\si}\xi_2}$$
\end{prop}
\proof
We denote by $\bar\xi_i\in\amod$ the point of type $\theta(\xi)$ 
with $\pos(\bar\xi_i,\simod)=\pos(\xi_i,\si)$.

Suppose first that $\xi_1\in\ol{B_{\si}\xi_2}$. 
Then there exists a sequence $(b_n)$ in $B_{\si}$ such that $b_n\xi_2\to\xi_1$. 
Let $a_n$ be apartments containing $\si$ and $b_n\xi_2$, 
and let $\al_n:\amod\to a_n$ be the apartment charts 
which restrict to the chamber chart of $\si$,
$\al_n|_{\simod}=\kappa_{\si}:\simod\to\si$. 
Then $\al_n(\bar\xi_2)=b_n\xi_2$. 
The Tits isometric embeddings $\al_n$ 
subconverge (with respect to the visual topology) to a limit map
$\al:\amod\to\geo X$. 
The map $\al$ is, in general, not an isometric embedding (chart),
but only a folding map extending $\kappa_{\si}$. 
It satisfies $\al(\bar\xi_2)=\xi_1$.
Monotonicity, cf.\ Lemma~\ref{lem:monot}(ii), yields 
$\pos(\xi_1,\si)\prec\pos(\bar\xi_2,\simod)=\pos(\xi_2,\si)$.

Vice versa, 
suppose now that $\bar\xi_1\prec\bar\xi_2$. 
By definition of the partial order $\prec$
there exists a folding map of $\amod$ fixing $\simod$
and carrying $\bar\xi_2\mapsto\bar\xi_1$. 
Furthermore there is an isometric embedding $\amod\to\geo X$ 
which extends the chamber chart $\kappa_{\si}$ 
and maps $\bar\xi_1\mapsto\xi_1$. 
By composition we obtain a folding map 
$\al:\amod\to\geo X$ 
which extends $\kappa_{\si}$ 
and maps $\al(\bar\xi_2)=\xi_1$. 
We want to find a sequence 
of isometric embeddings $\al_n:\amod\to\geo X$ extending $\kappa_{\si}$ 
such that $\al_n(\bar\xi_2)\to\al(\bar\xi_2)=\xi_1$. 
This will then imply that $\xi_1\in\ol{B_{\si}\xi_2}$. 
(Note that in general folding maps are not limits of isometric embeddings.) 

We may assume that the relative positions $\bar\xi_i$ are regular. 
(Otherwise, we may perturb them keeping the inequality 
$\bar\xi_1\prec\bar\xi_2$
and perturb the $\xi_i$ accordingly.)
We choose in $\amod$ a geodesic $\bar\ga$ of length $\pi$ 
starting in an interior point $\bar\eta_0$ of $\simod$ 
and passing through $\bar\xi_2$ 
while avoiding simplices of codimension $\geq2$. 
It crosses successively a sequence (gallery) of chambers 
$\bar\si_0=\simod,\bar\si_1,\dots,\bar\si_k=\hat\simod$
and intersects the intermediate panels 
$\bar\tau_i=\bar\si_i\cap\bar\si_{i-1}$ 
transversally in interior points $\bar\eta_i$. 
When applying the folding map $\al$, 
it may happen that successive chambers of the folded gallery coincide, 
i.e.\ that $\al(\bar\si_i)=\al(\bar\si_{i-1})$ for some $i$. 
(This happens if and only if $\al$ is not an isometric embedding.) 
One can arbitrarily well approximate 
(in the visual topology) the folded gallery by an embedded gallery
with the same initial chamber $\si$. 
To obtain such approximations it is convenient 
to use the $G$-action as follows. 
If $\al(\bar\si_i)=\al(\bar\si_{i-1})$ 
then one may pick an element $g\in G$ close to the identity, 
which fixes $\al(\bar\tau_i)$ 
and moves $\al(\bar\si_i)=\al(\bar\si_{i-1})$ away from itself, 
and apply it to the ``tail'' $\al(\bar\si_i),\dots,\al(\bar\si_k)$ 
of the gallery. 
Doing this inductively along the gallery, 
one obtains an arbitrarily good approximation 
of the folded gallery $\al(\bar\si_0)=\si,\dots,\al(\bar\si_k)$ 
by an embedded gallery $\si_0=\si,\dots,\si_k$, 
that is a sequence of chambers 
such that $\si_i\cap\si_{i-1}$ is precisely a panel for all $i$. 
This yields at the same time an approximation 
of the broken geodesic $\al(\bar\ga)$ in $\geo X$ by a true geodesic $\ga$ 
such that $\ga\cap\si_i$ and $\bar\ga\cap\bar\si_i$ 
are corresponding subsegments of the same type. 
Now we use the path $\bar\ga$ as a ``guiding line'' 
to extend the correspondence $\bar\si_i\mapsto\si_i$ of galleries 
to an isometric embedding $\al':\amod\to\geo X$ extending $\kappa_{\si}$: 
Since $\ga$ connects two antipodal regular points 
there exists a unique such $\al'$ extending the isometry $\bar\ga\to\ga$ 
and hence mapping $\bar\si_i$ to $\si_i$. 
By construction, $\al'(\bar\xi_2)$ approximates $\al(\bar\xi_2)=\xi_1$ 
arbitrarily well. 
So we find a sequence of apartment charts $\al_n$ with the desired properties. 
\qed

\medskip
The proposition readily generalizes to the simplex case 
(in the case when $G$ is a complex semisimple Lie group, 
a proof of the following proposition can be found in \cite[Prop.\ 3.13]{Mitchell2008}):
\begin{prop}
\label{prop:relposinttau}
For simplices $\tau\subset\geo X$ and ideal points $\xi_1,\xi_2$ in the same $G$-orbit $G\xi\subset\geo X$, we have:
$$\pos(\xi_1,\tau)\prec_{\theta(\tau)}\pos(\xi_2,\tau) \Leftrightarrow 
P_{\tau}\xi_1\subseteq\ol{P_{\tau}\xi_2}$$
\end{prop}
\proof
Let $\si\supset\tau$ be a chamber. 
Since the quotient space $P_{\tau}/B_{\si}$ is compact 
(it is the space of chambers containing $\tau$ as a face), 
the condition 
$\xi_1\in\ol{P_{\tau}\xi_2}$ 
is equivalent to the existence of an element $p\in P_{\tau}$ 
such that 
$p\xi_1\subseteq\ol{B_{\si}\xi_2}$. 
According to Proposition~\ref{prop:relposint}, 
this is equivalent to the existence of $p\in P_{\tau}$ such that 
$\pos(p\xi_1,\si)\prec\pos(\xi_2,\si)$. 
Since $P_{\tau}$ acts transitively on chambers containing $\tau$, 
we have that 
$\cup_{p\in P_{\tau}}\pos(p\xi_1,\si)=W_{\theta(\tau)}\pos(\xi_1,\si)
=\pos(\xi_1,\tau)$. 
This completes the proof.
\qed

\medskip
In other words, 
the proposition says that the $\tau$-Schubert cycles in $G/P_{\nu}$ 
correspond to the sublevels of the folding order $\prec_{\taumod}$ on  $\Wt\backslash W/\Wn$,
where $\taumod=\theta(\tau)$ and $\numod=\theta(\nu)$.

\medskip
Recall that the simplices opposite to simplices of type $\taumod$ 
have type $\iota\taumod$,
and that the action of $w_0$ induces the natural maps
\begin{equation*}
\Wt\backslash \amod
\buildrel w_0 \over\lra 
W_{\iota\taumod}\backslash \amod ,
\end{equation*}
and hence the maps 
\begin{equation*}
\Wt\backslash W/\Wn
\buildrel w_0 \over\lra 
W_{\iota\taumod}\backslash W/\Wn
\end{equation*}
of relative positions,
compare (\ref{eq:involrelpos}) and (\ref{eq:involrelposcos}).
\begin{dfn}[Complementary position]
\label{dfn:complemtau}
We define the {\em complementary position} by
$$\cpos:=w_0\pos .$$
\end{dfn}
This terminology is justified by
(cf.\ Def.~\ref{dfn:antip}(ii) for the notion of antipodality):

\begin{lem}
\label{lem:justrelpos}
Let $\tau,\hat\tau,\nu\subset\geo X$ be simplices contained in an apartment $a$,
and suppose that $\tau$ and $\hat\tau$ are antipodal.
Then $\pos(\nu,\hat\tau)=\cpos(\nu,\tau)$.
\end{lem}
\proof
Let $\al:\amod\to a$ be a chart such that 
$\al|_{\taumod}=\kappa_{\tau}$.
Then $\hat\tau=\al(\htaumod)=(\al\circ w_0)(\iota\taumod)$ 
with the simplex $\htaumod=w_0\iota\taumod\subset\amod$ opposite to $\taumod$.
Using the reparametrized chart $\al\circ w_0$, we obtain
$\pos(\nu,\hat\tau)=(\al\circ w_0)^{-1}(\nu)=w_0\al^{-1}(\nu)=w_0\pos(\nu,\tau)$.
\qed

\medskip
The relation of ``complementarity'' is clearly symmetric,
$\ccpos=\pos$.
Passing to complementary relative position reverses the partial order,
cf.\ Lemma~\ref{lem:compposrev}:
\begin{equation}
\label{ineq:revordpos}
\pos(\xi_1,\tau)\prec_{\theta(\tau)}\pos(\xi_2,\tau) 
\quad\Leftrightarrow\quad
\cpos(\xi_1,\tau)\succ_{\iota\theta(\tau)}\cpos(\xi_2,\tau) 
\end{equation}

\medskip
Points with smaller position relative to a simplex are closer to it 
in a {\em metric} sense. 
Namely, 
according to Lemma~\ref{lem:metestfold}
we have the inequality of Tits distances 
\begin{equation}
\label{ineq:ordrelpos1tau}
\pos(\xi_1,\tau)\prec_{\theta(\tau)}\pos(\xi_2,\tau) 
\quad\Ra\quad
\tangle(\xi_1,\cdot)|_{\tau}\leq\tangle(\xi_2,\cdot)|_{\tau} ,
\end{equation}
respectively, 
for simplices $\tau_1$ and $\tau_2$ of the same type $\taumod$, 
\begin{equation}
\label{ineq:ordrelpos2tau}
\pos(\xi_1,\tau_1)\prec_{\taumod}\pos(\xi_2,\tau_2) 
\quad\Ra\quad
\tangle(\xi_1,\cdot)\circ\kappa_{\tau_1}
\leq\tangle(\xi_2,\cdot)\circ\kappa_{\tau_2} .
\end{equation}
If $\tau_1$ and $\tau_2$ are simplices of opposite types, 
$\theta(\tau_1)=\iota\theta(\tau_2)$, 
and if the relative positions 
$\pos(\xi_1,\tau_1)$ and $\pos(\xi_2,\tau_2)$
are complementary, 
then
\begin{equation}
\label{eq:complementrelpostau}
\tangle(\xi_1,\cdot)\circ\kappa_{\tau_1}+
\tangle(\xi_2,\cdot)\circ\kappa_{\tau_2}\circ\iota\equiv\pi
\end{equation}
on $\theta(\tau_1)=\iota\theta(\tau_2)$. 
To see this, 
note that the formula reduces to the case when 
the simplices $\tau_1$ and $\tau_2$ are opposite to each other
and $\xi_1=\xi_2$ lies in an apartment containing them.

\medskip
The following triangle inequality extends Lemma~\ref{lem:justrelpos}:
\begin{lem}
\label{lem:triai}
Let $\tau,\hat\tau\subset\geo X$ be a pair of antipodal simplices 
and let $\nu\subset\geo X$ be an arbitrary simplex.
Then 
\begin{equation*}
\pos(\nu,\hat\tau)\succ_{\iota\theta(\tau)}\cpos(\nu,\tau)
\end{equation*}
with equality iff 
$\tau,\hat\tau,\nu$ are contained in an apartment.\footnote{Equivalently, 
$\nu$ lies in the spherical subbuilding $B(\tau,\hat\tau)\subset\geo X$
consisting of all apartments containing $\tau,\hat\tau$.}
\end{lem}
\proof
Let $a\subset\geo X$ be an apartment containing $\tau,\hat\tau$
and let 
$$\geo X\stackrel{r}{\lra}a$$
be a folding retraction,
i.e.\ a folding map such that $r|_a=\id_a$.
Such a retraction is given e.g.\ by the natural projection 
$\geo X\to\geo X/B_{\si}\cong a$ for a chamber $\si\subset a$.
By monotonicity, cf.\ Lemma~\ref{lem:monot}(i), we have
$$ \pos(\nu,\tau) \succ_{\theta(\tau)} \pos(r\nu,\tau) $$
and 
$$ \pos(\nu,\hat\tau) \succ_{\iota\theta(\tau)} \pos(r\nu,\hat\tau)=\cpos(r\nu,\tau) ,$$
cf.\ Lemma~\ref{lem:justrelpos}.
Since complementing position reverses the folding order, see (\ref{ineq:revordpos}), 
we obtain the desired inequality.

Suppose that equality holds,
$\pos(\nu,\hat\tau)=\cpos(\nu,\tau)$.
Let $\xi,\hat\xi,\eta$ be interior points of $\tau,\hat\tau,\nu$
such that $\xi,\hat\xi$ are antipodal.
Then 
$$\tangle(\xi,\eta)+\tangle(\eta,\hat\xi)=\pi ,$$
cf.\ Lemma~\ref{lem:justrelpos} again,
i.e.\ $\eta$ lies on a geodesic segment $\xi\hat\xi$.
It follows that there exists an apartment containing $\xi,\hat\xi,\eta$
and hence also the simplices spanned by these points. 
\qed

\subsection{Thickenings}
\label{sec:chthick}

\subsubsection{Thickenings in the Weyl group}

A {\em thickening} (of the neutral element) in $W$ is a subset $$\Th\subset W$$
which is a union of sublevels for the folding order,
i.e. which contains with every element $w$ also every element $w'$ satisfying $w'\prec w$.
In the theory of posets, such subsets are called {\em ideals}.

Unions and intersections of thickenings are again thickenings, 
and removing a maximal element from a thickening yields a thickening.
Furthermore, 
note that 
$$ \Th^c := w_0(W-\Th) = W -w_0\Th$$
is again a thickening,
because left multiplication with $w_0$ reverses the folding order, 
cf.\ Lemma~\ref{lem:compposrev}.
It holds that 
$$ W =\Th\sqcup w_0\Th^c $$
and we therefore call $\Th^c$ the thickening {\em complementary} to $\Th$.
\begin{dfn}[Fat and slim]
\label{dfn:fatslimthtau}
The thickening 
$\Th\subset W$ 
is called {\em fat} if $\Th\cup w_0\Th=W$, equivalently, $\Th\supseteq\Th^c$.
It is called {\em slim} if $\Th\cap w_0\Th=\emptyset$, equivalently, $\Th\subseteq\Th^c$.
It is called {\em balanced} if it is both fat and slim, equivalently, $\Th=\Th^c$.
\end{dfn}
For types $\bar\vartheta_0,\bar\vartheta\in\simod$ and a radius $r\in[0,\pi]$
we define the {\em metric thickening}
\begin{equation}
\label{eq:metth}
\Th_{\bar\vartheta_0,\bar\vartheta,r} := \{w\in W:d(w\bar\vartheta,\bar\vartheta_0)\leq r\},
\end{equation}
using the natural $W$-invariant spherical metric $d$ on $\amod$.
It is indeed a thickening by Lemma~\ref{lem:metestfold}.

Recall that for a face type $\taumod\subseteq\simod$, we denote by $\Wt$ its stabilizer in $W$. 
Furthermore, $\iota=-w_0:\simod\to\simod$ denotes the canonical involution of the model spherical Weyl chamber. 

\begin{lem}
\label{lem:metth}
(i) If $\bar\vartheta_0\in\taumod$, then 
$\Th_{\bar\vartheta_0,\bar\vartheta,r}$ is $\Wt$-left invariant,
$\Wt\Th_{\bar\vartheta_0,\bar\vartheta,r}=\Th_{\bar\vartheta_0,\bar\vartheta,r}$.

(ii) If also $\iota\bar\vartheta_0=\bar\vartheta_0$, then 
$\Th_{\bar\vartheta_0,\bar\vartheta,r}$ is fat for $r\geq\pihalf$ and slim for $r<\pihalf$.

(iii) If in addition 
$d(w\bar\vartheta,\bar\vartheta_0)\neq\pihalf$ for all $w\in W$,
then $\Th_{\bar\vartheta_0,\bar\vartheta,\pihalf}$ is balanced.
\end{lem}
\proof
(i) For $w'\in \Wt$, we have that $w'\bar\vartheta_0=\bar\vartheta_0$ and hence 
$$ d(w'w\bar\vartheta,\bar\vartheta_0) = 
d(w\bar\vartheta,\underbrace{{w'}^{-1}\bar\vartheta_0}_{\bar\vartheta_0}).$$
(ii) Since $w_0\bar\vartheta_0=-\iota\bar\vartheta_0=-\bar\vartheta_0$, we have 
$$ d(w_0w\bar\vartheta,-\bar\vartheta_0) = 
d(w\bar\vartheta,\underbrace{-w_0\bar\vartheta_0}_{\bar\vartheta_0})$$
and 
$$ d(w_0w\bar\vartheta,\bar\vartheta_0) = \pi- d(w\bar\vartheta,\bar\vartheta_0).$$
Hence 
$$ \Th_{\bar\vartheta_0,\bar\vartheta,r}^c =
W-w_0\Th_{\bar\vartheta_0,\bar\vartheta,r} := \{w\in W:d(w\bar\vartheta,\bar\vartheta_0)<\pi- r\}$$
which yields the assertion. 

(iii) Slimness holds because 
$\Th_{\bar\vartheta_0,\bar\vartheta,\pihalf}=\Th_{\bar\vartheta_0,\bar\vartheta,r}$ for radii $r$ slightly below $\pihalf$.
\qed

\medskip
The metric examples provide balanced thickenings with arbitrary left invariance:
\begin{cor}[Existence of balanced thickenings I]
\label{cor:exbaltau}
For every $\iota$-invariant face type $\taumod$ 
there exists a $\Wt$-left invariant balanced thickening $\Th\subset W$.
\end{cor}
\proof
Since $\iota\taumod=\taumod$,
there exists $\bar\vartheta_0\in\taumod$ 
such that $\iota\bar\vartheta_0=\bar\vartheta_0$.
Moreover, 
the set of types $\bar\vartheta\in\simod$ 
such that $d(\cdot\bar\vartheta,\bar\vartheta_0)\neq\pihalf$ on $W$
is the complement of a finite union of great spheres in $\amod$,
and hence open and dense.
\qed

\medskip
In order to obtain balanced thickenings with additional right invariance,
we modify the metric thickenings (\ref{eq:metth}) at their ``boundaries". 
The rigidity part of Lemma~\ref{lem:metestfold} implies that the elements of 
\begin{equation*}
\D\Th_{\bar\vartheta_0,\bar\vartheta,r} :=
\{w\in W:d(w\bar\vartheta,\bar\vartheta_0)=r\},
\end{equation*}
are pairwise $\prec$-incomparable
and maximal in $\Th_{\bar\vartheta_0,\bar\vartheta,r}$. 
Therefore every subset $\Th\subset W$ with
\begin{equation*}
\{w\in W:d(w\bar\vartheta,\bar\vartheta_0)<r\} \subseteq \Th\subseteq \Th_{\bar\vartheta_0,\bar\vartheta,r}
\end{equation*}
is a thickening.

Using these modified metric thickenings,
we can generalize our last existence result:
\begin{prop}[Existence of balanced thickenings II]
\label{prop:exbaltau}
Let $\taumod,\numod\subseteq\simod$ be face types and suppose that $\taumod$ is $\iota$-invariant.
Then a $\Wt$-left invariant and $\Wn$-right invariant balanced thickening $\Th\subset W$
exists if and only if (left multiplication by) $w_0$ has no fixed point on 
$\Wt\backslash W/\Wn$,
cf.\ (\ref{eq:involrelposcos}).
\end{prop}
\proof
If a balanced thickening exists, then $w_0$ cannot have a fixed point
as a consequence of the definition of balancedness.

Vice versa, 
let us 
assume that $w_0$ has no fixed point.
We choose $\bar\vartheta_0\in\taumod$ 
with $\iota\bar\vartheta_0=\bar\vartheta_0$
and $\bar\nu\in\numod$.
Then the fat thickening $\Th_{\bar\vartheta_0,\bar\nu,\pihalf}$ 
is $\Wt$-left and $\Wn$-right invariant,
and so is the ``great sphere"
$\D\Th_{\bar\vartheta_0,\bar\nu,\pihalf}$.
The latter is moreover preserved by the involution $w_0$ 
while the ``hemispheres bounded by it",
$\Th_{\bar\vartheta_0,\bar\nu,\pihalf}-\D\Th_{\bar\vartheta_0,\bar\nu,\pihalf}$
and $W-\Th_{\bar\vartheta_0,\bar\nu,\pihalf}$,
are exchanged, 
cf.\ the proof of Lemma~\ref{lem:metth}(ii).
Since $w_0$ has no fixed point,
$\D\Th_{\bar\vartheta_0,\bar\nu,\pihalf}$
decomposes as a union of pairs of double cosets $\Wt w\Wn$
which are swapped by $w_0$.
By removing from $\Th_{\bar\vartheta_0,\bar\nu,\pihalf}$
one double coset of each pair,
we therefore obtain a balanced thickening as desired.
\qed

\medskip
For instance, we can deduce:
\begin{cor}
\label{cor:bthwid}
If $w_0=-\id_{\amod}$, 
then a $\Wn$-right invariant balanced thickening 
exists for every face type $\numod$.
\end{cor}
\proof
We equivariantly identify the coset space $W/\Wn$
with an orbit $W\bar\nu\subset\amod$
for some $\bar\nu\in\inte(\numod)$.
By assumption, $w_0$ has no fixed point on $\amod$,
and hence none on $W/\Wn$. 
\qed

\begin{rem}
\label{rem:exbaltaunu}
(i) Note that 
$w_0=-\id_{\amod}$ if and only if 
all irreducible factors of $W$ are of type 
$A_1$, $B_{n\geq2}$, $D_{2k\geq4}$, $E_{7,8}$, $F_4$ or $G_2$, 
see \cite{Bourbaki}. 

(ii) $\Wt$-left and $\Wn$-right invariant balanced thickenings 
do not always exist.
For instance, in the $B_2$-case there are no $\Wt$-biinvariant thickenings 
for $\taumod$ a vertex type.
\end{rem}

In rank two, the balanced thickenings are easy to describe:
\begin{example}[Balanced thickenings in rank 2]
\label{ex:balth2w}
(i) If $W=W_{A_2}$,
then $\simod$ is an arc of length $\pithird$.
There is a unique balanced thickening $\Th\subset W$
described by the property that $\Th\cdot\simod\subset\amod$ 
is the $\pihalf$-ball centered at the midpoint of $\simod$.

(ii) If $W=W_{B_2}$ or $W_{G_2}$,
then $\simod$ is an arc of length $\piforth$ or $\pisixth$.
In these cases, there are two balanced thickenings.
Namely, for each vertex $\bar\xi$ of $\simod$
we have the $W_{\bar\xi}$-left invariant thickening $\Th\subset W$ for which 
$\Th\cdot\simod=\ol B(\bar\xi,\pihalf)$.
\end{example}

Below, we give two examples in higher rank.
First in the irreducible case:
\begin{example}[Some balanced thickenings  of type $A_n$]
\label{ex:balthhrw}
The spherical Coxeter complex $\amod$ can be modelled as the unit sphere in the hyperplane 
$$ x_0+\ldots +x_n=0 $$
in $\R^{n+1}$.
The Weyl group $W\cong S_{n+1}$ acts by permuting the coordinates,
and we choose the fundamental chamber $\simod\subset\amod$ as given by the inequalities
$x_0\geq\ldots\geq x_n$. 
It holds that
$$ (x_0,\ldots,x_n)\stackrel{w_0=-\iota}{\mapsto} (x_n,\ldots,x_0) ,$$
There are the $\iota$-invariant edge midpoints $\bar\vartheta_k\in\simod$ for $1\leq k\leq\frac{n}{2}$ 
represented by the vectors 
$$ (\underbrace{1,\ldots, 1}_{k},0,\ldots, 0, \underbrace{-1, \ldots, -1}_{k}) $$
and the unique $\iota$-invariant vertex $\bar\vartheta_{\frac{n+1}{2}}\in\simod$
if $n$ is odd.
The type $\bar\vartheta_1$ is the {\em unique root type}.

In incidence geometric terms, 
the Coxeter complex $\amod$ is the spherical building associated 
to the finite projective $n$-space 
$\P^n_{mod}=\{e_0,\dots,e_n\}$
consisting of $n+1$ points.
Every subset 
is a projective subspace (of dimension on less than the number of points in it)
and corresponds to a vertex of $\amod$.
Vertices are adjacent iff the corresponding subspaces are incident, 
i.e.\ one contains the other.
The element $w_{\pi}\in W$ corresponding to the permutation $\pi\in S_{n+1}$ acts by 
$e_i\stackrel{w_{\pi}}{\mapsto}e_{\pi(i)}$.
We let the fundamental chamber $\simod\subset\amod$
correspond to the full flag
$$\{e_0\}\subset\ldots\subset\{e_0,\ldots,e_i\}\subset\ldots\subset\{e_0,\ldots,e_{n-1}\}$$
The edge spanned by $\bar\vartheta_k$ then has the vertices 
$\{e_0,\ldots,e_{k-1}\}$ and $\{e_0,\ldots,e_{n-k}\}$.

We determine the metric thickenings 
$\Th_{\bar\vartheta_1,\bar\vartheta,\pihalf}$
for regular types $\bar\vartheta\in\inte(\simod)$:
The type $\bar\vartheta$ is represented by a vector $(t_i)$ with 
$t_0>\ldots>t_n$.
The element $w_{\pi}\in W$ 
carries $(t_i)$ to the vector $(t_{\pi^{-1}(i)})$.
Thus, 
$\tangle(w_{\pi}\bar\vartheta,\bar\vartheta_1)<\pihalf$ 
if and only if $t_{\pi^{-1}(0)}>t_{\pi^{-1}(n)}$ if and only if $\pi^{-1}(0)<\pi^{-1}(n)$,
and we obtain the balanced thickening:
\begin{equation}
\label{eq:anthrtty}
\Th_{\bar\vartheta_1,\bar\vartheta,\pihalf} = \{w_{\pi}\in W: \pi^{-1}(0)<\pi^{-1}(n) \} 
\end{equation}
Similarly, one can describe 
the thickenings $\Th_{\bar\vartheta_k,\bar\vartheta,\pihalf}$ for $k\geq2$.
They depend on the type $\bar\vartheta$ and are balanced for 
a dense open set of values.

To give an incidence geometric description of the thickening (\ref{eq:anthrtty}),
note that $\pi^{-1}(i)$ is the dimension of the smallest subspace in the flag $w_{\pi}\simod$ 
which contains $e_i$.
Hence $w_{\pi}\in\Th_{\bar\vartheta_1,\bar\vartheta,\pihalf}$
if and only if  to the flag $w_{\pi}\simod$ belongs a subspace $U$ 
which contains $e_0$ but not $e_n$,
equivalently 
\begin{equation}
\label{eq:anthrttyinc}
\{e_0\}\subseteq U \subseteq\{e_0,\ldots,e_{n-1}\} .
\end{equation}

Another interesting example is the $W_{\bar\vartheta_m}$-biinvariant thickening 
$\Th_{\bar\vartheta_m,\bar\vartheta_m,\pihalf}$ for $n=2m-1$.
We observe that 
$\tangle(w_{\pi}\bar\vartheta_m,\bar\vartheta_m)\leq\pihalf$ 
if and only if 
\begin{equation}
\label{eq:anthrmdinc}
\big|w_{\pi}\{e_0,\ldots,e_{m-1}\} \cap \{e_0,\ldots,e_{m-1}\} \big| \geq \frac{m}{2}
\end{equation}
with the equality cases corresponding to each other.
Equality cannot occur if $m$ is odd, 
and in this case the thickening 
$\Th_{\bar\vartheta_m,\bar\vartheta_m,\pihalf}$
is balanced.
\end{example}
\begin{rem}
With a bit more work one can classify all balanced thickenings in the $A_3$ case: 
There are 10 balanced thickenings. 
Two of them are $W_{\bar\vartheta_2}$-left invariant 
for the unique $\iota$-invariant vertex $\bar\vartheta_2\in\simod$,
and one is $W_{\epsmod}$-left invariant 
for the unique $\iota$-invariant edge $\epsmod\subset\simod$. 
\end{rem}

The next example concerns the reducible case:
\begin{example}[Some balanced thickenings  of type $A_1^n$]
\label{ex:product_case}
The spherical Coxeter complex $\amod$ can be modelled as the unit sphere in $\R^n$.
The Weyl group $W\cong \Z_2^n\cong\{\pm1\}^n$ 
acts by changing the signs of the coordinates $x_i$,
i.e.\ its canonical generators act by reflections at the coordinate hyperplanes.
We choose the fundamental chamber $\simod\subset\amod$ as given by the inequalities
$x_1,\ldots,x_n\geq0$. 

The longest element $w_0=(-1,\ldots, -1)$ acts as $-\id$.
The Bruhat order on $W$ is given by
$$ w_{\eps}\prec w_{\eps'} \quad \iff \quad \eps_i\ge \eps'_i \quad\forall\;i $$
where we denote the elements in $W$ by $w_{\eps}$ with $\eps=(\eps_i)$.

The $(k-1)$-simplices of the spherical Coxeter complex $\amod$ 
correspond to the $\{\pm1\}$-valued maps defined on subsets of $\{1,\ldots,n\}$ of cardinality $k$.
In particular, the chambers can be interpreted as the ordered $n$-point configurations on $\{\pm1\}$. 

Let $\bar\zeta\in\simod$ be the {\em central} type 
represented by the vector $(1,\ldots, 1)$.
We determine the metric thickenings of the form
$\Th_{\bar\zeta,\bar\vartheta,\pihalf}$
for the regular types $\bar\vartheta\in\inte(\simod)$:
The type $\bar\vartheta$ is represented by a vector $t=(t_i)$ with $t_i>0$.
The element $w_{\eps}\in W$ 
carries $t$ to the vector $(\eps_it_i)$.
Thus, 
$\tangle(w_{\eps}\bar\vartheta,\bar\zeta)<\pihalf$ 
if and only if 
$\eps\cdot t=\eps_1t_1+\ldots+\eps_nt_n >0$,
and 
$$
\Th_t=\{w_{\eps}\in W: \eps\cdot t> 0\}, \quad 
\ol{\Th}_t=\{w_{\eps}\in W: \eps\cdot t\geq 0\}
$$
are {\em metric thickenings}. 
The thickening $\Th_t$ is slim, while $\ol{\Th}_t$ is fat. 
We have that $\ol{\Th}_t=\Th_t$ is balanced,
iff $\eps\cdot t\neq0$ for all sign choices $\eps$,
which is the case for ``generic'' values of $t$.

To phrase it in terms of configurations,
we consider {\em weighted $n$-point configurations} on $\{\pm1\}$ 
with weights $t_i$.
Then the thickenings $\ol{\Th}_t$ and $\Th_t$ 
correspond to the sets of configurations with at least, respectively, strictly more than 
half of the total mass placed on $+1$.
\end{example}

\subsubsection{Thickenings at infinity}
\label{sec:thickinf}

From thickenings in the Weyl group, 
we derive {\em thickenings at infinity} as follows.

Given a $\Wt$-left invariant thickening $\Th\subset W$,
the induced thickening of a simplex $\tau\in\Flagt$ 
inside the Furstenberg boundary
$$ \ThF(\tau) := \{ \pos(\cdot,\tau)\in \Wt\backslash\Th \} \subset\DF X $$
is well-defined.
Furthermore, 
we define the thickening of $\tau$ inside the visual boundary 
as the union of the corresponding (closed) chambers
\begin{equation}
\label{eq:metthtau}
\Th(\tau) := \bigcup_{\si\in\ThF(\tau)} \si \subset\geo X .
\end{equation}
Due to the semicontinuity of relative position, cf.\ Lemma~\ref{lem:semcontrelpos},
the thickenings $\ThF(\tau)$ and $\Th(\tau)$ are compact.
The intersections $\Th(\tau)\cap G\xi$ with $G$-orbits $G\xi\subset\geo X$ 
are finite unions of Schubert cycles and hence projective subvarieties. 
For regular $G$-orbits $G\xi$, the intersection $\Th(\tau)\cap G\xi$ is naturally identified with $\ThF(\tau)$.

Note that 
if the thickening $\Th\subset W$ is $\Wn$-{\em right} invariant for a face type $\numod\subseteq\simod$,
then the thickenings of simplices are {\em unions of stars} of simplices of type $\numod$.

For a subset $A\subset\Flagt$, we define the induced thickenings
$$ \ThF(A)= \bigcup_{\tau\in A} \ThF(\tau) 
\quad\hbox{ and }\quad
\Th(A)= \bigcup_{\tau\in A} \Th(\tau) .$$
If $A$ is compact, then its thickenings are compact as well.

\medskip 
Below are several examples of thickenings,
based on the examples in the previous section.
\begin{example}[Rank 2]
\label{ex:balth2i}
We continue with Example~\ref{ex:balth2w}.

(i) If $W=W_{A_2}$,
then chambers in the visual boundary are arcs of Tits length $\pithird$.
For the unique balanced thickening $\Th\subset W$,
the associated thickening $\Th(\si)\subset\geo X$ of a chamber $\si\subset\geo X$ with midpoint $\zeta$
is the ball $\ol B(\zeta,\pihalf)$.
In incidence geometric terms,
regarding $\tits X$ as the spherical building associated to a projective plane $\Pi$, 
the chamber $\si$ corresponds to a flag $(l,p)$ consisting of a line $l\subset\Pi$ and a point $p\in l$.
The thickening $\ThF(\si)\subset\DF X$ inside the Furstenberg boundary 
consists of all flags $(l',p')$ such that $l'=l$ or $p'=p$.

(ii) If $W=W_{B_2}$ or $W_{G_2}$,
then chambers have length $\piforth$ or $\pisixth$.
For a vertex type $\bar\xi\in\simod$
and the unique $W_{\bar\xi}$-left invariant balanced thickening $\Th\subset W$,
the associated thickening of a vertex $\xi\in\geo X$ of type $\theta(\xi)=\bar\xi$ inside $\geo X$
is given by $\Th(\xi)=\ol B(\xi,\pihalf)$.
The thickening of a chamber $\si\subset\geo X$ equals
$\Th(\si)=\ol B(\xi_{\si},\pihalf)$ 
where $\xi_{\si}$ is the vertex of $\si$ with type $\bar\xi$.

For instance, if $G=O(n,2)$ with $n\geq2$ and hence $X$ has type $B_2$,
then $\tits X$ can be regarded from the incidence geometry perspective 
as the spherical building arising from isotropic flags in $\R^{n,2}=\R^n\oplus \R^2$.
A chamber corresponds to a flag $(L,P)$ consisting of an isotropic plane $P\subset\R^{n,2}$
and a(n isotropic) line $L\subset P$.
If $\bar\xi$ is the vertex type corresponding to isotropic {\em lines},
then the thickening $\ThF(L)\subset\DF X$ of an isotropic line $L\in\Flag_{\bar\xi}$ 
consists of all flags $(L',P')$ such that $P'\supset L$. 
On the other hand,
if $\bar\xi$ is the vertex type corresponding to isotropic {\em planes},
then the thickening $\ThF(P)$ of an isotropic plane $P\in\Flag_{\bar\xi}$ 
consists of all flags $(L',P')\in\DF X$ such that $L'\subset P$. 
\end{example}

\begin{example}[Type $A_n$]
\label{ex:balthhri}
We continue with Example~\ref{ex:balthhrw}.

Let $G=SL(n+1,\F)$.
We regard $\tits X$ as the spherical building associated to the projective $n$-space 
$\F P^n$.

Let $\taumod(\bar\vartheta_1)\subset\simod$ denote the edge type with midpoint $\bar\vartheta_1$.
Then $\Flag_{\taumod(\bar\vartheta_1)}={\mathcal F}_{1,n}$,
the manifold of 2-flags $(L,H)$ consisting of a hyperplane $H\subset\F^{n+1}$
and a line $L\subset H$.
For the $W_{\taumod(\bar\vartheta_1)}$-left invariant balanced thickening 
$\Th_{\bar\vartheta_1,\bar\vartheta,\pihalf}\subset W$
given by (\ref{eq:anthrtty}), 
the thickening $\Th_{\bar\vartheta_1,\bar\vartheta,\pihalf}^{F\ddot u}((L,H))$ of the flag $(L,H)\in{\mathcal F}_{1,n}$
in the full flag manifold $\DF X$
consists of all flags $U_1\subset\ldots\subset U_i\subset\ldots\subset U_n$ in $\F^{n+1}$
such that 
$$ L\subseteq U_i\subseteq H \quad\hbox{ for some $i$},$$
compare (\ref{eq:anthrttyinc}).

If $n=4l+1$,
then $\Flag_{\bar\vartheta_{2l+1}}={\mathcal F}_{2l+1}$ is the middle Grassmannian 
of $(2l+2)$-dimensional linear subspaces of $\F^{n+1}$. 
The balanced $W_{\bar\vartheta_{2l+1}}$-biinvariant thickening 
$\Th_{\bar\vartheta_{2l+1},\bar\vartheta_{2l+1},\pihalf}^{F\ddot u}(U)$ of a subspace $U\in{\mathcal F}_{2l+1}$
inside ${\mathcal F}_{2l+1}$
consists of all subspaces $U'\in{\mathcal F}_{2l+1}$ such that 
$$ \dim(U'\cap U) \geq l+1 ,$$
compare (\ref{eq:anthrmdinc}).
\end{example}

\begin{example}[Type $A_1^n$, configuration spaces and stability in the sense of Geometric Invariant Theory]
\label{ex:config}
We continue with Example~\ref{ex:product_case}. 

Let $X=Y^n$ be the $n$-fold product of a rank one symmetric space $Y$, e.g. $Y=\H^2$.
Then $\DF X\cong(\geo Y)^n$
and we will view chambers as ordered $n$-point configurations $\xi=(\eta_i)$ on $\geo Y$.
The relative position of two configurations $\xi,\xi'\in\DF X$ is given by:
$$\pos(\xi',\xi)=w_{\eps} \quad\hbox{ with $\eps_i=+1\Leftrightarrow\eta_i'=\eta_i$} $$ 
Thus, it records the entries $i$ where the configurations agree.

We fix a regular vector $t=(t_i)\in\inte(\De)\cong\R_+^n$
and assign the weight $t_i>0$ to the $i$-th point of a configuration. 
Then a chamber $\xi$, regarded now as a {\em weighted configuration} on $\geo Y$, 
defines the finite measure 
$$ \mu_{\xi}= t_1 \delta_{\eta_1}+\ldots+t_n \delta_{\eta_n} $$ 
on $\geo Y$,
where $\delta_{\eta_i}$ denotes the Dirac measure concentrated in the point $\eta_i$. 
(Masses add when points $\eta_i$ ``collide''). 
The total mass $|\mu_{\xi}|$ of $\mu_{\xi}$ equals
$$
M=t_1+\ldots+t_n. 
$$
In the language of Geometric Invariant Theory,
the finite measure $\mu_{\xi}$ 
(and the corresponding weighted configuration $\xi$) 
is called  {\em stable} if 
$\mu_{\xi}(\eta)<M/2$
for all points $\eta\in\geo Y$,
and {\em semistable} if 
$\mu_{\xi}(\eta)\leq M/2$ for all $\eta$.

Let $(\eta):=(\eta,\ldots,\eta)$ denote the configuration concentrated in the point $\eta$.
According to Example~\ref{ex:product_case}, 
the thickenings 
$(\ol{\Th}_t)_{F\ddot u}((\eta))$  and $(\Th_t)_{F\ddot u}((\eta))$ of $(\eta)$ 
consist of the weighted configurations 
where at least, respectively, strictly more than 
half of the total mass is concentrated in the point $\eta$.

Choose now $A\subset\DF X$ as the ``diagonal",
that is, as the compact antipodal subset of all configurations $(\eta)$ concentrated in one point.
Then the thickenings
$(\Th_t)_{F\ddot u}(A)$ and $(\ol{\Th}_t)_{F\ddot u}(A)$ of $A$ inside $\DF X$ 
equal the subsets of weighted configurations
which are {\em not semistable}, respectively, {\em not stable}. 
In the case when $\Th_t$ is balanced, both notions agree: 
``stable=semistable''.

The sets of stable and semistable configurations depend on the weights $t$. 
For instance, if $t_i>M/2$ for some $i$, 
then there are no semistable weighted configurations,
equivalently,
$$
(\Th_t)_{F\ddot u}(A)= \DF X.
$$
For instance, for $n=2$ and any $t_1\ne t_2$, there are no semistable configurations.   

In contrast, if $n\geq3$ and $t_i< M/2$ for all $i$, then there are always stable configurations,
for instance, the configurations where no two point coincide.
Equivalently,
$$
(\Th_t)_{F\ddot u}(A) \neq \DF X
$$
in this case. 
\end{example}

We return to the general discussion of thickenings in $\geo X$. 
Our motivation for introducing the notion of slimness is the following observation:
\begin{lem}[Disjointness of slim thickenings]
\label{lem:sldisj}
(i) Let $\taumod\subseteq\simod$ be an $\iota$-invariant face type,
and let $\Th\subset W$ be a slim $\Wt$-left invariant 
thickening.
Then for any two antipodal simplices $\tau,\hat\tau\in\Flagt$
it holds that 
$$\ThF(\tau)\cap\ThF(\hat\tau)=\emptyset .$$

(ii) More generally,
suppose that $\numod\subseteq\simod$ is another face type
and that the slim thickening $\Th$ is also $\Wn$-right invariant.
Then for any $G$-orbit $G\xi\subset\geo X$ of type $\bar\xi=\theta(\xi)\in\inte(\numod)$
and any two antipodal simplices $\tau,\hat\tau\in\Flagt$
it holds that 
$$\Th(\tau)\cap\Th(\hat\tau)\cap G\xi=\emptyset .$$
\end{lem}
\proof
Part (i) follows from Lemma~\ref{lem:triai} and the definition of slimness.
Indeed, suppose that $\ThF(\tau)\cap\ThF(\hat\tau)$ contains a chamber $\si$.
Then $\pos(\si,\tau),\pos(\si,\hat\tau)\in\Wt\backslash\Th$.
By Lemma~\ref{lem:triai},
$\pos(\si,\hat\tau)\succ_{\taumod}\cpos(\si,\tau)$.
Hence also 
$\cpos(\si,\tau)\in\Wt\backslash\Th$,
equivalently,
$\pos(\si,\tau)\in\Wt\backslash w_0\Th$.
It follows that $\Th\cap w_0\Th\neq\emptyset$, contradicting slimness.

Part (ii) follows because the thickenings are unions of stars of simplices of type $\numod$.
Indeed, by (i), the intersection $\Th(\tau)\cap\Th(\hat\tau)$ contains no chamber,
and hence it cannot contain the star of a simplex of type $\numod$. 
\qed

\section{Asymptotic geometric notions in symmetric spaces}

\subsection{Shadows at infinity and strong asymptoticity of Weyl cones}

For a simplex $\tau_-\subset\geo X$ of type $\iota\tau_{mod}$
and a point $x\in X$,
we consider the function 
\begin{equation}
\label{eq:distfrpar}
\tau\mapsto d(x,P(\tau_-,\tau))
\end{equation}
on the open Schubert stratum $C(\tau_-)\subset\Flagt$.
We denote by $\tau_+\in C(\tau_-)$ the chamber $x$-opposite to $\tau_-$.
\begin{lem}
\label{lem:contpr}
The function (\ref{eq:distfrpar}) is continuous and proper.
\end{lem}
\proof
This follows from the fact that 
$C(\tau_-)$ and $X$ are homogeneous spaces for the parabolic subgroup $P_{\tau_-}$.
Indeed, 
continuity follows from the continuity of the function 
$$g\mapsto d(x,P(\tau_-,g\tau_+))=d(g^{-1}x,P(\tau_-,\tau_+))$$
on $P_{\tau_-}$ which factors through the orbit map $P_{\tau_-}\to C(\tau_-),g\mapsto g\tau_+$.

Regarding properness,
note that a simplex $\tau\in C(\tau_-)$ is determined by any point $y$ contained in the parallel set $P(\tau_-,\tau)$,
namely as the simplex $y$-opposite to $\tau_-$.
Thus,
if $P(\tau_-,\tau)\cap B(x,R)\neq\emptyset$ for some fixed $R>0$,
then there exists $g\in P_{\tau_-}$ such that $\tau=g\tau_+$ and $d(x,gx)<R$. 
In particular, $g$ is bounded.
This implies properness.
\qed

\medskip
Moreover, 
the function (\ref{eq:distfrpar}) has a unique minimum zero in $\tau_+$. 

We define the following open subsets of $C(\tau_-)$ 
which can be regarded as {\em shadows} of balls in $X$ with respect to $\tau_-$. 
For $x\in X$ and $r>0$, we put
\begin{equation*}
U_{\tau_-,x,r}:=\{\tau\in C(\tau_-) | d(x,P(\tau_-,\tau))<r\} .
\end{equation*}
The next fact expresses the uniform strong asymptoticity
of asymptotic Weyl cones. 
\begin{lem}
\label{lem:expconvsect}
For $r,R>0$ exists $d=d(r,R)>0$ such that:

If $y\in V(x,\st(\tau_-))$ with $d(y,\D V(x,\st(\tau_-)))\geq d(r,R)$, 
then 
$U_{\tau_-,x,R}\subset U_{\tau_-,y,r}$.
\end{lem}
\proof
If $U_{\tau_-,x,R}\not\subset U_{\tau_-,y,r}$ 
then there exists $x'\in B(x,R)$ such that $d(y,V(x',\st(\tau_-)))\geq r$. 
Thus, if the assertion is wrong, 
there exist a sequence $x_n\to x_{\infty}$ in $\ol B(x,R)$ 
and a sequence $y_n\to\infty$ in $V(x,\st(\tau_-))$ 
such that $d(y_n,\D V(x,\st(\tau_-)))\to+\infty$
and $d(y_n,V(x_n,\st(\tau_-)))\geq r$.  

Let $\rho:[0,+\infty)\to V(x,\tau_-)$ be a geodesic ray with initial point $x$ and asymptotic to an interior point of $\tau_-$.
Then the sequence $(y_n)$ eventually enters every Weyl cone $V(\rho(t),\st(\tau_-))$. 
Since the distance function $d(\cdot,V(x_n,\st(\tau_-)))$ is convex and bounded, and hence non-increasing 
along rays asymptotic to $\st(\tau_-)$, 
we have that 
\begin{equation*}
R\geq d(x,V(x_n,\st(\tau_-))) 
\geq d(\rho(t),V(x_n,\st(\tau_-)))\geq d(y_n,V(x_n,\st(\tau_-)))\geq r
\end{equation*}
for $n\geq n(t)$. 
It follows that 
\begin{equation*}
R\geq d(\rho(t),V(x_{\infty},\st(\tau_-)))\geq r
\end{equation*}
for all $t$. 
However, the ray $\rho$ is strongly asymptotic to $V(x_{\infty},\st(\tau_-))$, a contradiction. 
\qed

\subsection{Asymptotic properties of sequences and subgroups}

We first consider sequences in the model euclidean Weyl chamber $\De$.
\begin{dfn}\label{def:pure_and_regular}
We say that a sequence $(\de_n)$ in $\De$ is 

(i) {\em $\taumod$-pure} if it is contained in a tubular neighborhood of the sector $V(0,\taumod)$
and drifts away from its boundary $\D V(0,\taumod)=V(0,\D\taumod)$,
$$ d(\de_n,V(0,\D\taumod)) \to+\infty .$$
(ii) {\em $\taumod$-regular} if
$$ d(\de_n,V(0,\D\st(\taumod))) \to+\infty .$$
\end{dfn}
These properties are stable under bounded perturbation of the sequence,
due to the triangle inequality
$|d_{\De}(x,y)-d_{\De}(x',y')|\leq d(x,x')+d(y,y')$.
Therefore the following definitions for sequences in $X$ and $G$ are sensible:
\begin{dfn}[Pure and weakly regular]
(i) We say that a sequence $(x_n)$ in $X$ is {\em $\taumod$-pure}, respectively, {\em $\taumod$-regular}
if for some (any) base point $o\in X$ the sequence of $\De$-distances $d_{\De}(o,x_n)$ 
has this property.

(ii) We say that a sequence $(g_n)$ in $G$ is {\em $\taumod$-pure}, respectively, {\em $\taumod$-regular}
if for some (any) point $x\in X$ the orbit sequence $(g_nx)$ in $X$ has this property.

(iii) We say that a subgroup $\Ga<G$ is {\em $\taumod$-regular}
if all sequences of pairwise distinct elements in $\Ga$ have this property.
\end{dfn}
The face type of a pure sequence is uniquely determined.
Moreover,
a $\taumod$-regular sequence is $\taumod'$-regular for every face type $\taumod'\subseteq\taumod$,
because $\ost(\taumod')\supseteq\ost(\taumod)$.

Note that $\taumod$-regular subgroups are in particular discrete.
If $\rank(X)=1$, then discreteness is equivalent to ($\simod$-)regularity.
In higher rank, regularity can be considered as a strengthening of discreteness:
A discrete subgroup $\Ga<G$ may not be $\taumod$-regular for any face type $\taumod$;
this can happen e.g.\ for free abelian subgroups of transvections of rank $\geq2$.

We observe furthermore:
\begin{lem}
\label{lem:obspureg}
(i) $\taumod$-pure sequences are $\taumod$-regular.

(ii) Every sequence, which diverges to infinity, 
contains a $\taumod$-pure subsequence for some face type $\taumod\subseteq\simod$.
\end{lem}
\proof
Assertion (i) is a direct consequence of the definitions,
and (ii) follows by induction on face types.
\qed

\medskip
Note also that a sequence, which diverges to infinity, 
is $\taumod$-regular if and only if it contains $\numod$-pure subsequences only for face types $\numod\supseteq\taumod$. 
(We will not use this fact.)

\begin{rem}
[Relation to Finsler compactifications]
\label{rem:regfinscomp}
There is a close relation between the regularity of sequences and the asymptotic geometry of certain $G$-invariant 
{\em Finsler metrics} on $X$, see \cite[\S 8.1.2]{bordif}.
For instance,
a sequence in $X$ is ($\simod$-)regular if and only if it accumulates 
at the Furstenberg boundary inside the regular Finsler compactification.
\end{rem}

\section{Some topological dynamics}

\subsection{(Proper) discontinuity and dynamical relation}
\label{sec:propdynrel}

Let $Z$ be a compact metrizable space,
and let $\Ga\subset\Homeo(Z)$ be a countably infinite subgroup.
We consider the action $\Ga\acts Z$.
\begin{dfn}[Discontinuous]
\label{dfn:wand}
A point $z\in Z$ is called {\em wandering} 
with respect to the $\Ga$-action 
if the action is {\em discontinuous} at $z$,
i.e.\ if $z$ has a neighborhood $U$ such that 
$U\cap\ga U\neq\emptyset$ for at most finitely many $\ga\in\Ga$.
\end{dfn}
Nonwandering points are called {\em recurrent}. 

\begin{dfn}[Domain of discontinuity]
\label{dfn:domdisc}
We call the set 
$$\Om_{disc}\subset Z$$ 
of wandering points 
the {\em wandering set} or {\em domain of discontinuity} for the action $\Ga\acts Z$.
\end{dfn}
Note that $\Om_{disc}$ is open and $\Ga$-invariant.

\begin{dfn}[Proper]
The action of $\Ga$ on an open 
subset $U\subset Z$
is called {\em proper} if
for every compact subset $K\subset U$ it holds that 
$K\cap\ga K\neq\emptyset$ for at most finitely many $\ga\in\Ga$.
\end{dfn}
In particular, 
the action of $\Ga$ on $U$ is then discontinuous, $U\subseteq\Om_{disc}$,
and is therefore called {\em properly discontinuous}. 
\begin{dfn}[Domain of proper discontinuity]
\label{dfn:dompropdisc}
If $\Ga$ is a group, 
we call a $\Ga$-invariant open subset $\Om\subseteq\Om_{disc}$ on which $\Ga$ acts properly
a {\em domain of proper discontinuity} for $\Ga$.
\end{dfn}
The orbit space $\Om/\Ga$ is then Hausdorff. 
Note that in general there is {\em no unique} maximal {\em proper} domain of discontinuity. 

\begin{example}[Nonunique maximal domain of proper discontinuity]
\label{ex:nonundd}
Consider the infinite cyclic group $\Ga\cong\Z$ acting projectively on $Z=\R P^2$, so that a generator $\ga$ of $\Ga$ acts as the projectivization of a diagonal matrix with distinct positive eigenvalues $\la_1>\la_2>\la_3$. Let $e_1, e_2, e_3\in Z$ be the three fixed points of $\ga$ (eigenspaces for $\la_1, \la_2, \la_3$ respectively). 
Let $E_{ij}\subset Z$ 
denote the projective lines spanned by $e_i$ and $e_j$, $i< j$. 
Then $\Om_{disc}=Z-\{e_1,e_2,e_3\}$,
and both sets
$U_1=Z-(E_{23} \cup\{e_1\})$ and $U_3=Z-(E_{12}\cup\{e_3\})$ 
are maximal domains of proper discontinuity for $\Gamma$. 
(The maximality follows from the fact that the points on $E_{12}$ 
are dynamically related to the points on $E_{23}$.) 
Observe also that in this example both $U_1/\Ga$ and $U_3/\Ga$ are compact. 
\end{example}

Discontinuity and proper discontinuity can be nicely expressed 
using the notion of dynamical relation.
The following definition is due to Frances \cite[Def.\ 1]{Frances}:
\begin{dfn}[Dynamically related]
Two points $z,z'\in Z$ 
are called 
{\em dynamically related} with respect to a sequence $(h_n)$ in $\Homeo(Z)$,
$$ z\stackrel{(h_n)}{\sim}z' $$
if there exists a sequence $z_n\to z$ in $Z$ such that $h_nz_n\to z'$. 

The points $z,z'$ are called 
{\em dynamically related} with respect to the $\Ga$-action, 
$$ z\stackrel{\Ga}{\sim}z' $$
if there exists a sequence $\ga_n\to\infty$ in $\Ga$ such that $z\stackrel{(\ga_n)}{\sim}z'$.
\end{dfn}
Here, for a sequence $(\ga_n)$ in $\Ga$ we write $\ga_n\to\infty$
if every element of $\Ga$ occurs at most finitely many times in the sequence. 

One verifies (see e.g. \cite{manicures}):

(i) Dynamical relation is a closed relation in $Z\times Z$. 

(ii) Points in different $\Ga$-orbits are dynamically related 
if and only if their orbits cannot be separated 
by disjoint $\Ga$-invariant open subsets. 

The concept of dynamical relation is useful for our discussion of discontinuity, because:

(i) A point is nonwandering if and only if it is dynamically related to itself. 

(ii) The action is proper on an open subset $U\subset Z$ 
if and only if no two points in $U$ are dynamically related.

\subsection{Accumulation and proper discontinuity}
\label{sec:accac}

In this paper,
we derive proper discontinuity of actions from a certain accumulation behavior
which is a relaxation of {\em convergence dynamics}.

Let $Z$ and $\Ga$ be as above.
Let $(Y_n)$ be a sequence of subsets of $Z$.
We denote by $\Acc((Y_n))\subset Z$ the closed subset consisting of the accumulation points of all 
sequences $(y_n)$ of points $y_n\in Y_n$.
\begin{dfn}[Accumulation]
We say that the sequence of subsets $Y_n\subset Z$ 
{\em accumulates at} a subset $S\subset Z$,
$$ Y_n\acc S ,$$
if $\Acc((Y_n))\subseteq S$. 
\end{dfn}
If $S\subset Z$ is closed,
then the sequence $(Y_n)$ accumulates at $S$ 
if and only if every neighborhood $U$ of $S$ contains all but finitely many of the subsets $Y_n$.

\medskip
We first consider the dynamics of a sequence $(h_n)$ in $\Homeo(Z)$.
\begin{dfn}[Accumulating sequence]
\label{def:accseq}
For compact subsets $A_{\pm}\subset Z$
we say that the sequence $(h_n)$ {\em accumulates} at $A_+$ outside $A_-$, 
briefly, {\em $(A_-,A_+)$-accumulates}, if
\begin{equation}
\label{eq:acc}
h_nK \acc A_+
\end{equation} 
for all compacta $K$ disjoint from $A_-$.
\end{dfn}
Property (\ref{eq:acc}) is a statement about the locally uniform accumulation of the $(h_n)$-orbits 
initiating outside the exceptional subset $A_-$
and can be rephrased in terms of dynamical relations between points in $Z$ 
with respect to the $(h_n)$-action.
Namely, equivalently,
for all points $z,z'\in Z$ it holds that:
\begin{equation}
\label{eq:accdynr}
z\stackrel{(h_{n_k})}{\sim}z'\hbox{ for some subsequence $(h_{n_k})$ }
\quad\Ra\quad
z\in A_-\hbox{ or }z'\in A_+
\end{equation} 
Note that the dynamical relation condition 
$z\stackrel{(h_{n_k})}{\sim}z'$
is equivalent to the dual condition
$z'\stackrel{(h_{n_k}^{-1})}{\sim}z$,
and consequently we have the symmetry:
$$ \hbox{$(h_n)$ is $(A_-,A_+)$-accumulating} 
\quad\Leftrightarrow\quad 
\hbox{$(h_n^{-1})$ is $(A_+,A_-)$-accumulating} $$
Note that if $A_{\pm}\subset A'_{\pm}$,
then $(A_-,A_+)$-accumulation implies $(A'_-,A'_+)$-accumulation.

\medskip
Now we consider the action $\Ga\acts Z$.
\begin{dfn}[Accumulating action I]
We say that the action $\Ga\acts Z$ is {\em $(A_-,A_+)$-ac\-cu\-mu\-la\-ting}
if every sequence $\ga_n\to\infty$ in $\Ga$ has an $(A_-,A_+)$-accumulating subsequence. 
\end{dfn}
According to (\ref{eq:accdynr}) we obtain for dynamical relations:
\begin{lem}[Dynamical relations I]
\label{lem:dynrel1}
If the action $\Ga\acts Z$ is $(A_-,A_+)$-accumulating, then for any two points $z,z'\in Z$ it holds that:
\begin{equation}
\label{eq:accdynrgp}
z\stackrel{\Ga}{\sim}z'
\quad\Ra\quad
z\in A_-\hbox{ or }z'\in A_+
\end{equation} 
\end{lem}
\proof
This is a direct consequence of (\ref{eq:accdynr}).
\qed

\medskip
We conclude:
\begin{prop}[Proper discontinuity I]
\label{prop:pd1}
If the subsets $A_{\pm}$ are $\Ga$-invariant
and if the action $\Ga\acts Z$ is $(A_-,A_+)$-accumulating, then the action 
$$ \Ga\acts Z-(A_-\cup A_+) $$
is properly discontinuous.
\end{prop}
\proof
By the lemma, there are no dynamical relations between points outside $A_-\cup A_+$.
\qed

\medskip
Suppose that ${\mathcal A}_{\pm}$ are $\Ga$-invariant compact (with respect to the Hausdorff topology) 
families of compact subsets $A_{\pm}\subset Z$.
\begin{dfn}[Limit family]
The {\em forward limit family} of $\Ga$ 
with respect to $({\mathcal A}_-,{\mathcal A}_+)$
is the family ${\mathcal L}_+\subset{\mathcal A}_+$ consisting of all subsets $A_+\in{\mathcal A}_+$
for which there exists a sequence $\ga_n\to\infty$ in $\Ga$ 
which is $(A_-,A_+)$-accumulating for some subset $A_-\in{\mathcal A}_-$.
Similarly, we define the {\em backward limit family} ${\mathcal L}_-\subset{\mathcal A}_-$.
\end{dfn}
The limit families ${\mathcal L}_{\pm}$ are $\Ga$-invariant. 
Due to the compactness of the families ${\mathcal A}_{\pm}$, they are closed and hence compact themselves:
\begin{lem}
\label{lem:limfcpt}
${\mathcal L}_{\pm}$ is closed.
\end{lem}
\proof
Suppose, for instance, that $(A^k_+)$ is a sequence in ${\mathcal L}_+$ 
such that $A^k_+\to A_+\in{\mathcal A}_+$.
There exist sequences $\ga^k_n\to\infty$ in $\Ga$ 
which $(A^k_-,A^k_+)$-accumulate for some $A^k_-\in{\mathcal A}_-$
(in fact $\in{\mathcal L}_-$).
After passing to a subsequence, we may assume that $A^k_-\to A_-\in{\mathcal A}_-$.
A diagonal argument yields an $(A_-,A_+)$-accumulating sequence 
$\left(\ga^{k(m)}_{n(m)}\right)_m$ in $\Ga$.
Hence $A_+\in{\mathcal L}_+$ and ${\mathcal L}_+$ is closed.
\qed

\medskip
As a consequence of the lemma, the $\Ga$-invariant subsets
$$ T_{\pm} := \bigcup_{A_{\pm}\in{\mathcal L}_{\pm}}A_{\pm} \subset Z $$
are compact.

\begin{dfn}[Accumulating action II]
We say that the action $\Ga\acts Z$ is {\em $({\mathcal A}_-,{\mathcal A}_+)$-ac\-cu\-mu\-la\-ting}
if every sequence $\ga_n\to\infty$ in $\Ga$ has a subsequence which is $(A_-,A_+)$-accumulating
for some subsets $A_{\pm}\in{\mathcal A}_{\pm}$.
\end{dfn}
For such accumulating actions, the limit families are closely related to their dynamics:
\begin{lem}[Dynamical relations II]
\label{lem:dynrel2}
If the action $\Ga\acts Z$ is $({\mathcal A}_-,{\mathcal A}_+)$-accumulating, 
then for any two points $z,z'\in Z$ it holds that:
\begin{equation*}
z\stackrel{\Ga}{\sim}z'
\quad\Ra\quad
z\in T_-\hbox{ or }z'\in T_+
\end{equation*} 
\end{lem}
\proof
If two points are dynamically related with respect to the $\Ga$-action,
then they are dynamically related with respect to an $(A_-,A_+)$-accumulating sequence in $\Ga$
with $A_{\pm}\in{\mathcal L}_{\pm}$ and, hence, $A_{\pm}\subset T_{\pm}$.
The assertion therefore follows from (\ref{eq:accdynr}).
\qed

\medskip
We conclude as before:
\begin{prop}[Proper discontinuity II]
\label{prop:pd2}
If the action $\Ga\acts Z$ is $({\mathcal A}_-,{\mathcal A}_+)$-accumulating, 
then the action 
$$ \Ga\acts Z-(T_-\cup T_+) $$
is properly discontinuous.
\end{prop}

\begin{rem}[Convergence actions]
\label{rem:convact}
The action $\Ga\acts Z$ is a {\em convergence action}
(see e.g. \cite{Bowditch_config})
if and only if it is $({\mathcal A}_-,{\mathcal A}_+)$-accumulating
with ${\mathcal A}_{\pm}$ the family of one point subsets.
The limit families ${\mathcal L}_{\pm}$ then become the {\em limit set} $\La\subset Z$ of the action. 
The action on its complement is properly discontinuous, 
compare Proposition~\ref{prop:pd2}.
We recover the {\em dynamical decomposition}
$$ Z = \Om_{disc}\sqcup \La $$
and that the action on the domain of discontinuity is proper.
Hence, for convergence actions there exists a unique maximal domain of proper discontinuity.

The main example of convergence actions comes from the following fact: 
Every discrete group $\Ga$ of isometries of a proper Gromov hyperbolic geodesic metric space $Y$ 
acts as a convergence group 
on the visual compactification $\bar Y=Y\cup\geo Y$,
and in particular on the Gromov boundary $\geo Y$ of $Y$. 
\end{rem}

\begin{rem}[Accumulation phenomena in nonpositive curvature]
Convergence type behavior in the sense of accumulation 
has been studied by Karlsson, Papasoglu and Swenson 
in the general context of nonpositive curvature.
They showed that for proper isometric actions $\Ga\acts Y$ on CAT(0) spaces
the induced action $\Ga\acts\geo Y$ on the visual boundary is
$({\mathcal B}(\theta),{\mathcal B}(\pi-\theta))$-accumulating for $0<\theta<\pi$,
where ${\mathcal B}(\theta)$ is the family of closed balls of Tits radius $\theta$ in $\geo Y$,
see \cite[Thm.\ 1]{Karlsson} and \cite[Thm.\ 4]{PapaSwen}.
Some of our results can be viewed as combinatorial versions of this (Tits) metric result
for actions on CAT(0) {\em model} spaces of higher rank,
see e.g.\ Corollary~\ref{cor:orbacc} and Lemma~\ref{lem:discracc} below.
\end{rem}

\subsection{Expansion and cocompactness}
\label{sec:trexpand}

In this section,
let $(Z,d)$ be a compact metric space
and let $\Ga \acts Z$ be a continuous action of a discrete group.

The following notion 
is due to Sullivan \cite[\S 9]{Sullivan}:
\begin{dfn}[Expanding action]
\label{dfn:expact}
We say that the action $\Ga\acts Z$ is {\em expanding} 
at the {\em point} $z\in Z$ 
if there exists an element $\ga\in\Ga$ 
which is {\em uniformly expanding} on a neighborhood $U$ of $z$, 
i.e.\ for some constant $c>1$ and all points $z_1,z_2\in U$ we have  
\begin{equation*}
d(\ga z_1,\ga z_2)\geq c\cdot d(z_1,z_2) .
\end{equation*}
We say that the action of $\Ga$ is {\em expanding} 
at a compact $\Ga$-invariant {\em subset} $E\subset Z$ 
if it is expanding at all points $z\in E$. 
\end{dfn}

\begin{rem}
\label{rem:arbstrexp}
If the action $\Ga\acts Z$ is expanding at $E$, then it is {\em arbitrarily strongly} expanding there,
i.e.\ for every point $z\in E$ exist a sequence $(\ga_n)$ in $\Ga$ and a sequence of (shrinking) 
neighborhoods $U_n$ of $z$ 
such that the $\ga_n|_{U_n}$ are uniformly expanding with expansion factors 
$c_n\to+\infty$.
This follows directly from the definition by iterating locally expanding elements.
Note that, as a consequence, the action is expanding at $E$ also with respect to any bilipschitz equivalent metric on $Z$.
\end{rem}
We will need the following more general notion of {\em partial expansion}. 
We suppose that the action $\Ga\acts Z$ 
has the following structure: 
There is a $\Ga$-invariant compact subset $E\subset Z$ 
and a continuous map $\pi: E\to\La$ 
onto a compact topological space $\La$
(e.g.\ a fiber bundle), 
such that the restricted action $\Ga\acts E$ is fiber preserving, 
i.e.\ it descends to a continuous action $\Ga\acts\La$. 
We set $E_{\la}:=\pi^{-1}(\la)$. 

\begin{dfn}[Transversely expanding action I]
We say that the action $\Ga\acts Z$ is 
{\em expanding transversely to $\pi$ at the fiber $E_{\la}$}
if there exist an element $\ga\in\Ga$ 
and a neighborhood $U\subset Z$ of $E_{\la}$ 
such that for some constant $c>1$ we have  
\begin{equation}
\label{ineq:trexpan}
d(\ga z,E_{\ga\la'})\geq c\cdot d(z,E_{\la'})
\end{equation}
for all points $z\in U$ and fibers $E_{\la'}\subset U$. 

We say that the action $\Ga\acts Z$ is 
{\em expanding at $E$ transversely to $\pi$} 
if it is expanding at all fibers $E_{\la}$. 
\end{dfn}
The action $\Ga\acts Z$ is expanding at $E$
if and only if it is expanding at $E$ transversely to $\id_E$.

The concept of expansion is important to us 
due to the following observation: 
\begin{prop}[Transversely expanding implies cocompact on the complement I]
\label{prop:trexpcoco}
If the action $\Ga\acts Z$ is expanding at $E$ transversely to $\pi$, 
then the action $\Ga\acts Z-E$ is cocompact. 
\end{prop}
\proof
We claim that for some constant $c>1$,  
\begin{equation}
\label{eq:maxprinctr}
\sup d(\cdot,E)|_{\Ga z} > c\cdot d(z,E)
\end{equation}
for all $z\in Z-E$ sufficiently close to $E$. 
Otherwise, 
there would exist a sequence $(z_n)$ in $Z-E$ accumulating at $E$ 
and a sequence of constants $c_n\to1$ such that
\begin{equation*}
d(\ga z_n,E)\leq c_n\cdot d(z_n,E)
\end{equation*}
for all $n\in\N$ and $\ga\in\Ga$. 
Since $E$ is compact, 
we may assume, after passing to a subsequence, 
that $(z_n)$ accumulates at a fiber, 
$z_n\to E_{\la}$. 
Due to expansion, 
there exists an element $\ga_{\la}\in\Ga$ 
which satisfies the expansion property (\ref{ineq:trexpan})
on a neighborhood $U_{\la}\subset Z$ of $E_{\la}$ 
with some expansion factor $c_{\la}>1$. 
Let $E_{\ga_{\la}\la_n}$ be the fiber closest to $\ga_{\la}z_n$, 
$d(\ga_{\la}z_n,E_{\ga_{\la}\la_n})=d(\ga_{\la}z_n,E)$. 
Then $\la_n\to\la$. 
Since $z_n\in U_{\la}$ and $E_{\la_n}\subset U_{\la}$ for large $n$, 
it follows that 
\begin{equation*}
c_{\la}\cdot d(z_n,E) \leq c_{\la}\cdot d(z_n,E_{\la_n})
\leq d(\ga_{\la}z_n,E_{\ga_{\la}\la_n})
= d(\ga_{\la}z_n,E)\leq c_n\cdot d(z_n,E) ,
\end{equation*}
a contradiction confirming our claim. 

Let $U\subset Z$ be an open tubular 
neighborhood of $E$ 
where (\ref{eq:maxprinctr}) holds. 
Thus, no $\Ga$-orbit is entirely contained in 
$U-E$ and, therefore, every $\Ga$-orbit in $Z-E$ 
meets the compact subset $Z-U\subset Z-E$. 
\qed

\medskip
The above argument (from \cite[sec.\ 2.2]{coco13}) leads actually to a more general result.

Let us suppose, more generally, that the action $\Ga\acts Z$ 
has the following structure: 
There is a $\Ga$-invariant compact subset $E\subset Z$ 
which is represented as the (not necessarily disjoint) union 
of a $\Ga$-invariant collection ${\mathcal E}=\{E_{\la}:\la\in\La\}$
of compact subsets $E_{\la}\subset Z$
parametrized by some set $\La$.

\begin{dfn}[Transversely expanding action II]
We say that the action $\Ga\acts Z$ is 
{\em expanding transversely to ${\mathcal E}$ at a point $z$} 
if there exist an element $\ga\in\Ga$, a neighborhood $U\subset Z$ of $z$ 
and a constant $c>1$ such that we have  
\begin{equation}
\label{ineq:trexpan2}
d(\ga u, \ga E_{\la})\geq c\cdot d(u,E_{\la})
\end{equation}
for all points $u\in U-E$ and all $E_{\la}$ which have nonempty intersection with $U$. 
We say that the action $\Ga\acts Z$ is 
{\em expanding transversely to ${\mathcal E}$} 
if it is expanding at all points $z\in E$. 
\end{dfn}

\begin{prop}[Transversely expanding implies cocompact on the complement II]
\label{prop:trexpcoco1}
If the action $\Ga\acts Z$ is expanding transversely to ${\mathcal E}$,
then the action $\Ga\acts Z-E$ is cocompact.
\end{prop}
\proof
We claim that for some constant $c>1$,  
\begin{equation}
\label{eq:maxprinctr1}
\sup d(\cdot,E)|_{\Ga u} > c\cdot d(u,E)
\end{equation}
for all $u\in Z-E$ sufficiently close to $E$. 
Otherwise, there would exist a sequence $(u_n)$ in $Z-E$ accumulating at $E$ 
and a sequence of constants $c_n\to1$ such that
\begin{equation*}
d(\ga u_n,E)\leq c_n\cdot d(u_n,E)
\end{equation*}
for all $n\in\N$ and $\ga\in\Ga$. 
Since $E$ is compact, 
we may assume, after passing to a subsequence, 
that $(u_n)$ converges to some point $z\in E_{\la}$ for some $\la\in\La$. 
Due to expansion, 
there exists an element $\ga_z\in\Ga$ 
which satisfies the expansion property (\ref{ineq:trexpan2})
on a neighborhood $U\subset Z$ of $z$ 
with some expansion factor $c > 1$. 
Let $\ga_z E_{\la_n}\in{\mathcal E}$ be the set in the collection ${\mathcal E}$ closest to $\ga_z u_n$, 
$$d(\ga_z u_n, E)=d(\ga_z u_n, \ga_z E_{\la_n})=d(\ga_z u_n, \ga_z z_n)$$
with $z_n\in E_{\la_n}$.
Then $z_n\to z$ because $\ga_z u_n\to\ga_z z\in E$,
which implies that $E_{\la_n}\cap U\ne \emptyset$ for all sufficiently large $n$. 
It follows that (for large $n$), 
\begin{equation*}
c\cdot d(u_n,E) \leq c\cdot d(u_n,E_{\la_n})
\leq d(\ga_z u_n, \ga_z E_{\la_n})= d(\ga_z u_n, E)\leq c_n\cdot d(u_n, E) ,
\end{equation*}
a contradiction confirming our claim. 

Let $U\subset Z$ be an open tubular neighborhood of $E$ 
where (\ref{eq:maxprinctr1}) holds. 
Thus, no $\Ga$-orbit is entirely contained in 
$U-E$ and, therefore, every $\Ga$-orbit in $Z-E$ 
meets the compact subset $Z-U\subset Z-E$. 
\qed

\section{Accumulation dynamics on flag manifolds and proper discontinuity}
\label{sec:accfl}

We now study the dynamics of $G$ and its discrete subgroups $\Ga<G$ 
on its associated flag manifolds, equivalently, 
on (the $G$-orbits in) the visual boundary $\geo X$.
In this section, we will discuss a certain dynamical behavior, which is a relaxed version of convergence dynamics,
and use it to construct domains of proper discontinuity for discrete subgroups.

\subsection{Weakly contracting sequences}
\label{sec:wcseq}

Let $(g_n)$ be a sequence in $G$,
and let $\taumod\subseteq\simod$ be a face type.

We consider the following contraction property for the dynamics of $(g_n)$ on $\Flagt$.
An equivalent notion had been studied in \cite{Benoist}, see \S 3.5 there.
\begin{dfn}[$\taumod$-Contracting sequence]\label{def:contracting_sequence}
We say that the sequence $(g_n)$ is {\em $\taumod$-con\-trac\-ting} 
if there exist simplices $\tau_{\pm}$ of type $\pmtaumod$ such that 
\begin{equation}
\label{eq:contrtau}
g_n|_{C(\tau_-)}\to\tau_+
\end{equation} 
uniformly on compacta as $n\to+\infty$.
\end{dfn}
We recall that $C(\tau_-)$ is a dense open subset of $\Flagt$.

Property (\ref{eq:contrtau}) means that $(g_n)$ is $(\Flagt-C(\tau_-),\tau_+)$-accumulating,
cf.\ Definition~\ref{def:accseq}.
It can be rephrased in terms of dynamical relations between points in $\Flagt$
with respect to the $(g_n)$-action.
Namely, equivalently,
for all simplices $\tau,\tau'\in\Flagt$ it holds that, cf.\ (\ref{eq:accdynr}):
\begin{equation}
\label{eq:contrtaudynr}
\tau\stackrel{(g_{n_k})}{\sim}\tau'\hbox{ for some subsequence $(g_{n_k})$ }
\quad\Ra\quad
\tau\not\in C(\tau_-)\hbox{ or }\tau'=\tau_+
\end{equation} 

The conclusion of the last implication can be expressed in terms of relative positions:
$$ \pos(\tau,\tau_-)\hbox{ maximal }
\quad\Ra\quad
\pos(\tau',\tau_+)\hbox{ minimal } $$
We observe that the last implication follows from the combinatorial inequality
\begin{equation}
\label{ineq:dynrelsimp}
\pos(\tau',\tau_+)\prec\cpos(\tau,\tau_-) .
\end{equation}
The next result shows that this inequality holds for dynamically related points 
on all flag manifolds, thought of as $G$-orbits in $\geo X$.
It is the key step in our study of proper discontinuity.
\begin{prop}[Dynamical relation inequality]
\label{prop:dynreltaun}
The following are equivalent:

(i) Property (\ref{eq:contrtau})

(ii) For any two points $\xi,\xi'\in\geo X$ it holds that:
\begin{equation}
\label{eq:dynrelposcondtaun}
\xi\stackrel{(g_{n_k})}{\sim}\xi'\hbox{ for some subsequence $(g_{n_k})$ }
\quad\Ra\quad
\pos(\xi',\tau_+)\prec\cpos(\xi,\tau_-) 
\end{equation}
\end{prop}
\proof
Suppose first that property (\ref{eq:contrtau}) holds and that $\xi\stackrel{(g_n)}{\sim}\xi'$.
Then $\xi$ and $\xi'$ lie in the same $G$-orbit,
$G\xi=G\xi'$, 
and there exists a sequence $(\xi_n)$ in this $G$-orbit such that 
$\xi_n\to\xi$ and $g_n\xi_n\to\xi'$. 
Let $a\subset\geo X$ be an apartment containing $\tau_-$ and $\xi$. 
Nearby apartments $a_n$ containing $\xi_n$ 
can be obtained by using isometries $h_n\to e$ in $G$
with $\xi_n=h_n\xi$
and putting $a_n=h_n a$.
Let $\hat\tau_-\subset a$ be the simplex opposite to $\tau_-$,
and let $\tau_n=h_n\hat\tau_-\subset a_n$.
Then $\tau_n\to\hat\tau_-$.
Since $\hat\tau_-\in C(\tau_-)$,
the locally uniform convergence (\ref{eq:contrtau}) 
implies that $g_n\tau_n\to\tau_+$. 
We obtain
\begin{equation*}
\pos(\xi',\tau_+)\prec\pos(g_n\xi_n,g_n\tau_n)
=\pos(\xi_n,\tau_n)
=\pos(h_n\xi,h_n\hat\tau_-)
=\pos(\xi,\hat\tau_-)=\cpos(\xi,\tau_-)
\end{equation*}
where the first inequality follows from 
the semicontinuity of relative position
(Lemma~\ref{lem:semcontrelpos}). 

Conversely, suppose that (ii) holds.
Since inequality (\ref{ineq:dynrelsimp}) is a special case of the inequality in the implication of (\ref{eq:dynrelposcondtaun}),
it follows that (\ref{eq:contrtaudynr}) holds, equivalently, (\ref{eq:contrtau}).
\qed

\medskip
We observe a symmetry:
Condition (\ref{eq:dynrelposcondtaun}) is equivalent to the dual condition
\begin{equation}
\label{eq:dynrelposcondtaundual}
\xi'\stackrel{(g_{n_k}^{-1})}{\sim}\xi\hbox{ for some subsequence $(g_{n_k}^{-1})$ of $(g_n^{-1})$ }
\quad\Ra\quad
\pos(\xi,\tau_-)\prec\cpos(\xi',\tau_+)
\end{equation}
because both dynamical relation hypotheses are equivalent, as are the combinatorial inequality conclusions. 
Therefore the proposition implies that (\ref{eq:contrtau}) is equivalent to the dual property on $\Flagmt$ that 
\begin{equation*}
g_n^{-1}|_{C(\tau_+)}\to\tau_-
\end{equation*} 
uniformly on compacta as $n\to+\infty$.

Note that the simplices $\tau_{\pm}$ in (\ref{eq:contrtau}) are well-defined,
because this is clear for $\tau_+$ and follows for $\tau_-$ by symmetry. 

\medskip
Inequality (\ref{eq:dynrelposcondtaun}) can be (re)converted 
into a statement about the asymptotic behavior of arbitrary $(g_n)$-orbits in $\geo X$.
We can in general not expect that these orbits converge, but we obtain information where they accumulate.
For individual orbits, 
it follows that for a point $\xi\in\geo X$ the orbit $(g_n\xi)$ accumulates in $G\xi\subset\geo X$ 
at the Schubert cycle
\begin{equation*}
\{ \pos(\cdot,\tau_+) \prec\cpos(\xi,\tau_-)  \}
\end{equation*}
A locally uniform statement can be conveniently formulated using the language of thickenings:
\begin{cor}[Orbit accumulation]
\label{cor:orbacc}
If property (\ref{eq:contrtau}) holds, 
and if $\Th\subset W$ is a $\Wt$-left invariant thickening,
then the sequence $(g_n)$ is $(\Th^c(\tau_-),\Th(\tau_+))$-accumulating
(cf.\ Def.~\ref{def:accseq}).
\end{cor}
\proof
Otherwise, there is a dynamical relation
$\xi\stackrel{(g_{n_k})}{\sim}\xi'$ 
with $\xi\not\in\Th^c(\tau_-)$ and $\xi'\not\in\Th(\tau_+)$,
compare (\ref{eq:accdynr}),
i.e.\ $\pos(\xi,\tau_-)\not\in\Th^c$ and $\pos(\xi',\tau_+)\not\in\Th$.
Moreover (\ref{eq:contrtau}) implies (\ref{eq:dynrelposcondtaun}),
and hence the inequality 
$\pos(\xi',\tau_+)\prec\cpos(\xi,\tau_-)$. 
It follows that 
$$ \pos(\xi',\tau_+) \prec \cpos(\xi,\tau_-) = w_0 \pos(\xi,\tau_-)\in\Th $$
and hence $\pos(\xi',\tau_+)\in\Th$,
a contradiction. 
\qed

\subsection{Weak convergence subgroups}

Let $\Ga<G$ be a discrete subgroup.
\begin{dfn}[$\taumod$-Limit set]
\label{dfn:lims}
We define the {\em forward/backward $\taumod$-limit set} of $\Ga$
as the set 
$$\Lat^{\pm}=\Lat^{\pm}(\Ga)\subset\Flagpmt$$ 
of all simplices $\tau_{\pm}$ as in (\ref{eq:contrtau}) 
for all $\taumod$-contracting sequences $\ga_n\to\infty$ in $\Ga$.
\end{dfn}
Note that passing to a finite index subgroup does not change the limit sets. 

The limit sets $\Lat^{\pm}$ are $\Ga$-invariant and compact, cf.\ Lemma~\ref{lem:limfcpt}.
Moreover, one has the symmetry
$$ \Lamt^{\pm}(\Ga) = \Lat^{\mp}(\Ga) .$$
In particular,
if $\taumod$ is $\iota$-invariant
we can define the {\em $\taumod$-limit set} 
$$\Lat(\Ga):=\Lat^{\pm}(\Ga) .$$
To any $\Wt$-left invariant thickening $\Th\subset W$,
we associate the $\Ga$-invariant compact families of compact subsets 
$$ {\mathcal A}^-_{\taumod,Th} := \{ \Th^c(\tau_-) : \tau_-\in\Lat^- \} 
\quad\hbox{ and }\quad
{\mathcal A}^+_{\taumod,Th} := \{ \Th(\tau_+) : \tau_+\in\Lat^+ \} $$
The structure of the dynamics of the action $\Ga\acts\geo X$
is closely related to the limit sets 
if it enjoys contraction behavior in the following sense:
\begin{dfn}[$\taumod$-Convergence action]
\label{dfn:convact}
The action $\Ga\acts\geo X$ is called a {\em $\taumod$-con\-ver\-gence action} 
if every sequence $\ga_n\to\infty$ in $\Ga$ 
has a $\taumod$-contracting subsequence. 
The subgroup $\Ga<G$ is then called a {\em $\taumod$-convergence subgroup}.
\end{dfn}
\begin{rem}[Rank one]
If $\rank(X)=1$,
this property is equivalent to being a {\em convergence action} 
and is satisfied for all discrete subgroups $\Ga<G$,
compare Remark~\ref{rem:convact}.
\end{rem}

Corollary~\ref{cor:orbacc} implies:
\begin{prop}[Weak contraction implies accumulation]
If $\Ga<G$ is a $\taumod$-con\-ver\-gen\-ce subgroup
and if $\Th\subset W$ is a $\Wt$-left invariant thickening,
then the action $\Ga\acts\geo X$ is $({\mathcal A}^-_{\taumod,Th},{\mathcal A}^+_{\taumod,Th})$-accumulating.
\end{prop}
We obtain our main result for proper discontinuity:
\begin{thm}[Domains of proper discontinuity for $\taumod$-convergence subgroups]
\label{thm:pdwconv}
If $\Ga<G$ is a $\taumod$-con\-ver\-gen\-ce subgroup,
then for any $\Wt$-left invariant thickening $\Th\subset W$
the action 
$$ \Ga\acts\geo X-(\Th^c(\Latm)\cup\Th(\Latp))$$
is properly discontinuous. 
In particular,
if $\taumod$ is $\iota$-invariant and $\Th$ is fat,
then the action 
$$ \Ga\acts\geo X-\Th(\Lat) $$
is properly discontinuous. 
\end{thm}
\proof
The first assertion follows from the last proposition by applying 
Proposition~\ref{prop:pd2} with ${\mathcal L}_{\pm}={\mathcal A}^{\pm}_{\taumod,Th}$.
The second assertion follows because $\Th^c\subseteq\Th$
due to fatness. 
\qed

\medskip
Note that the thickenings of limit sets $\Th(\Lat^{\pm}(\Ga))$ are $\Ga$-invariant and compact. 

For examples of thickenings, 
we refer to section~\ref{sec:thickinf}.

\subsection{Weakly regular subgroups}
\label{sec:wreg}

The properties of contraction, defined in terms of the dynamics at infinity 
(Definition \ref{def:contracting_sequence}), 
and regularity, defined in terms of the asymptotics of orbits in $X$ 
(Definition \ref{def:pure_and_regular}),
are equivalent in a suitable sense,
compare the discussion in \cite[\S 5.2]{morse}. 
The most relevant aspect for the purposes of this paper is that regularity implies contraction
in a suitable sense.

We first consider sequences of isometries:
\begin{prop}
\label{prop:regimplcontr}
Every $\taumod$-regular sequence in $G$ contains a $\taumod$-contracting subsequence. 
\end{prop}
\proof
Compare the proof of \cite[Proposition 5.14]{morse}.

Suppose that the sequence $(g_n)$ in $G$ is $\taumod$-regular.
Let $x\in X$.
There exist simplices $\tau_n^{\pm}\in\Flagpmt$
(unique for large $n$)
such that 
$$ g_n^{\pm1}x \in V(x,\st(\tau_n^{\pm})) .$$
After passing to a subsequence,
we may assume convergence $$\tau_n^{\pm}\to\tau_{\pm}$$ in $\Flagpmt$,
because the flag manifolds are compact.

Since $x\in g_nV(x,\st(\tau_n^-))=V(g_nx,\st(g_n\tau_n^-))$,
it follows together with $g_nx \in V(x,\st(\tau_n^+))$
that the Weyl cones $V(g_nx,\st(g_n\tau_n^-))$ and $V(x,\st(\tau_n^+))$ 
lie in the same parallel set, namely in $P(g_n\tau_n^-,\tau_n^+)$, and face in opposite directions. 
In particular, the simplices $g_n\tau_n^-$ and $\tau_n^+$ are $x$-opposite,
and thus $g_n\tau_n^-$ converges to the simplex $\hat\tau_+$ $x$-opposite to $\tau_+$,
$$ g_n\tau_n^- \to\hat\tau_+ .$$
Since the sequence $(g_nx)$ is $\taumod$-regular,
it holds that 
$$ d(g_n^{-1}x_n,\D V(x,\st(\tau_n^-))) \to+\infty  $$
According to Lemma~\ref{lem:expconvsect},
for any $r,R>0$
the inclusion of shadows
$$ U_{\tau_n^-,x,R} \subset U_{\tau_n^-,g_n^{-1}x,r} $$
holds for $n\geq n(r,R)$.
Therefore there exist positive numbers $R_n\to+\infty$ and $r_n\to0$ such that 
$$ U_{\tau_n^-,x,R_n} \subset U_{\tau_n^-,g_n^{-1}x,r_n} $$
for large $n$,
equivalently
\begin{equation}
\label{eq:mpfshd}
g_nU_{\tau_n^-,x,R_n} \subset U_{g_n\tau_n^-,x,r_n} .
\end{equation}
Since $\tau_n^-\to\tau_-$ and $R_n\to+\infty$,
the sequence of shadows $U_{\tau_n^-,x,R_n}\subset C(\tau_n^-)\subset\Flagt$ {\em exhausts} $C(\tau_-)$
in the sense that every compactum in $C(\tau_-)$ is contained in $U_{\tau_n^-,x,R_n}$ for large $n$.
Indeed, for fixed $R>0$ we have Hausdorff convergence 
$U_{\tau_n^-,x,R}\to U_{\tau_-,x,R}$ in $\Flagt$,
which immediately follows e.g.\ using symmetry,
i.e.\ from the transitivity of the action $K_x\acts\Flagit$
of the maximal compact stabilizer $K_x<G$ of $x$.
Furthermore,
the shadows $U_{\tau_-,x,R}$ exhaust $C(\tau_-)$ as $R\to+\infty$,
cf.\ the continuity part of Lemma~\ref{lem:contpr}.

On the other hand,
since $g_n\tau_n^-\to\hat\tau_+$ and $r_n\to0$,
the shadows 
$U_{g_n\tau_n^-,x,r_n}$ {\em shrink}, i.e.\ Hausdorff converge to the point $\tau_+$.
Indeed, $U_{g_n\tau_n^-,x,r}\to U_{\hat\tau_+,x,r}$ in $\Flagt$ for fixed $r>0$,
and $U_{\hat\tau_+,x,r}\to\tau_+$ as $r\to0$,
using again the continuity part of Lemma~\ref{lem:contpr}
and the fact that the function (\ref{eq:distfrpar}) assumes the value zero only in $\tau_+$.

Together with these observations on exhaustion and shrinking of shadows, 
(\ref{eq:mpfshd}) shows that 
$$ g_n|_{C(\tau_-)} \to \tau_+ $$
uniformly on compacta,
i.e.\ the (sub)sequence $(g_n)$ is $\taumod$-contracting. 
\qed

\begin{rem}
The converse, that $\taumod$-contracting sequences in $G$ are $\taumod$-regular,
was shown in \cite[Theorem 5.23]{morse}.
\end{rem}
We conclude for groups of isometries:
\begin{cor}
\label{cor:regimplcontrgp}
$\taumod$-Regular subgroups are $\taumod$-convergence subgroups.
\end{cor}

\begin{rem}
For $\taumod$-regular subgroups,
the notion of $\taumod$-limit set introduced in Definition~\ref{dfn:lims}
is equivalent
to the notion of $\taumod$-limit set introduced in \cite[Def.\ 5.32]{morse}, see Proposition 5.29 of \cite{morse}. 
\end{rem}

Based on the corollary,
we can translate our proper discontinuity result for convergence subgroups 
(Theorem~\ref{thm:pdwconv})
into one for regular subgroups:
\begin{thm}[Domains of proper discontinuity for $\taumod$-regular subgroups]
\label{thm:pdwreg}
Let $\taumod\subseteq\simod$ be an arbitrary face type.
If $\Ga<G$ is a $\taumod$-regular subgroup,
then for every $\Wt$-left invariant thickening $\Th\subset W$
the action 
$$ \Ga\acts\geo X-(\Th^c(\Latm)\cup\Th(\Latp)) $$
is properly discontinuous. 
In particular,
if $\taumod$ is $\iota$-invariant and $\Th$ is fat,
then the action 
$$ \Ga\acts\geo X-\Th(\Lat) $$
is properly discontinuous. 
\end{thm}

\subsection{Discrete subgroups}
\label{sec:discr}

The general construction of domains of proper discontinuity in section~\ref{sec:accac}
applies equally to arbitrary discrete subgroups $\Ga<G$.
There are several ways to proceed. 
The most immediate possibility is the following.

Choose for every face type $\taumod\subseteq\simod$ 
a $\Wt$-left invariant thickening $\Th_{\taumod}$,
and define the $\Ga$-invariant compact families
$${\mathcal A}_{\pm} := \bigcup_{\taumod\subseteq\simod}{\mathcal A}^{\pm}_{\taumod,Th_{\taumod}} .$$
\begin{lem}
\label{lem:discracc}
The action $\Ga\acts\geo X$ is $({\mathcal A}_-,{\mathcal A}_+)$-accumulating.
\end{lem}
\proof
According to Lemma~\ref{lem:obspureg},
every sequence $\ga_n\to\infty$ in $\Ga$ contains a $\taumod$-regular (even $\taumod$-pure) subsequence,
and hence a $\taumod$-contracting subsequence
for some face type $\taumod$. 
The assertion follows therefore from Corollary~\ref{cor:orbacc}.
\qed

\medskip
Thus Proposition~\ref{prop:pd2} yields in this case:
\begin{prop}[Domains of proper discontinuity for discrete subgroups I]
\label{prop:pddisc}
If $\Ga<G$ is a discrete subgroup,
then the action 
\begin{equation}
\label{eq:dompd}
\Ga\acts\geo X-\bigcup_{\taumod}(\Th_{\taumod}^c(\Latm)\cup\Th_{\taumod}(\Latp)) 
\end{equation}
is properly discontinuous. 
\end{prop}
In general, this domain of proper discontinuity can be further enlarged
by only removing the thickenings of the limit simplices arising from {\em pure} sequences in the group: 
Define the {\em pure forward/backward $\taumod$-limit set} 
$$\Lat^{pure,\pm}\subseteq\Lat^{\pm}$$
as the closure of the set 
of all simplices $\tau_{\pm}$ as in (\ref{eq:contrtau}) 
for all $\taumod$-pure $\taumod$-contracting sequences $(\ga_n)$ in $\Ga$.
As above, we conclude:
\begin{prop}[Domains of proper discontinuity for discrete subgroups II]
\label{prop:pddisclar}
If $\Ga<G$ is a discrete subgroup,
then the action 
\begin{equation}
\label{eq:dompdlar}
\Ga\acts\geo X-
\bigcup_{\taumod}(\Th_{\taumod}^c(\Lat^{pure,-})\cup\Th_{\taumod}(\Lat^{pure,+}))
\end{equation}
is properly discontinuous. 
\end{prop}
Since the domain in (\ref{eq:dompd}) is in general smaller than the domain in (\ref{eq:dompdlar}),
one cannot expect the $\Ga$-action on it to be cocompact.

On the other hand, 
if $\Ga$ is $\taumod$-regular, 
then it contains $\numod$-pure sequences only for the face types $\numod\supseteq\taumod$,
and hence only these limit sets $\Lan^{\pm}$ can be nonempty.
Since $\Wn\leq \Wt$, 
we may choose $\Th_{\numod}=\Th_{\taumod}$ for these face types,
and then the domain in (\ref{eq:dompdlar}) coincides with the domain in Theorem~\ref{thm:pdwreg}.

\section{Cocompactness}
\label{sec:cocompact}

\subsection{Nearby simplex thickenings}
\label{sec:nearbythickenings}

For incident faces $\ups\subset\tau\subset\geo X$, 
the parabolic subgroups fixing them are contained in each other, 
$P_{\ups}\supset P_{\tau}$. 
Correspondingly, 
for incident face types $\upsmod\subset\taumod\subseteq\simod$
there is the natural forgetful map 
\begin{equation*}
\pi_{\upsmod\taumod}:\Flagt\to\Flagu
\end{equation*}
assigning to a face $\tau$ of type $\taumod$ 
its face $\ups$ of type $\upsmod$.  
It is a $G$-equivariant smooth fibration 
with compact base and fiber. 

We fix auxiliary Riemannian metrics 
on all partial flag manifolds $\Flagt$. 
Thereby 
also the $G$-orbits $G\xi\subset\geo X$ 
are equipped with Riemannian metrics 
by equivariantly identifying them with the appropriate flag manifolds. 

The fibrations
$\pi_{\upsmod\taumod}$
are then Lipschitz continuous by compactness. 
Vice versa, we have:
\begin{lem}
[Controlled lifts]
\label{lem:liftsimpnear}
Let $\tau$ and $\ups'$ be simplices of types 
$\taumod$ and $\upsmod$, 
$\upsmod\subset\taumod$, 
and let $\ups\subset\tau$ be the face of type $\upsmod$. 
Then there exists a simplex $\tau'\supset\ups'$ of type $\taumod$ 
such that 
\begin{equation*}
d(\tau',\tau)\leq C_0\cdot d(\ups',\ups)
\end{equation*}
with a uniform constant $C_0\geq1$ only depending 
on the chosen Riemannian metrics. 
\end{lem}
\proof
The Riemannian metrics on 
$\Flagt$ and $\Flagu$ 
can be chosen so that 
$\pi_{\upsmod\taumod}$ 
becomes a Riemannian submersion. 
With respect to these metrics, there exists $\tau'$ so that 
$d(\tau,\tau')=d(\ups,\ups')$. 
For other choices of the metrics, a multiplicative constant enters. 
\qed

\medskip
The lemma generalizes (by induction along galleries) to:
\begin{lem}
\label{lem:liftgallnear}
Let $\tau,\tau'$ be simplices of type $\taumod$ 
and let $\tilde\tau$ be a simplex
of type $\tiltaumod$. 
Then there exists another simplex $\tilde\tau'$ of type $\tiltaumod$ 
with relative position 
$\pos(\tilde\tau',\tau')=\pos(\tilde\tau,\tau)$
such that 
\begin{equation*}
d(\tilde\tau',\tilde\tau) \leq C_1\cdot d(\tau',\tau)
\end{equation*}
with a uniform constant $C_1\geq1$ only depending 
on the chosen Riemannian metrics. 
\end{lem}

Let now $G\xi\subset\geo X$ 
be a $G$-orbit at infinity,
which we think of as identified with the appropriate flag manifold.
We fix a $\Wt$-left invariant thickening $\Th\subset W$.
Then the distance between simplices $\tau$ in $\Flagt$
and the Hausdorff distance between their thickenings $\Th(\tau)\cap G\xi$ in $G\xi$ control each other,
and through an ideal point in $G\xi$ outside a simplex thickening 
exists a simplex thickening at controlled distance:
\begin{lem}[Nearby simplex thickenings]
\label{lem:findsi}
The following assertions hold 
with a uniform constant $C\geq1$ only depending 
on the chosen Riemannian metrics: 

\noindent 
(i)
The Hausdorff distance between the thickenings 
of any two simplices $\tau',\tau\in\Flagt$
is controlled by
\begin{equation*}
d_H(\Th(\tau')\cap G\xi,\Th(\tau)\cap G\xi)\leq C\cdot d(\tau',\tau) .
\end{equation*}

\noindent
(ii) 
For a point $\xi'\in G\xi$ and a simplex $\tau\in\Flagt$
there exists a simplex $\tau'\in\Flagt$ such that 
$\xi'\in\Th(\tau')$ and 
\begin{equation*}
d(\tau',\tau)\leq C\cdot d(\xi',\Th(\tau)\cap G\xi) .
\end{equation*}
\end{lem}
\proof
(i) 
If $\eta\in\Th(\tau)\cap G\xi$ is arbitrary, 
then applying Lemma~\ref{lem:liftgallnear} 
(to the flag manifold identified with $G\xi$)
yields a point $\eta'\in G\xi$ 
with $\pos(\eta',\tau')=\pos(\eta,\tau)\in\Th(\taumod)$,
i.e.\ $\eta'\in\Th(\tau')$,
and controlled distance 
$d(\eta',\eta)\leq C\cdot d(\tau,\tau')$. 

(ii) 
Suppose that $\xi\in\Th(\tau)\cap G\xi$ is the point closest to $\xi'$, 
i.e.\ $d(\xi',\xi)=d(\xi',\Th(\tau)\cap G\xi)$. 
Lemma~\ref{lem:liftgallnear} 
yields a simplex $\tau'\in\Flagt$
with $\pos(\xi',\tau')=\pos(\xi,\tau)\in\Th$ 
and controlled distance 
$d(\tau',\tau)\leq C\cdot d(\xi',\xi)$, 
whence the second inequality. 
\qed

\subsection{From expansion to transverse expansion}
\label{sec:exptotrexp}

Let $\taumod\subseteq\simod$ be a face type,
and let $\Th\subset W$ be a $\Wt$-left invariant thickening. 
In this section, we work on a fixed but arbitrary $G$-orbit $G\eta\subset\geo X$.

We start with an observation concerning the topology of thickenings in flag manifolds.
\begin{lem}[Fibration of thickenings]
\label{lem:thickfib}
Let $A\subset\Flagt$ be compact, 
and suppose that 
the thickenings $\Th(\tau)\cap G\eta$ of the simplices $\tau\in A$
are pairwise disjoint. 
Then the natural map 
$$ \Th(A)\cap G\eta\stackrel {\pi}{\lra} A$$
is a continuous fibration with compact fiber.
\end{lem}
\proof
Suppose that $\xi_n\to\xi$ in $\Th(A)\cap G\eta$ and $\tau_n\to\tau$ in $A$
with $\xi_n\in\Th(\tau_n)$.
Then $\xi\in\Th(\tau)$ by 
semicontinuity of relative position, 
cf.\ Lemma~\ref{lem:semcontrelpos}.
The assumption on the disjointness of fibers implies that $\pi(\xi)=\tau$. 
Thus, $\pi$ is continuous.

In order to show that $\pi$ is a fiber bundle, we need to construct local trivializations. 
Fix $\tau_0\in A$.
There exists a compact subset $S\subset G$ 
which is mapped by $s\mapsto s\tau_0$ homeomorphically onto a compact neighborhood of $\tau_0$ in $A$.
(Such a subset can be found in a slice through $e$ transverse to $P_{\tau_0}$.)
Restricting the action $G\acts\geo X$,
we obtain a topological embedding 
$$ S\times(\Th(\tau_0)\cap G\eta)\to\Th(A)\cap G\eta$$
and a local trivialization of $\pi$ over a neighborhood of $\tau_0$ in $A$.
\qed

\medskip
Now we turn to dynamics.

Let $(g_n)$ be a sequence of isometries in $G$ which preserve $A$,
$g_nA=A$. 
We consider the action of $(g_n)$ on $G\eta$
and derive transverse expansion from expansion on $\Flagt$:
\begin{lem}
[Expansion implies transverse expansion]
\label{lem:expimpltrexp}
Suppose that $(g_n)$ is on $\Flagt$ arbitrarily expanding at $\tau_+\in A$,
i.e.\ there exist neighborhoods $V_n$ of $\tau_+$ in $\Flagt$ 
and constants $c_n\to+\infty$
such that $g_n|_{V_n}$ is expanding with expansion factor $c_n$.

Then there exist neighborhoods $W_n$ of $\Th(\tau_+)\cap G\eta$ 
and constants $C_n\to+\infty$ such that 
\begin{equation}
\label{ineq:trexpanA}
d(g_n\xi,g_n\Th(\tau)\cap G\eta)\geq C_n\cdot d(\xi,\Th(\tau)\cap G\eta)
\end{equation}
for all $\xi\in W_n$ and $\tau\in A$ with $\Th(\tau)\cap G\eta\subset W_n$, 
compare inequality (\ref{ineq:trexpan}).
\end{lem}
\proof
To simplify notation, we write (only in this proof) $\Th(\cdot)$ instead of $\Th(\cdot)\cap G\eta$.

Let $W_n$ be some (small) neighborhood of the compact subset $\Th(\tau_+)$,
and let $\xi,\tau$ be as in inequality (\ref{ineq:trexpanA}).
We work near $\Th(g_n\tau_+)$.
According to Lemma~\ref{lem:findsi}(ii), 
we can choose $g_n\tau'\in\Flagt$ such that $g_n\xi\in\Th(g_n\tau')$, 
and $g_n\tau'$ has controlled distance from $g_n\tau$,
$$ d(g_n\tau',g_n\tau) \leq C\cdot d(g_n\xi,\Th(g_n\tau))$$
with a uniform constant. 

After shrinking the neighborhood $g_nW_n$ of $\Th(g_n\tau_+)$, 
we may assume that $g_n\tau$ is close to $g_n\tau_+$, 
using that $\Th(A)$ fibers over $A$, cf.\ Lemma~\ref{lem:thickfib}.
Moreover,
that $g_n\xi$ is close to $\Th(g_n\tau)$.
Thus, after shrinking $W_n$ sufficiently, we may assume that 
$\tau',\tau\in V_n$.

Then 
$$ d(g_n\tau',g_n\tau) \geq c_n\cdot d(\tau',\tau) .$$
Since we also have uniform control 
$$ d(\xi,\Th(\tau)) \leq d_H(\Th(\tau'),\Th(\tau)) \leq C\cdot d(\tau',\tau)$$
by Lemma~\ref{lem:findsi}(ii), 
it follows that 
$$ d(g_n\xi,\Th(g_n\tau)) \geq C^{-2}c_n\cdot d(\xi,\Th(\tau)) ,$$
that is, our assertion with $C_n=C^{-2}c_n$.
\qed

\medskip
We apply the above discussion to discrete group actions on flag manifolds.
\begin{prop}[Transverse expansion at slim thickenings]
\label{prop:conicimplexpatrtau}
Let $\Ga<G$ be a discrete subgroup 
and suppose that:

(i) The action $\Ga\acts\Flagt$ is expanding at $\Latp$.

(ii) The thickenings $\Th(\tau)\cap G\eta$ of the simplices $\tau\in\Latp$ are pairwise disjoint. 

Then the action $\Ga\acts G\eta$ is expanding at $\Th(\Latp)\cap G\eta$
transversely to the natural fibration $\Th(\Latp)\cap G\eta\to\Latp$
given by Lemma~\ref{lem:thickfib}.
\end{prop}
\proof
Since the action $\Ga\acts\Flagt$ is expanding at $\Latp$, 
it is arbitrarily strongly expanding there, cf.\ Remark~\ref{rem:arbstrexp},
i.e.\ for every limit simplex $\tau_+\in\Latp$ exists a sequence 
$(\ga_n)$ in $\Ga$ and a sequence of neighborhoods $V_n$ of $\tau_+$ 
such that the $\ga_n|_{V_n}$ are uniformly expanding with expansion factors 
$c_n\to+\infty$.
Lemma~\ref{lem:expimpltrexp} then implies that the action 
$\Ga\acts G\eta$
is (arbitrarily) expanding at $\Th(\Latp)\cap G\eta$ transversely to $\pi$.
\qed

\medskip
Using that transverse expansion implies cocompactness on the complement 
(Proposition~\ref{prop:trexpcoco}), 
we derive our main cocompactness result:
\begin{thm}[Cocompact domains]
\label{thm:cocomain}
Let $\Ga$ and $\Th$ be as in the previous proposition. 
Then 
the action 
$$\Ga\acts G\eta-\Th(\Latp)$$ 
is cocompact. 
\end{thm}

The following is a special case of the theorem.
Here, 
for $\iota$-invariant $\taumod$,
a subset of $\Flagt$ is called {\em antipodal} if the simplices in it are pairwise opposite,
cf.\ Definition~\ref{dfn:antip}(ii).

\begin{cor}[Cocompactness outside slim thickenings]
\label{cor:cocomain} 
Let $\taumod\subseteq\simod$ be an $\iota$-invariant face type. 
Suppose that $\Ga<G$ is a discrete subgroup such that 
$\Lat$ is antipodal 
and the action $\Ga\acts\Flagt$ is expanding at $\Lat$.

(i) Then for any slim $\Wt$-left invariant thickening $\Th\subset W$
the action 
$$\Ga\acts \DF X-\ThF(\Lat)$$ 
is cocompact. 

(ii) More generally, 
suppose that $\numod\subseteq\simod$ is another face type
and that the thickening $\Th$ is also $\Wn$-right invariant.
Then for any $G$-orbit $G\eta\subset\geo X$ of type $\bar\eta=\theta(\eta)\in\inte(\numod)$ 
the action 
$$\Ga\acts G\eta-\Th(\Lat)$$ 
is cocompact. 
\end{cor}
\proof
We have that $\Latpm=\Lat$, because $\taumod$ is $\iota$-invariant.
Since $\Th$ is slim and the simplices $\tau$ in $\Lat$ are pairwise antipodal,
their thickenings $\ThF(\tau)$ in $\DF X$ are pairwise disjoint,
cf.\ Lemma~\ref{lem:sldisj}. 
Thus, the hypotheses of the theorem are satisfied.
\qed

\begin{rem}[Rank one]
If $\rank(X)=1$, this follows from part of a basic result for Kleinian groups
characterizing convex-cocompactness.
Namely, the following properties are equivalent for a discrete subgroup $\Ga<G$:

(i) $\Ga$ is {\em convex-cocompact}.

(ii) The action $\Ga\acts\ol X$, equivalently, the action $\Ga\acts\geo X$,  is expanding at $\La$.

(iii) The (properly discontinuous) action $\Ga\acts \ol X-\La$ is cocompact. 

\no
In particular, then the action $\Ga\acts \geo X-\La$ is cocompact. 
\end{rem}

\subsection{Cocompact domains of proper discontinuity}

We consider the following class of discrete subgroups
(see Definition~\ref{dfn:ceaintro} in the introduction):
\begin{dfn}[CEA subgroup]
\label{dfn:cea}
For a $\iota$-invariant face type $\taumod\subseteq\simod$ 
we call a $\taumod$-con\-ver\-gen\-ce subgroup $\Ga<G$ a {\em $\taumod$-CEA subgroup}
(convergence, expanding, antipodal)
if $\Lat$ is antipodal and if the action $\Ga\acts\Flagt$ is expanding at $\Lat$.
\end{dfn}
\begin{rem}[CEA versus Anosov]
The class of $\taumod$-CEA subgroups coincides with the class of $P_{\taumod}$-Anosov subgroups,
see \cite[\S 6.5]{morse}.
Here, $P_{\taumod}$ refers to the conjugacy class of parabolic subgroups of $G$ 
corresponding to the face $\taumod$ of the spherical Weyl chamber $\simod$.
\end{rem}

Combining our main results on proper discontinuity (Theorem~\ref{thm:pdwconv})
and cocompactness (Corollary~\ref{cor:cocomain}),
we obtain:
\begin{thm}[Cocompact domains of proper discontinuity]
\label{thm:cocodisconttau}
Suppose that $\Ga<G$ is a $\taumod$-CEA subgroup. 

(i) Then for any balanced $\Wt$-left invariant thickening $\Th\subset W$
the action 
$$\Ga\acts \DF X-\ThF(\Lat)$$ 
is properly discontinuous and cocompact. 

(ii) More generally, 
suppose that $\numod\subseteq\simod$ is another face type
and that the thickening $\Th$ is also $\Wn$-right invariant.
Then for every $G$-orbit $G\eta\subset\geo X$ of type $\bar\eta=\theta(\eta)\in\inte(\numod)$
the action 
$$\Ga\acts G\eta-\Th(\Lat)$$ 
is properly discontinuous and cocompact. 
\end{thm}
\begin{rem}
\label{rem:ddreg}
According to Corollary~\ref{cor:exbaltau},
balanced $\Wt$-left invariant thickenings always exist, 
and Theorem~\ref{thm:cocodisconttau} therefore provides 
cocompact domains of discontinuity at least in the Furstenberg boundary 
$\DF X$. 
\end{rem}

The question whether these domains are nonempty will be addressed in section~\ref{sec:nonempty}.

\subsection{A relation with Mumford's Geometric Invariant Theory}
\label{sec:git} 

We continue the discussion in Example~\ref{ex:config},
now looking at actions (of Lie subgroups) on configuration spaces.
(See \cite{KM1, KLM} for a more detailed  
discussion of Geometric Invariant Theory in the context of weighted configurations.) 

Let $H=\Isom_o(Y)$. 
We consider the diagonal action $H\acts\DF X$ on configurations.
As we discussed in Example~\ref{ex:config},
the choice of a regular vector $t=(t_i)\in\inte(\De)$
determines subsets 
$$(\DF X)_{st,t}=\DF X-(\ol{\Th}_t)_{F\ddot u}(A)
\quad\hbox{ and }
(\DF X)_{sst,t}=\DF X-(\Th_t)_{F\ddot u}(A)$$
of stable, respectively, semistable weighted configurations in $\geo Y$.
Mumford's GIT \cite{Mumford} defines the {\em Mumford quotient} 
$$
\DF X//_t H= (\DF X)_{sst,t}//H.
$$  
by suitably extended orbit equivalence.
In the case when the thickening $\Th_t$ is balanced, 
all semistable points are even stable,
and one has
$$
\DF X//_t H= (\DF X)_{sst,t}//H= (\DF X)_{st,t}/H,
$$
the latter being a quotient in the usual sense. 

A nice exercise is to prove directly that the space $\DF X//_t H$ is compact and
Hausdorff in this case. 
For instance, if $H=PSL(2,\R)$, $Y=\H^2$, $n=3$ and $t=(1,1,1)$, then
$\DF X//_{t} H$ consists of exactly two points represented by configurations
of three distinct points on the circle with different cyclic orders.
Continuing with $Y=\H^2$ and letting $n=4$, one verifies that for
${t}=(2,1,1,1)$ the Mumford quotient is homeomorphic to $S^1$, while for
${t}=(5,4,3,1)$ the Mumford quotient is homeomorphic to the disjoint union
of two circles. Taking $n=5$, one obtains that for 
${t}=(1,1,1,1,1)$ the Mumford quotient is the genus 4 oriented surface,
while for ${t}= (5,4,1,1,1)$ the quotient is the disjoint union of two
2-spheres. Thus, we see that quotients are not homeomorphic for distinct
choices of ${t}$'s. 

More generally, one can describe dependence of the topology of the Mumford quotient 
$\DF X//_t H$ on the parameter $t$ as follows. 

The hyperplanes $\sum_{i\in I} t_i = \sum_{j\notin I} t_j$ (also called
{\em interior walls}), where $I$ runs over subsets of $\{1,\ldots,n\}$,
partition the chamber 
$$
\Delta=\{(t_1,\ldots,t_n): t_i>0\}
$$
into open convex subsets, also called {\em chambers}. The topology of
$\DF X//_{t} H$ does not change as long as ${t}$ varies in a single chamber;
permuting the chambers does not change the topology either; however, {\em
crossing through a wall} amounts to a certain Morse surgery on the
manifold. This can be seen by identifying the quotients $\DF X//_{t} H$ with certain
moduli spaces of polygons with fixed side-length: In the case when
$H=PSL(2,\R)$, these are polygons in the Euclidean plane, cf.\ \cite{KM1}. 

It was conjectured by Kevin Walker that if ${t}, {t}'$ belong to chambers
in distinct $S_n$-orbits then the Mumford quotients are not homeomorphic. This
conjecture was proven 20 years later in ``most'' cases by Farber, Hausmann
and Sch\"utz \cite{FHS} and in full generality by  Sch\"utz \cite{Schutz}.
Similar results hold when the circle is replaced by a $k$-sphere. In fact,
different quotients are distinguished by their cohomology rings.

We will now see how the dependence of the topology of $\DF X//_t H$ on the parameter $t$ described above 
leads to the change of the topology of quotients by {\em discrete} group actions. 

\begin{example}
\label{ex:product_actions}
We continue with the notation of Example \ref{ex:config}. For concreteness, we
assume that $Y=\H^2$, $H=PSL(2,\R)$ and  
$\Ga< H$ is a torsion-free uniform lattice (a closed hyperbolic
surface subgroup). 
The embedding $H<G=H\times \ldots\times H$ is diagonal
and we view $\Ga$ as a subgroup of $G$.
Then $\Ga$ preserves the diagonally embedded totally-geodesic hyperbolic plane $\H^2\subset X$
and acts cocompactly on it.
Thus, 
$\Las=\geo\H^2 \subset\DF X$,  the diagonally embedded circle,
and $\La=\geo\H^2\subset\geo X$ for the visual limit set.
The ideal boundary points in $\geo\H^2\subset\geo X$
are contained in the central regular $G$-orbit $\theta^{-1}(\bar\zeta)\subset\geo X$
of type $\bar\zeta\in\inte(\simod)$
represented by the vector $(1,\ldots, 1)\in\inte(\De)$.
It follows that the subgroup $\Ga<G$
is uniformly $\simod$-regular (see \cite{morse} for the precise definition).
More precisely, 
it is $\{\bar\zeta\}$-regular.
Moreover,
the group $\Ga$ is obviously quasi-isometrically embedded in
$H$, and hence also in $G$. 
We conclude that $\Ga<G$ is a $\simod$-CEA subgroup 
(e.g.\ as a consequence of \cite[Theorem 1.5]{mlem}). 

Given a balanced metric thickening $\Th=\Th_t\subset W$, 
the domain $\Om_{\Th}=\DF X-\ThF(\Las)$
considered in Theorem~\ref{thm:cocodisconttau}(i) 
equals the set $(\DF X)_{st, t}$ of stable weighted $n$-point
configurations on $\geo Y\cong S^1$ (stability being defined with respect to the weights $t$).
The group $H$ acts on $(\DF X)_{st, t}$ freely and we have a principal $H$-bundle
$$
H\to (\DF X)_{st, t}\to (\DF X)_{st,t}/H= \DF X//_{t} H. 
$$
Dividing $(\DF X)_{st,t}$ by $\Ga$ instead of $H$ we obtain a fiber bundle
$$
H/\Ga \to (\DF X)_{st,t}/\Ga \to (\DF X)_{st,t}/H .
$$
In particular, 
by taking non-homeomorphic Mumford quotients, we obtain non-homeomorphic
quotients $\Om_{\Th}/\Ga$. For instance, for $n=4$ we obtain three distinct topological
types of quotients: The empty quotient, a connected nonempty quotient (a bundle over the circle with
the fiber $H/\Ga$) and a disconnected quotient which is the disjoint
union of two copies of an $H/\Ga$-bundle over $S^1$. 
\end{example}

\section{Nonemptiness}
\label{sec:nonempty}

\subsection{Thickenings and packings}

We will use the following notion of ball packing for the visual boundary 
(using its structure as a topological spherical building).
\begin{dfn}[Packing]
A {\em packing} of $\geo X$ by $\pihalf$-balls 
is a family ${\cal B}$ of disjoint open $\pihalf$-balls (with respect to the Tits metric)
the union of whose closures equals $\geo X$. 
We call the packing {\em compact}
if the set of centers of these balls is compact with respect to the visual topology.
\end{dfn}
Note that the set of centers of the balls is necessarily antipodal,
cf.\ Definition~\ref{dfn:antip}(i),
and hence the centers must have the same $\iota$-invariant type. 
We call it the {\em type} of the packing.
We call the packing {\em simplicial} if the balls are simplicial subcomplexes of $\geo X$.
The simplicial $\pihalf$-balls 
are precisely the $\pihalf$-balls centered at points of root type,
and hence a packing is simplicial if and only if it is of root type. 

\medskip
We will show that compact packings often do not exist.

\begin{dfn}[Non-packing type]
We say that the symmetric space $X$ is of 

(i) {\em non-packing type} 
if $\geo X$ admits no compact packing by $\pihalf$-balls. 

(ii) {\em non-$\bar\vartheta$-packing type} for an $\iota$-invariant type $\bar\vartheta\in \simod$ 
if $\geo X$ admits no compact packing by $\pihalf$-balls of type $\bar\vartheta$. 

(iii) {\em non-root packing type} 
if it is of non-$\bar\vartheta$-packing type for {\em some} root type $\bar\vartheta\in\simod$.
\end{dfn}

\medskip
Our motivation for proving the nonexistence of packings is 
that it implies via the nonfullness of thickenings,
as is made precise by the next result,
the nonemptiness of domains of proper discontinuity,
see Proposition~\ref{prop:neccdd} below.

Let $\taumod\subseteq\simod$ be an $\iota$-invariant face type.
Suppose that $$A\subset\Flagt$$ is an antipodal compact subset.
It determines for every $\iota$-invariant type $\bar\vartheta_0\in\taumod$
the antipodal compact subset $C\subset\geo X$ consisting of the points 
$\zeta_{\tau,\bar\vartheta_0}=\tau\cap\theta^{-1}(\bar\vartheta_0)$ of type $\bar\vartheta_0$
in the simplices $\tau\in A$,
and hence the family of disjoint open $\pihalf$-balls 
$$ {\mathcal B}(A,\bar\vartheta_0) = \Bigl\{B\bigl(\zeta_{\tau,\bar\vartheta_0},\pihalf\bigr):\tau\in A\Bigr\} .$$
Note that the union of the closed balls is a compact subset of $\geo X$.

\begin{prop}[Full thickenings yield packings]
\label{prop:emptdpack}
Let $\taumod$ and $\bar\vartheta_0\in\taumod$ be $\iota$-invariant,
and let $A\subset\Flagt$ be an antipodal compact subset.
Suppose that $$\ThF(A)=\DF X$$ for all balanced $\Wt$-left invariant thickenings $\Th\subset W$
of the form 
$\Th=\Th_{\bar\vartheta_0,\bar\vartheta,\pihalf}$
(as defined by (\ref{eq:metth})).
Then the family of balls 
${\mathcal B}(A,\bar\vartheta_0)$ is a packing of $\geo X$.
\end{prop}
\proof
Suppose that ${\mathcal B}(A,\bar\vartheta_0)$ is not a packing.
The union of the corresponding closed balls is compact in $\geo X$
(with respect to the visual topology),
and its complement therefore open.
Let $\xi$ be a point in the complement, and denote $\theta(\xi)=\bar\vartheta$.
After perturbing $\xi$,
we may assume that $\xi$ is regular and that 
the (always fat) $\Wt$-left invariant metric thickening $\Th_{\bar\vartheta_0,\bar\vartheta,\pihalf}\subset W$ is balanced,
cf.\ Lemma~\ref{lem:metth} and the proof of Corollary~\ref{cor:exbaltau}.
By the construction of metric thickenings, 
see (\ref{eq:metth}) and (\ref{eq:metthtau}),
$$\Th_{\bar\vartheta_0,\bar\vartheta,\pihalf}(\tau) \cap G\xi = \ol B\bigl(\zeta_{\tau,\bar\vartheta_0},\pihalf) \cap G\xi $$
for $\tau\in A$.
It follows that 
$\xi\not\in\Th_{\bar\vartheta_0,\bar\vartheta,\pihalf}(A)$
and hence $(\Th_{\bar\vartheta_0,\bar\vartheta,\pihalf})_{F\ddot u}(A)\neq\DF X$. 
\qed

\begin{rem}[Full thickenings yield fibrations]
\label{rem:emptdfib}
Note that the hypothesis of the proposition implies in particular the existence of the following kind of fibrations 
of the Furstenberg boundary:
In view of Lemma~\ref{lem:thickfib}, 
it follows from $\ThF(A)=\DF X$ 
that there is a fiber bundle
$$ \DF X\lra A $$
whose fibers are finite unions of Schubert cycles
(namely the thickenings $\Th(\tau)$ for $\tau\in A$).
Any two fibers are equivalent modulo the $G$-action on $\DF X$.
If $A=\Lat$ for a subgroup $\Ga<G$, then the fibration is $\Ga$-equivariant. 
\end{rem}

\subsection{Nonexistence of packings}

We show in this section that compact packings of the visual boundary by $\pihalf$-balls 
do not exist for most Weyl groups.
Note that the discussion applies more generally to packings of {\em compact topological spherical buildings}.

\subsubsection{Type $A_2$}

Suppose that the symmetric space $X$ has type $A_2$.
The spherical model chamber $\simod$ is then an arc $\bar\xi\bar\eta$ of length $\pithird$,
with $\bar\xi,\bar\eta\in\simod$ the two vertex types. 
The involution $\iota$ of $\simod$ is the reflection at the midpoint $\bar\zeta$,
which is therefore the only $\iota$-invariant type.
We denote by 
$$\Flags\stackrel{\pi_{\bar\xi}}{\lra}\Flag_{\bar\xi}=\theta^{-1}(\bar\xi)
\quad\hbox{ and }\quad
\Flags\stackrel{\pi_{\bar\eta}}{\lra}\Flag_{\bar\eta}=\theta^{-1}(\bar\eta)$$
the canonical projections from the full flag manifold (of chambers)
to the partial flag manifolds (of vertices of fixed type).

A packing ${\cal B}$ of $\geo X$ by $\pihalf$-balls is necessarily of type $\bar\zeta$ and hence simplicial.
A $\pihalf$-ball in $\geo X$ with center of type $\bar\zeta$ consists of a central chamber 
and all chambers adjacent to it,
i.e. it is the $\pithird$-neighborhood of its central chamber. 
Thus, the packing corresponds to a set $C\subset\Flags$ of pairwise opposite chambers 
such that every other chamber is adjacent to a chamber in $C$. 

We denote by $C_{\bar\xi}=\pi_{\bar\xi}(C)$ and $C_{\bar\eta}=\pi_{\bar\eta}(C)$ 
the sets of vertices of the chambers in $C$,
and by $O_{\bar\xi}=\Flag_{\bar\xi}-C_{\bar\xi}$ and $O_{\bar\eta}=\Flag_{\bar\eta}-C_{\bar\eta}$ 
their complements.
The complement of the union of the chambers in $C$ 
is the union of the open $\pithird$-balls centered at the points in $O_{\bar\xi}\cup O_{\bar\eta}$,
i.e. the chambers not in $C$ are the chambers with a vertex in $O_{\bar\xi}$ or $O_{\bar\eta}$.
We therefore have the disjoint decomposition
$$ \Flags=C\sqcup\pi_{\bar\xi}^{-1}(O_{\bar\xi})\sqcup\pi_{\bar\eta}^{-1}(O_{\bar\eta}) .$$
We observe that, if a chamber has a vertex in $O_{\bar\xi}$,
then its other vertex lies in $C_{\bar\eta}$.
Vice versa, 
every vertex in $C_{\bar\eta}$ belongs to a chamber whose other vertex lies in $O_{\bar\xi}$.
This means that 
\begin{equation}
\label{eq:chdecoa2}
C_{\bar\eta} = \pi_{\bar\eta}(\pi_{\bar\xi}^{-1}(O_{\bar\xi}))
\end{equation}
So far, our discussion applies to packings of arbitrary spherical buildings of type $A_2$.
Now we take into account the visual topology.
\begin{thm}
\label{thm:nopacka2}
If $X$ has type $A_2$,
then it is of non-packing type.
\end{thm}
\proof
We keep the notation from the previous discussion. 
Suppose that ${\cal B}$ is a compact packing of $\geo X$,
i.e.\ $C$ is compact 
and therefore also its images $C_{\bar\xi}$ and $C_{\bar\eta}$
under the projections $\pi_{\bar\xi}$ and $\pi_{\bar\eta}$.
Then $O_{\bar\xi}$ is open. 
Since the projection $\pi_{\bar\eta}$ is open,
(\ref{eq:chdecoa2}) implies that $C_{\bar\eta}$ is also open,
i.e.\ it is clopen. 
Since it is a nonempty proper subset,
it follows that $\Flag_{\bar\eta}$ is disconnected,
and consequently also $\Flags$.
This is absurd,
because $\Flags$ is a homogeneous space of $\Isom_o(X)$ and therefore connected.
\qed

\subsubsection{Irreducible case of rank $\geq3$}

\medskip
For most irreducible Weyl groups, 
the question of the nonexistence of {\em simplicial}
packings can be reduced to the $A_2$-case.
\begin{thm}[Nonexistence of simplicial packings in rank $\geq3$]
\label{thm:nopackr3}
If $X$ is irreducible of $\rank(X)\geq3$, 
then it is of non-root packing type.
\end{thm}
\proof
We make use of the spherical building geometry of $\tits X$, 
see \cite{qirigid} for a detailed discussion. 
The question of nonexistence can be reduced to lower rank 
by observing that packings of spherical buildings by $\pihalf$-balls
induce such packings of their spaces of directions.

The space of directions $\Si_{\xi}\tits X$ of a point $\xi\in\geo X$ 
carries again a natural spherical building structure. 
We will use the notation $(S_\xi,W_\xi)$ for the associated Coxeter complex. 
More precisely, $\Si_{\xi}\tits X$ is naturally identified with the Tits building of the symmetric subspace $X'\subset X$, 
which appears in the decomposition $P(l)= X'\times l$ of the parallel set 
of a geodesic $l\subset X$ asymptotic to $\xi$. 

Furthermore, 
the spaces of directions of closed $\pihalf$-balls $\ol B(\zeta,\pihalf)$ 
at boundary points $\xi\in\D\ol B(\zeta,\pihalf)$ are again $\pihalf$-balls,
$$\Si_{\xi}\ol B(\zeta,\pihalf)=\ol B(\oa{\xi\zeta},\pihalf) .$$
This follows from the first variation formula in $S^2$, 
because for any point $\eta$ sufficiently close to $\xi$ 
the three points $\xi,\eta,\zeta$ are the vertices 
of an embedded spherical triangle. 
If two open balls $B(\zeta_i,\pihalf)$ are disjoint 
and if $\xi$ is a point in the intersection of their boundaries, 
then the spaces of directions $\Si_{\xi}\ol B(\zeta_i,\pihalf)$ have disjoint interiors. 

Let now ${\mathcal B}$ be a compact packing of $\geo X$ by $\pihalf$-balls. 
Then ${\mathcal B}$ induces packings ${\mathcal B}_{\xi}$ by $\pihalf$-balls 
of the spaces of directions $\Si_{\xi}\tits X$ 
for all boundary points $\xi$ of the packing balls;
the family ${\mathcal B}_{\xi}$ consists of the balls $B(\oa{\xi\zeta},\pihalf)$
for which $B(\zeta,\pihalf)\in{\mathcal B}$ and $\xi\in\D B(\zeta,\pihalf)$. 
If the packing ${\mathcal B}$ is simplicial, then so are the packings ${\mathcal B}_{\xi}$.

To see that the compactness of ${\mathcal B}$ 
implies the compactness of the induced families ${\mathcal B}_{\xi}$, 
consider a convergent sequence $\zeta_n\to\zeta$ 
of centers of packing balls in ${\mathcal B}$ 
such that $\tangle(\xi,\zeta_n)=\pihalf$ for all $n$. 
Then $\tangle(\xi,\zeta)\leq\pihalf$ 
by the semicontinuity of Tits distance. 
However, 
strict inequality is impossible, 
because then $\xi$ would be an interior point of the packing ball 
$B(\zeta,\pihalf)$, 
which is absurd. 
Thus also $\tangle(\xi,\zeta)=\pihalf$ 
and it follows that 
$\oa{\xi\zeta_n}\to\oa{\xi\zeta}$.
Hence the centers $\oa{\xi\zeta_n}$ of packing balls in ${\mathcal B}_{\xi}$ 
converge to the center of such a ball. 
Thus the families ${\mathcal B}_{\xi}$ are compact.

Let us now focus on simplicial packings. 
Suppose that $\geo X$ admits compact simplicial packings by $\pihalf$-balls for all (at most two) root types.
Consider for all such packings of $\geo X$ the induced simplicial packings by $\pihalf$-balls 
of the spaces of directions $\Si_{\xi}\tits X$ 
for all vertices $\xi$ in the boundaries of packing balls.
Then for these, the type $\theta(\xi)\in\simod$ runs through all possible vertex types,
and for every fixed vertex type $\theta(\xi)$ 
the type of the packing ${\mathcal B}_{\xi}$ runs through all possible root types.
(The vertex type $\theta(\xi)$ and the root type of the packing ${\mathcal B}_{\xi}$
uniquely determine the root type of the packing ${\mathcal B}$ which has to be used.)

The type $\theta(\xi)$ of a vertex $\xi$ 
corresponds to a wall of the fundamental Weyl chamber, 
and the Dynkin diagram for the link $\Si_{\xi}\tits X$ 
is obtained from the Dynkin diagram for $\tits X$ 
by removing the corresponding node. 
By examining Dynkin diagrams of irreducible root systems,  
we note that every irreducible root system of rank $\ge 3$ has a simple edge
and can hence be reduced to the $A_2$ root system 
by successively removing nodes without disconnecting it.
Thus, if $\rank(X)\geq3$,
it follows that there exists a symmetric space $X'$ of type $A_2$ 
whose visual boundary admits a compact packing by $\pihalf$-balls. 
This contradicts Theorem~\ref{thm:nopacka2}.
\qed

\medskip
Note that root types are $\iota$-invariant.

Regarding the irreducible case,
Theorems~\ref{thm:nopacka2} and~\ref{thm:nopackr3}
leave open the cases of type $B_2$ and $G_2$ in rank 2.
We will prove some partial results for the $B_2$-case in section~\ref{sec:nonexb2}.

\subsubsection{Reducible case}

We reduce to the irreducible case using the observation:

\begin{lem}\label{lem:reduction}
Suppose that $X$ decomposes as the product $X=X_1\times X_2$ of symmetric spaces.
If $X_1$ is of non-root packing type, then so is $X$.
\end{lem}
\proof The model chamber $\simod$ of $X$ splits as the spherical join 
$$
\simod=\simod^1\circ \simod^2
$$
of the model chambers of the factors. 
The same applies to the visual  boundaries:
$$
\geo X= \geo X_1 \circ \geo X_2 
$$
A root type $\bar\vartheta_1\in\simod^1$ remains a root type in $\simod$ under the inclusion $\simod^1\subset\simod$,
and a $\pihalf$-ball $\ol B^{\geo X}(\zeta_1,\pihalf)\subset\geo X$ centered at a point 
$\zeta_1\in\geo X_1\subset\geo X$ of type $\bar\vartheta_1$ 
splits as the spherical join 
$$ \ol B^{\geo X}(\zeta_1,\pihalf)=\ol B^{\geo X_1}(\zeta_1,\pihalf)\circ\geo X_2$$
of the ball $\ol B^{\geo X_1}(\zeta_1,\pihalf)\subset\geo X_1$ 
with the full visual boundary $\geo X_2$. 
Hence, 
$\geo X$ admits a compact packing by $\pihalf$-balls of type $\bar\vartheta_1$
if and only if $\geo X_1$ does.
The assertion follows.
\qed

\medskip
Combining Theorems~\ref{thm:nopacka2} and~\ref{thm:nopackr3} with Lemma~\ref{lem:reduction}, we obtain:

\begin{thm}[Nonexistence of simplicial packings]
\label{thm:nopack}
If $X$ has at least one de Rham factor not of the type $A_1, B_2$ or $G_2$,
then it is of non-root packing type.
\end{thm}
\proof
The assumptions imply that $X$ has a de Rham factor of type $A_2$ or with rank $\geq3$.
\qed

\subsubsection{Type $B_2$}
\label{sec:nonexb2}

Suppose now that the symmetric space $X$ has type $B_2$.
We obtain only partial results on the nonexistence of packings.

The model spherical chamber $\simod$ is an arc $\bar\xi\bar\eta$ of length $\piforth$,
and $\bar\xi,\bar\eta\in\simod$ are the two vertex types. 
Moreover, $\iota=\id_{\simod}$ and all types in $\simod$ are $\iota$-invariant.
We again denote by 
$$\Flags\stackrel{\pi_{\bar\xi}}{\lra}\Flag_{\bar\xi}=\theta^{-1}(\bar\xi)
\quad\hbox{ and }\quad
\Flags\stackrel{\pi_{\bar\eta}}{\lra}\Flag_{\bar\eta}=\theta^{-1}(\bar\eta)$$
the canonical projections from the full flag manifold (of chambers)
to the partial flag manifolds (of vertices of fixed type).

Simplicial $\pihalf$-balls in $\tits X$ are centered at vertices,
and hence simplicial packings by $\pihalf$-balls are of vertex type. 

A packing ${\cal B}_{\bar\vartheta_0}$ of regular type $\bar\vartheta_0\in\inte(\simod)$
gives rise to a continuous family of packings ${\cal B}_{\bar\vartheta}$, $\bar\vartheta\in\simod$,
by simultaneously ``sliding'' its centers along the chambers containing them.
Namely, 
we choose as the centers of ${\cal B}_{\bar\vartheta}$ 
the points of type $\bar\vartheta$ in those chambers which contain the centers of ${\cal B}_{\bar\vartheta_0}$. 
A regular packing thus gives rise to singular packings of {\em both} vertex types.

Consider now a packing ${\cal B}$ of $\geo X$ by $\pihalf$-balls of type $\bar\xi$
with set of centers $C_{\bar\xi}\subset\Flag_{\bar\xi}$. 
Note that $\Flag_{\bar\eta}$ is {\em partitioned} by the $\piforth$-spheres around the points in $C_{\bar\xi}$,
and that the map 
\begin{equation*}
\pi_{\bar\xi}^{-1}C_{\bar\xi} \stackrel{\pi_{\bar\eta}}{\lra} \Flag_{\bar\eta}
\end{equation*}
is bijective. 
If $C_{\bar\xi}$ is compact, then this map is a homeomorphism,
and its inverse is a {\em section} of $\pi_{\bar\eta}$ whose image is {\em $\pi_{\bar\xi}$-saturated},
i.e.\ is a union of $\pi_{\bar\xi}$-fibers.
Conversely, each section of $\pi_{\bar\eta}$ whose image is $\pi_{\bar\xi}$-saturated yields 
a compact packing of $\geo X$ by $\pihalf$-balls of type $\bar\xi$. 

\begin{ques}
For which groups $G$ of type $B_2$ 
the projection $\pi_{\bar\xi}$ resp.\ $\pi_{\bar\eta}$
admits a section whose image is $\pi_{\bar\eta}$- resp.\ $\pi_{\bar\xi}$-saturated?
\end{ques}

\begin{example}
\label{ex:b2on2}
Let $G=SO(n,2)$ with $n\geq2$. 
The partial flag manifolds in this case are the 
Grassmannian ${\mathcal L}$ of isotropic lines $L$
and the Grassmannian ${\mathcal P}$ of isotropic planes $P$,
and the full flag manifold is the manifold ${\mathcal F}$ of isotropic flags $(L,P)$.

Fix an orthogonal splitting $\R^{n,2}=\R^n\oplus \R^2$ 
so that the quadratic form $q=x_1^2+\cdots +x_n ^2-x_{n+1}^2-x_{n+2}^2$ is definite on each factor. 
Then the isotropic planes in $\R^{n,2}$ are the graphs of isometries $\Phi: (\R^2, -q|_{\R^2})\to (\R^n, q|_{\R^n})$.
The isotropic lines in $\R^{n,2}$ are the graphs of isometries $\phi: (l, -q|_l)\to (\R^n, q|_{\R^n})$ 
defined on lines $l\subset\R^2$.
A full isotropic flag corresponds to a pair $(\Phi,l)$,
its isotropic line corresponding to the restriction $\Phi|_l$.
Thus, we have the product splitting ${\mathcal F}\cong{\mathcal P}\times \R P^1$.

The projection $\pi_{\mathcal P}:{\mathcal F}\to{\mathcal P}$ is the projection to the first factor.
It admits ``constant'' sections $s_l$ by fixing $l$.
Their images are $\pi_{\mathcal L}$-saturated,
namely $s_l({\mathcal P})=\pi_{\mathcal L}^{-1}{\mathcal L}_l$
where 
${\mathcal L}_l$
denotes the set of isotropic lines contained in the hyperplane
$\R^n\oplus l\cong\R^{n,1}\subset\R^{n,2}$.
The subset ${\mathcal L}_l\subset{\mathcal L}$ is compact and antipodal,
and hence constitutes the set of centers of a packing of $\geo X$ by $\pihalf$-balls.
In incidence geometric terms,
the antipodality corresponds to the fact that the hyperplane $\R^n\oplus l$ contains no isotropic plane,
and the packing to the (equivalent) fact that every isotropic plane intersects $\R^n\oplus l$ in a line.
The hyperplane $\R^n\oplus l$ is the orthogonal complement, in $\R^{n,2}$,
of the line $l^{\perp}\subset\R^2$ orthogonal to $l$.
Accordingly, ${\mathcal L}_l\subset{\mathcal L}$ is the orbit of a subgroup $\cong SO(n,1)$ of $SO(n,2)$,
namely of the one which fixes $l^{\perp}$.

The projection $\pi_{\mathcal L}:{\mathcal F}\to{\mathcal L}$ 
is given by $(\Phi,l)\mapsto\Phi|_l$.
Let ${\mathcal L}_1\subset{\mathcal L}$
denote the subset of isotropic lines 
which project to $l_1=\R e_1\subset\R^2$,
i.e.\ for which the isometry $\phi$ is defined on $l_1$.
A section of $\pi_{\mathcal L}$ over ${\mathcal L}_1$ 
would associate with each unit vector $\phi(e_1)$ in $\R^n$
a unit vector orthogonal to it, namely $\Phi(e_2)$ for the extension $\Phi$ of $\phi$ 
determined by the section.
It would thus yield a unit vector field on $S^{n-1}\subset\R^n$.
Such a vector field does not exist if $n$ is odd.

On the other hand, if $n$ is even, 
then we can use the standard identification $\R^{n,2}\cong\C^{\frac{n}{2},1}$
and consider the 
subset ${\mathcal P}_c\subset{\mathcal P}$
of isotropic planes which are invariant under the complex structure,
i.e.\ which are complex lines.
Every isotropic line is contained in a unique such isotropic plane by complexification,
which means that 
\begin{equation*}
\pi_{\mathcal P}^{-1}{\mathcal P}_c \stackrel{\pi_{\mathcal L}}{\lra} {\mathcal L}
\end{equation*}
is a homeomorphism 
and, accordingly,
${\mathcal P}_c$ is the set of centers of a packing.
It is the orbit of the subgroup $SU(\frac{n}{2},1)\subset SO(n,2)$.
\end{example}

Our discussion and the example imply:
\begin{thm}
\label{thm:nopackb2orth}
(i) If $G=O(2k+1, 2)$ with $k\geq1$, 
then $\geo X$ admits no compact packing by $\pihalf$-balls of regular type,
neither of the singular type corresponding to isotropic planes. 

(ii) If $G=O(2k, 2)$ with $k\geq1$, 
then $\geo X$ admits a compact packing by $\pihalf$-balls of the singular type corresponding to isotropic planes
whose set of centers is an orbit of $U(k,1)<G$.

(iii) If $G=O(n, 2)$ with $n\geq2$, 
then $\geo X$ admits a compact packing by $\pihalf$-balls of the singular type corresponding to isotropic lines
whose set of centers is an orbit of $O(n,1)<G$.
\end{thm}
\proof
(i) According to our above discussion, 
a packing of regular type gives rise to packings of both singular types.
We assume therefore that there exists a 
packing of the singular type corresponding to isotropic planes, 
i.e.\ with centers in ${\mathcal P}\subset\geo X$.
It yields a section of the fiber bundle $\pi_{\mathcal L}:{\mathcal F}\to{\mathcal L}$.
However, such a section does not exist, cf.\ the example, contradiction.

(ii)+(iii) See the previous example.
\qed

\medskip
The theorem leaves open the question whether packings of regular type exist if $n$ is even.

\subsection{Nonemptiness of domains of proper discontinuity}
\label{sec:ccddne}

We now apply our results on packings to discrete subgroups.
Proposition~\ref{prop:emptdpack} yields:
\begin{prop}[Empty domains yield packings]
\label{prop:neccdd}
Let $\taumod$ and $\bar\vartheta_0\in\taumod$ be $\iota$-invariant.
Suppose that 
$\Ga<G$ is a discrete subgroup
such that $\Lat$ is antipodal 
and $$\ThF(\Lat)=\DF X$$ for all balanced $\Wt$-left invariant thickenings $\Th\subset W$
of the form $\Th=\Th_{\bar\vartheta_0,\bar\vartheta,\pihalf}$
(as defined by (\ref{eq:metth})).
Then the family of balls 
${\mathcal B}(\Lat,\bar\vartheta_0)$ is a packing of $\geo X$.
\end{prop}
Applying our nonexistence results for packings (Theorem~\ref{thm:nopack}),
we conclude that some of the domains of proper discontinuity constructed earlier (cf.\ Theorem~\ref{thm:pdwconv})
are nonempty.
For instance, we obtain in the regular case $\taumod=\simod$:
\begin{thm}[Nonemptiness of domains of proper discontinuity]
\label{thm:nedoman}
Suppose that $X$ has at least one de Rham factor not of the type $A_1, B_2$ or $G_2$,
and let $\Ga<G$ be a $\simod$-convergence subgroup
with antipodal limit set $\Las$.
Then for some balanced thickening $\Th\subset W$
the domain of proper discontinuity $\DF X-\ThF(\Las)$ for the $\Ga$-action 
(provided by Theorem~\ref{thm:pdwconv})
is nonempty.
Moreover, the thickening can be chosen of the form $\Th=\Th_{\bar\vartheta_0,\bar\vartheta,\pihalf}$
(as defined by (\ref{eq:metth}))
with $\bar\vartheta_0\in\simod$ a root type.
\end{thm}
\proof
Otherwise,
by the proposition,
$\geo X$ admits compact packings by $\pihalf$-balls 
for all (of the at most two) root types.
However,
this contradicts 
Theorem~\ref{thm:nopacka2}, respectively, Theorem~\ref{thm:nopackr3}.
(Note that root types are $\iota$-invariant.)
\qed

\medskip
In the $B_2$ case, we can only treat a family of examples:
\begin{add}[Nonemptiness of domains of proper discontinuity, $B_2$ case]
\label{add:nedomanb2orth}
Suppose that $G=O(2k+1, 2)$ with $k\geq1$, 
and let $\Ga<G$ be a $\pimod$-convergence subgroup
with antipodal limit set $\Lap$ 
for the vertex type $\pimod\in\simod$ corresponding to isotropic planes.
Then for the balanced $\Wp$-left invariant thickening $\Th\subset W=W_{B_2}$
the domain of proper discontinuity $\DF X-\ThF(\Lap)$ for the $\Ga$-action 
(provided by Theorem~\ref{thm:pdwconv})
is nonempty.
\end{add}
\proof
Otherwise,
the proposition yields a compact packing of $\geo X$ by $\pihalf$-balls of type $\pimod$,
contradicting Theorem~\ref{thm:nopackb2orth}(i).
\qed

\medskip
This leaves open the question whether, 
in the case of $G=O(2k,2)$ for $k\geq2$, 
there are $\simod$-convergence subgroups with antipodal limit sets, which have 
empty domains of proper discontinuity 
for {\em arbitrary} balanced thickenings $\Th$. 
We note that Example \ref{ex:product_actions} provides examples of 
$\simod$-CEA subgroups with empty domains of proper discontinuity in $\DF X$  
for {\em some} choices of balanced thickenings.

\begin{rem}
Theorem~\ref{thm:nedoman} is both weaker and stronger than the 
nonemptiness results in \cite[Thms.\ 1.11, 1.12 and 9.10]{GW}. 
It is stronger in the sense that it applies to hyperbolic groups $\Ga$ without assumptions on their cohomological dimension,
unlike the results in \cite{GW} which require small cohomological dimension; 
furthermore, it applies to domains of discontinuity in various partial flag manifolds (always including $G/B=\DF X$), unlike the results in \cite{GW} which work (in general) only for domains of discontinuity in $G/AN$ 
(which is a certain fiber bundle over $G/B$). 
On the other hand, it is weaker in the sense that it addresses only the regular case ($\taumod=\simod$). 
We also note that some examples of Anosov subgroups for which some discontinuity domains are empty are given in \cite[Remark 8.5]{GW}.
\end{rem}

\bigskip 
Addresses:

\noindent M.K.: Department of Mathematics, \\
University of California, Davis\\
CA 95616, USA\\
email: kapovich@math.ucdavis.edu

\noindent B.L.: Mathematisches Institut\\
Universit\"at M\"unchen \\
Theresienstr. 39\\ 
D-80333, M\"unchen, Germany\\ 
email: b.l@lmu.de

\noindent J.P.: Departament de Matem\`atiques,\\
Universitat Aut\`onoma de Barcelona,\\ 
E-08193 Bellaterra, Spain\\
email: porti@mat.uab.cat

\end{document}